\documentclass[a4paper,10pt]{article}
\usepackage[utf8]{inputenc}
\usepackage[T1]{fontenc}
\usepackage{amsthm}
\usepackage{amsfonts}
\usepackage{amssymb}
\usepackage{mathrsfs}	
\usepackage{amsmath}
\allowdisplaybreaks
\usepackage{mathtools}
\usepackage[english]{babel}
\usepackage[table]{xcolor}
\usepackage{blkarray}
\usepackage{slashed}
\usepackage{enumerate}
\usepackage{graphicx}
\usepackage{systeme}
\usepackage{dsfont}
\usepackage{relsize}
\usepackage[margin=0.90in]{geometry}
\usepackage{float}
\usepackage{tikz}
\usepackage[labelformat=simple]{subcaption}

\usepackage{graphicx}
\usepackage{epstopdf}
\usepackage{bbm}
\usepackage{etoolbox}
\usepackage{titling}
\usepackage{authblk}

\usepackage{blindtext,graphicx}
\usepackage[absolute]{textpos}
\usepackage[hang,flushmargin]{footmisc}

\usepackage{xfrac}
\usepackage{nicefrac}
\usepackage{esvect}
\usepackage{cases}
\usepackage{empheq}
\usepackage{stmaryrd}
\usepackage{cancel}
\usepackage[linguistics]{forest}
\usepackage{url}
\usepackage[font=small]{caption}
\usepackage{hhline}
\RequirePackage[table]{xcolor}
\definecolor{tssteelblue}{RGB}{70,130,180}
\definecolor{tsorange}{RGB}{255,138,88}
\definecolor{tsblue}{RGB}{23,74,117}
\definecolor{tsforestgreen}{RGB}{21,122,81}
\definecolor{tsyellow}{RGB}{255,185,88}
\definecolor{tsgrey}{RGB}{200,200,200}
\definecolor{ttsorange}{RGB}{255,120,0}
\definecolor{forestgreen}{RGB}{0.13, 0.55, 0.13}
\definecolor{cof}{RGB}{219,144,71}
\definecolor{pur}{RGB}{186,146,162}
\definecolor{greeo}{RGB}{91,173,69}
\definecolor{greet}{RGB}{52,111,72}
\definecolor{myblue}{RGB}{80,191,227}
\definecolor{myblue2}{RGB}{80,141,227}
\newtheoremstyle{mystyle}
{}
{}
{}
{}
{\bfseries}
{.}
{ }
{}
\theoremstyle{plain}
\theoremstyle{mystyle}
\newtheorem{rem}{Remark}
\newenvironment{rema}
{\pushQED{\qed}\rem}
{\popQED\endrem}
\theoremstyle{plain}
\theoremstyle{mystyle}
\newtheorem{thm}{Theorem}
\theoremstyle{plain}
\theoremstyle{mystyle}
\newtheorem{defn}{Definition}
\usetikzlibrary{shapes.geometric,shapes.multipart}
\usetikzlibrary{decorations.pathmorphing,decorations.text,decorations.pathreplacing,decorations.markings}
\usetikzlibrary{calc,intersections,through,backgrounds}
\usetikzlibrary{patterns,snakes}
\usetikzlibrary{arrows.meta}
\usepackage{pgfplots}
\usetikzlibrary{tikzmark}
\usetikzlibrary{matrix}
\pgfdeclarelayer{edgelayer,nodelayer}
\usetikzlibrary{patterns,shadows,calc,arrows}
\tikzstyle{none}=[inner sep=0pt]
\RequirePackage{pgfplots}\pgfplotsset{compat=1.13}
\usetikzlibrary{backgrounds,scopes,cd,positioning,shadows,shapes.multipart}
\RequirePackage[tikz]{bclogo} 
\usetikzlibrary{cd}
\tikzcdset{
	arrow style=tikz,
	diagrams={>={Straight Barb[scale=0.8]}},math to/.tip={Glyph[glyph math command=rightarrow]},
	loop/.tip={Glyph[glyph math command=looparrowleft, swap]},
	loop'/.tip={Glyph[glyph math command=looparrowleft]},
	weird/.tip={Glyph[glyph math command=Rrightarrow, glyph length=1.5ex]},
	pi/.tip={Glyph[glyph math command=pi, glyph length=1.5ex, glyph axis=0pt]},
}
\RequirePackage{empheq}
\allowdisplaybreaks

\tikzset{
	math to/.tip={Glyph[glyph math command=rightarrow]},
	loop/.tip={Glyph[glyph math command=looparrowleft, swap]},
	loop'/.tip={Glyph[glyph math command=looparrowleft]},
	weird/.tip={Glyph[glyph math command=Rrightarrow, glyph length=1.5ex]},
	pi/.tip={Glyph[glyph math command=pi, glyph length=1.5ex, glyph axis=0pt]},
}
\newlength{\hatchspread}
\newlength{\hatchthickness}
\newlength{\hatchshift}
\newcommand{\hatchcolor}{}
\tikzset{hatchspread/.code={\setlength{\hatchspread}{#1}},
	hatchthickness/.code={\setlength{\hatchthickness}{#1}},
	hatchshift/.code={\setlength{\hatchshift}{#1}},
	hatchcolor/.code={\renewcommand{\hatchcolor}{#1}}}
\tikzset{hatchspread=3pt,
	hatchthickness=0.4pt,
	hatchshift=0pt,
	hatchcolor=black}
\tikzset{
	labl/.style={anchor=south, rotate=270, inner sep=.5mm}}
\pgfdeclarepatternformonly[\hatchspread,\hatchthickness,\hatchshift,\hatchcolor]
{custom north west lines}
{\pgfqpoint{\dimexpr-2\hatchthickness}{\dimexpr-2\hatchthickness}}
{\pgfqpoint{\dimexpr\hatchspread+2\hatchthickness}{\dimexpr\hatchspread+2\hatchthickness}}
{\pgfqpoint{\dimexpr\hatchspread}{\dimexpr\hatchspread}}
{
	\pgfsetlinewidth{\hatchthickness}
	\pgfpathmoveto{\pgfqpoint{0pt}{\dimexpr\hatchspread+\hatchshift}}
	\pgfpathlineto{\pgfqpoint{\dimexpr\hatchspread+0.15pt+\hatchshift}{-0.15pt}}
	\ifdim \hatchshift > 0pt
	\pgfpathmoveto{\pgfqpoint{0pt}{\hatchshift}}
	\pgfpathlineto{\pgfqpoint{\dimexpr0.15pt+\hatchshift}{-0.15pt}}
	\fi
	\pgfsetstrokecolor{\hatchcolor}
	\pgfusepath{stroke}
}
\RequirePackage{array, boldline, rotating}
\RequirePackage{colortbl,bigstrut}
\usepackage{tabularx}
\usepackage{array}
\RequirePackage{tabu}
\newcolumntype{M}[1]{>{\centering\arraybackslash}m{#1}}
\newcolumntype{N}{@{}m{0pt}@{}}
\usepackage{multirow}
\usepackage[colorinlistoftodos,prependcaption,textsize=tiny]{todonotes}
\newcommand{\skalarProd}[2]{\big\langle#1,#2\big\rangle}

\newcommand{\norm}[1]{\lVert #1 \rVert}

\global\long\def\sA{\mathscr{A}}
\global\long\def\sB{\mathscr{B}}

\global\long\def\sD{\mathscr{D}}	

\global\long\def\sF{\mathscr{F}}

\global\long\def\sH{\mathscr{H}}

\global\long\def\bC{\mathbb{C}}

\global\long\def\bG{\mathbb{G}}
\global\long\def\bH{\mathbb{H}}

\global\long\def\bL{\mathbb{L}}

\global\long\def\bN{\mathbb{N}}

\global\long\def\bR{\mathbb{R}}

\global\long\def\boA{\boldsymbol{A}}
\global\long\def\boB{\boldsymbol{B}}

\global\long\def\boD{\boldsymbol{D}}

\global\long\def\boP{\boldsymbol{P}}
\global\long\def\boQ{\boldsymbol{Q}}

\global\long\def\boU{\boldsymbol{U}}

\global\long\def\bob{\boldsymbol{b}}
\global\long\def\boe{\boldsymbol{e}}
\global\long\def\bog{\boldsymbol{g}}
\global\long\def\bom{\boldsymbol{m}}
\global\long\def\bou{\boldsymbol{u}}
\global\long\def\bov{\boldsymbol{v}}

\global\long\def\boox{\boldsymbol{x}}
\global\long\def\bovarphi{\boldsymbol{\varphi}}

\global\long\def\gR{\mathfrak{R}}

\global\long\def\go{\mathfrak{o}}

\global\long\def\gs{\mathfrak{s}}

\global\long\def\fr{\rightarrow}

\global\long\def\id{\mathbbm{1}}	
\global\long\def\defi{\vcentcolon =}

\global\long\def\Sym{\textrm{Sym}\,}
\global\long\def\sym{\textrm{sym}\,}
\global\long\def\skew{\textrm{skew}\,}

\global\long\def\dev{\textrm{dev}\,}

\global\long\def\Curl{\textrm{Curl}\,}
\global\long\def\curl{\textrm{curl}\,}
\global\long\def\div{\textrm{div}\,}
\global\long\def\Div{\textrm{Div}\,}
\global\long\def\DIV{\textrm{DIV}\,}
\global\long\def\tr{\textrm{tr}\,}
\global\long\def\supp{\textrm{supp}\,}
\global\long\def\SO{\textrm{SO}(3)}
\global\long\def\so{\mathfrak{so}(3)}
\global\long\def\sso(n){\mathfrak{so}\,(n)}
\global\long\def\sl(n){\mathfrak{sl}\,(n)}

\global\long\def\muc{\mu_{\textrm{c}}}
\global\long\def\mum{\mu_{\textrm{micro}}}
\global\long\def\muM{\mu_{\textrm{macro}}}
\global\long\def\lam{\lambda_{\textrm{micro}}}

\global\long\def\kam{\kappa_{\textrm{micro}}}
\global\long\def\me{\mu_{\text{e}}}
\global\long\def\mue{\mu_{\text{e}}}
\global\long\def\le{\lambda_{\text{e}}}
\global\long\def\pu{\boldsymbol{u}}
\global\long\def\pv{\boldsymbol{v}}
\global\long\def\pP{\boldsymbol{P}}
\global\long\def\pQ{\boldsymbol{Q}}
\global\long\def\pnu{\boldsymbol{\nu}}
\global\long\def\gbm{\boldsymbol{\mathfrak{m}}}
\global\long\def\gbu{\boldsymbol{\mathfrak{u}}}

\global\long\def\pD{\text{D}}
\newcommand{\D}{\operatorname{D}}
\global\long\def\bosigma{\boldsymbol{\sigma}}
\DeclareMathOperator{\Anti}{\textrm{Anti}}
\DeclareMathOperator{\axl}{\textrm{axl}}

\global\long\def\LC{\boldsymbol{\epsilon}}
\usepackage[giveninits=true,sortcites=true,date=year,maxbibnames=99,doi=false,isbn=false,url=false,eprint=false]{biblatex}
\renewbibmacro{in:}{}
\DeclareFieldFormat{pages}{#1}
\AtEveryBibitem{%
  \clearlist{language}%
}
\usepackage{csquotes}

\usepackage{setspace}
\title{\vspace{-3cm}The consistent coupling boundary condition for the classical micromorphic model: existence, uniqueness and interpretation of parameters}
\date{}

\begin{document}
	\author{Marco Valerio d\textquoteright Agostino\thanks{Marco Valerio d'Agostino, corresponding author, marco-valerio.dagostino@insa-lyon.fr,
			GEOMAS, INSA-Lyon, Université de Lyon, 20 avenue Albert Einstein,
			69621, Villeurbanne cedex, France},
		$\;\;$Gianluca Rizzi\thanks{Gianluca Rizzi, gianluca.rizzi@insa-lyon.fr, Technische Universität Dortmund, August-Schmidt-Str. 8, 44227 Dortmund, Germany},
		$\;\;$Hassam Khan\thanks{Hassam Khan,  hassam.khan@stud.uni-due.de,  Fakultät für Mathematik, Universität Duisburg-Essen,
			Mathematik-Carrée, Thea-Leymann-Straße 9, 45127 Essen, Germany},
		$\;\;$Peter Lewintan\thanks{Peter Lewintan, peter.lewintan@uni-due.de, Chair for Nonlinear
			Analysis and Modelling, Fakultät für Mathematik, Universität Duisburg-Essen,
			Mathematik-Carrée, Thea-Leymann-Straße 9, 45127 Essen, Germany},
		\\
		$\;\;$Angela Madeo\thanks{Angela Madeo, angela.madeo@insa-lyon.fr, Technische Universität Dortmund, August-Schmidt-Str. 8, 44227 Dortmund, Germany}
		$\;$and$\;$  Patrizio Neff\,\thanks{Patrizio Neff, patrizio.neff@uni-due.de, Head of Chair for Nonlinear
			Analysis and Modelling, Fakultät für Mathematik, Universität Duisburg-Essen,
			Mathematik-Carrée, Thea-Leymann-Straße 9, 45127 Essen, Germany}}
	\thanksmarkseries{arabic}
	\maketitle
	\vspace{-17mm}
	\begin{abstract}
		We consider the classical Mindlin-Eringen linear micromorphic model with
		a new strictly weaker set of displacement boundary conditions. The new consistent
		coupling condition aims at minimizing spurious influences from arbitrary
		boundary prescription for the additional microdistortion field $\pP$. In
		effect, $\pP$ is now only required to match
		the tangential derivative of the classical displacement $\pu$ which
		is known at the Dirichlet-part of the boundary.
		
		We derive the full boundary condition, in adding the missing Neumann condition
		on the Dirichlet-part. We show existence and uniqueness of the static
		problem for this weaker boundary condition. These results are based
		on new coercive inequalities for incompatible tensor fields with prescribed
		tangential part.
		
		Finally, we show that compared to classical Dirichlet conditions on $\bou$ and $\boP$,
		the new boundary condition modifies the interpretation of the constitutive
		parameters.
	\end{abstract}
	\addtocounter{footnote}{5}
	\tableofcontents{}
	\section{Introduction}
	Classical linear elasticity is sufficient to describe the material response of structures in which no microstructural effects are taken into account. Recent advances in science and technology necessitate, however, to describe materials in which the effects of the microstructure cannot be omitted.
	
	The presence of a microstructure show itself in a size-effect (smaller is stiffer) for the static case and in non-classical dispersion curves, extra vibration modes, etc. in the dynamic case.
	One of the best known models incorporating a rigid substructure is the Cosserat model \cite{cosserat1909theorie,jeong2008existence,jeong2009numerical,neff2010stable,neff2009subgrid,neff2006cosserat,taliercio2009some,Taliercio2010,lakes1995experimental,Lakes2015,rueger2018strong,masiani1996cosserat,fantuzzi2018some,fantuzzi2020material}. Here, the kinematics includes the standard translational degrees of freedom and the rigid rotation of the substructure. A generalization can be obtained by successively enriching the kinematics further. For example, the micro-stretch model \cite{neff2009mean,neff2014unifying,scalia2000extension,de1997torsion,kirchner2007mechanics} allows for extra rotation and volume increase of the substructure. Finally, the  micromorphic approach \cite{eringen1968mechanics,mindlin1964micro} allows for a substructure that has 9 dofs to describe a full affine motion (see also \cite{d2020effective,dagostino2016panorama,barbagallo2016transparent,barbagallo2018relaxed,aivaliotis2019microstructure,madeo2016reflection,madeo2015wave,madeo2014band,ghiba2014relaxed,owczarek2019nonstandard,neff2014unifying,neff2016real,Forest2018,Forest2013_CISM,shaat2020review,hutter2017homogenization,hutter2019micro,mariano2009ground,mariano2005computational,trovalusci2017multiscale,trovalusci2015scale,trovalusci1998continuum,wallen2021willis,mahnken2020goal,scherer2020lagrange,ALBERDI2021104540,GANGHOFFER1998125} ). 
	
	The common problem in these so-called generalized continua is the increased number of constitutive coefficients 
	that are not easy to be determined together with the possible choice of boundary conditions for the enriched kinematics. In this paper we want to mainly concentrate on the latter question of boundary conditions of the linear isotropic micromorphic model. 
	
	The model is always formulated as a minimization problem in some appropriate function spaces. Our goal is to propose an alternative boundary condition, called consistent coupling, that will improve some of the problematic issues in generalized continua. An understanding of this issue is fundamental for any progress in parameter identification in order to be relevant for applications.
	
	This paper is structured as follows. First we discuss the linear isotropic micromorphic model as a minimization problem together with the traditional Dirichlet boundary condition for the enriched kinematics. Then we introduce our new weaker consistent coupling boundary condition and show that it still leads to an existence and uniqueness result. We provide some mechanical analogies to motivate the new condition and derive some simple analytical solutions for the traditional and new boundary condition which allows us to better understand their differences. Then we are naturally lead to investigate rigid body modes as compared to zero energy modes for generalized continua and how boundary conditions can be used to set these modes to zero. We do this for a whole range of generalized continua. Here, we are touching a lesser known concept of redundancy in generalized continua, introduced by Romano et al. \cite{romano2016micromorphic}. We are finally entering a discussion of additional invariance conditions beyond frame-indifference and isotropy, loosely related to the foregoing discussion.
	In the appendix we provide the detailed derivation for all combinations of boundary conditions.
	\section{Notations and technical preliminaries}
	Let $n\geq2$. For vectors $\boldsymbol{a},\boldsymbol{b}\in\bR^n$ we consider the scalar product  $\skalarProd{\boldsymbol{a}}{\boldsymbol{b}}\coloneqq\sum_{i=1}^n a_i\,b_i \in \bR$, the (squared) norm  $\norm{\boldsymbol{a}}^2\coloneqq\skalarProd{\boldsymbol{a}}{\boldsymbol{a}}$ and the dyadic product  $\boldsymbol{a}\otimes \boldsymbol{b} \defi \left(a_i\,b_j\right)_{i,j=1,\ldots,n}\in \bR^{n\times n}$. Similarly, the scalar product for matrices $\pP,\boldsymbol{Q}\in\bR^{n\times n}$ is given by $\skalarProd{\pP}{\boldsymbol{Q}} \defi \sum_{i,j=1}^n P_{ij}\,Q_{ij} \in \bR$ and the (squared) Frobenius-norm by $\norm{\pP}^2\defi\skalarProd{\pP}{\pP}$.
	Moreover, $\pP^T\coloneqq (P_{ji})_{i,j=1,\ldots,n}$ stands for the transposition of the matrix $\pP=(\pP_{ij})_{i,j=1,\ldots,n}$. We make use of the orthogonal decomposition of the latter into the symmetric part $\sym \pP \coloneqq \frac12\left(\pP+\pP^T\right)$ and the skew-symmetric part $\skew \pP \coloneqq \frac12\left(\pP-\pP^T\right)$. We denote by $\Sym(n)\coloneqq\{\boldsymbol{X}\in\bR^{n\times n}\mid  \boldsymbol{X}^T=\boldsymbol{X}\}$ the vector space of symmetric matrices and by $\sso(n)\coloneqq \{\boldsymbol{A}\in\bR^{n\times n}\mid \boldsymbol{A}^T = -\boldsymbol{A}\}$ the Lie-Algebra of skew-symmetric matrices. For the identity matrix we write $\id$, so that the trace of a square matrix $\pP$ is $\tr \pP \coloneqq \skalarProd{\pP}{\id}$. The deviatoric (trace-free) part of $\pP$ is given by $\dev_{\!n} \,\pP\coloneqq \pP -\frac1n\tr(\pP)\,\id$ and in three dimensions we will suppress its index, i.e.,~we write $\dev$ instead of $\dev_{\!3}$. Third order tensor fields will be denoted by bold-frak letters.
	We denote by $\!\cdot\!$ the single contraction of tensors and by $:$ the double contraction. Specifically, for the cases we are concerned with, we have
	\[
	\begin{aligned}
	\forall\left(\pP,\pu\right) & \in\bR^{3\times3}\times\bR^{3}, & \pP\!\!\cdot\!\pu & \in\bR^{3}, &  & \textrm{and} & \left(\pP\!\cdot\!\pu\right)_{i} & =\textstyle{\sum_{j=1}^3\pP_{ij}\pu_{j}},
	\\
	\forall\left(\pP,\pQ\right) & \in\bR^{3\times3}\times\bR^{3\times3}, & \pP\!\!\cdot\!\pQ & \in\bR^{3\times3}, &  & \textrm{and} & \left(\pP\!\cdot\!\pQ\right)_{ij} & =\textstyle{\sum_{h=1}^3}\pP_{ih}\pQ_{hj},
	\\
	\forall\left(\gbm,\pu\right) & \in\bR^{3\times3\times3}\times\bR^{3}, & \gbm\!\cdot\!\pu & \in\bR^{3\times3}, &  & \textrm{and} & \left(\gbm\!\cdot\!\pu\right)_{ij} & =\textstyle{\sum_{h=1}^3}\gbm_{ijh}\pu_{h},
	\\
	\forall\left(\gbm,\pP\right) & \in\bR^{3\times3\times3}\times\bR^{3\times3}, & \gbm\!\!:\!\!\pP & \in\bR^{3}, &  & \textrm{and} & \left(\gbm\!:\!\pP\right)_{i} & =\textstyle{\sum_{j,h=1}^3}\gbm_{ijh}\pP_{jh}.
	\end{aligned}
	\]
	The double and triple contractions between fourth order tensors and second order tensors, and sixth order tensors and third order tensors will be denoted just by a juxtaposition, as well as the generic action of a linear function. 
	Let $\sD'(\Omega, \bR^3),\sD'(\Omega, \bR^{3\times 3})$ and $\sD'(\Omega, \bR^{3\times 3\times 3})$ respectively be the spaces of the vector, matrix and third order tensor valued distributions over a Lipschitz domain $\Omega$.  Let us introduce the following Sobolev spaces: 
	\begin{align}
	H\!\left(\div;\Omega,\bR^{3}\right) & \defi\left\{ \bou\in L^{2}\!\left(\Omega,\bR^{3}\right)\left.\right|\div\bou\in L^{2}\!\left(\Omega,\bR\right)\right\}, \nonumber\\
	H\!\left(\Div;\Omega,\bR^{3\times3}\right) & \defi\left\{ \bosigma\in L^{2}\!\left(\Omega,\bR^{3\times3}\right)\left.\right|\Div\bosigma\in L^{2}\!\left(\Omega,\bR^{3}\right)\right\}, \\
	H\!\left(\DIV;\Omega,\bR^{3\times3\times3}\right) & \defi\left\{ \gbm\in L^{2}\!\left(\Omega,\bR^{3\times3\times3}\right)\left.\right|\DIV\gbm\in L^{2}\!\left(\Omega,\bR^{3\times3}\right)\right\}, \nonumber
	\end{align}
	where 
	\begin{equation}
	\Div\bosigma=
	\begin{pmatrix}
	\div\!\left(\sigma_{11},\sigma_{12},\sigma_{13}\right)\\
	\div\!\left(\sigma_{21},\sigma_{22},\sigma_{23}\right)\\
	\div\!\left(\sigma_{31},\sigma_{32},\sigma_{33}\right)
	\end{pmatrix}
	\end{equation} 
	and
	\begin{equation}
	\DIV\gbm=
	\begin{pmatrix}
	\div\!\left(\mathfrak{m}_{111},\mathfrak{m}_{112},\mathfrak{m}_{113}\right) & \div\!\left(\mathfrak{m}_{121},\mathfrak{m}_{122},\mathfrak{m}_{123}\right) & \div\!\left(\mathfrak{m}_{131},\mathfrak{m}_{132},\mathfrak{m}_{133}\right)
	\\
	\div\!\left(\mathfrak{m}_{211},\mathfrak{m}_{212},\mathfrak{m}_{213}\right) & \div\!\left(\mathfrak{m}_{221},\mathfrak{m}_{222},\mathfrak{m}_{223}\right) & \div\!\left(\mathfrak{m}_{231},\mathfrak{m}_{232},\mathfrak{m}_{233}\right)
	\\
	\div\!\left(\mathfrak{m}_{311},\mathfrak{m}_{312},\mathfrak{m}_{313}\right) & \div\!\left(\mathfrak{m}_{321},\mathfrak{m}_{322},\mathfrak{m}_{323}\right) & \div\!\left(\mathfrak{m}_{331},\mathfrak{m}_{332},\mathfrak{m}_{333}\right)
	\end{pmatrix}.
	\end{equation}
	Setting $\left(\pP\right)_{i}=(\pP_{ij})_{j}$ we introduce the cross-product operation between a matrix valued tensor field $\pP$ and a vector field $\bou$ always row-wise, i.e., $\left(\pP\times\bou\right)_{ij}=\left(\left(\pP\right)_{i}\times\bou\right)_{j}$ for every $i,j\in\left\{ 1,2,3\right\} $. The last space we need is the Sobolev space $H\!\left(\Curl;\Omega,\bR^{3\times3}\right)\defi\left\{ \pP\in L^{2}\!\left(\Omega,\bR^{3\times3}\right)\left.\right|\Curl\pP\in L^{2}\!\left(\Omega,\bR^{3\times3}\right)\right\} $ defined via the matrix-generalization $\Curl$ of the classical $\curl$ operator in the following way
	\begin{align}
	\left\langle \Curl\pP\left.\right|\delta\pP\right\rangle _{\sD'(\Omega,\bR^{3\times3}),\sD(\Omega,\bR^{3\times3})} & =\sum_{i=1}^{3}\left\langle \curl\!\left(\pP\right)_{i}\left.\right|\left(\delta\pP\right)_{i}\right\rangle _{\sD'(\Omega,\bR^{3}),\sD(\Omega,\bR^{3})}\nonumber=\sum_{i=1}^{3}\left\langle \left(\pP\right)_{i}\left.\right|\curl\!\left(\delta\pP\right)_{i}\right\rangle _{\sD'(\Omega,\bR^{3}),\sD(\Omega,\bR^{3})}
	\end{align}
	for all $\left(\pP,\delta\pP\right)\in\sD'\left(\Omega,\bR^{3\times3}\right)\times\sD\left(\Omega,\bR^{3\times 3}\right)$, giving,
	\begin{equation}\label{Curl definition}
	\Curl\pP=
	\begin{pmatrix}
	\curl\!\left(P_{11},P_{12},P_{13}\right)\\
	\curl\!\left(P_{21},P_{22},P_{23}\right)\\
	\curl\!\left(P_{31},P_{32},P_{33}\right)
	\end{pmatrix} \, .
	\end{equation}
	For every $u\in H^1(\Omega,\bR^{3})$ we have the bound
	\begin{align}
	&\Vert \curl u \Vert _{L^{2}\left(\Omega\right)}^{2}  =\left\Vert \left(u_{3,2}-u_{2,3},u_{1,3}-u_{3,1},u_{2,1}-u_{1,2}\right)\right\Vert _{L^{2}\left(\Omega\right)}^{2}=\left\Vert \left(u_{3,2},u_{1,3},u_{2,1}\right)-\left(u_{2,3},u_{3,1},u_{1,2}\right)\right\Vert _{L^{2}\left(\Omega\right)}^{2}
	\\
	& \leqslant\left(\left\Vert \left(u_{3,2},u_{1,3},u_{2,1}\right)\right\Vert _{L^{2}\left(\Omega\right)}+\left\Vert \left(u_{2,3},u_{3,1},u_{1,2}\right)\right\Vert _{L^{2}\left(\Omega\right)}\right)^{2}\vphantom{\sum_{i,j=1}^{3}}
	\!\!\!\leqslant2\left(\left\Vert \left(u_{3,2},u_{1,3},u_{2,1}\right)\right\Vert _{L^{2}\left(\Omega\right)}^{2}+\left\Vert \left(u_{2,3},u_{3,1},u_{1,2}\right)\right\Vert _{L^{2}\left(\Omega\right)}^{2}\right) \nonumber
	\\
	& 
	=2\sum_{\substack{i,j=1\\
			i\neq j
		}
	}^{3}\left\Vert u_{i,j}\right\Vert _{L^{2}\left(\Omega\right)}^{2}\leqslant2\sum_{i,j=1}^{3}\left\Vert u_{i,j}\right\Vert _{L^{2}\left(\Omega\right)}^{2}=2\left\Vert \text{D} u\right\Vert _{L^{2}\left(\Omega\right)}^{2}. \nonumber
	\end{align}
	giving the trivial tensor-valued analogous
	\begin{align}\label{Curl Bound}
	\left\Vert \Curl P\right\Vert _{L^{2}\left(\Omega\right)}^{2} & =\sum_{i=1}^{3}\left\Vert \curl(P)_{i}\right\Vert _{L^{2}\left(\Omega\right)}^{2}\leqslant2\sum_{i=1}^{3}\left\Vert \text{D}(P)_{i}\right\Vert _{L^{2}\left(\Omega\right)}^{2}=2\left\Vert \text{D} P\right\Vert _{L^{2}\left(\Omega\right)}^{2}.
	\end{align}
	Here,	$\Omega$ will denote a bounded Lipschitz domain in $\bR^3$ with boundary $\partial \Omega$ whose unit normal (external) vector field (defined almost everywhere) is $\pnu$. We always will account for $\Gamma\subseteq\partial\Omega$ relatively open in $\partial\Omega$ and such that $\partial\Omega\setminus\overline{\Gamma}$ will also be relatively open in $\partial\Omega$.
	Given a space of functions $\Xi$ (in general an open subset of a Banach space), the Gâteaux derivative of a functional $\sF:\Xi\fr\bR$ at $u$ in the direction $\delta u$ is 
	$$
	\delta\sF\left[u,\delta u\right]\defi\left.\frac{d}{dt}\right|_{t=0}\sF\left[u+t\,\delta u\right].
	$$
	\section{The classical Mindlin-Eringen micromorphic model}
\begin{figure}[H]
		\begin{centering}
			\begin{tabular}{  cc }
				\begin{tabular}{c}\includegraphics[scale=1.3]{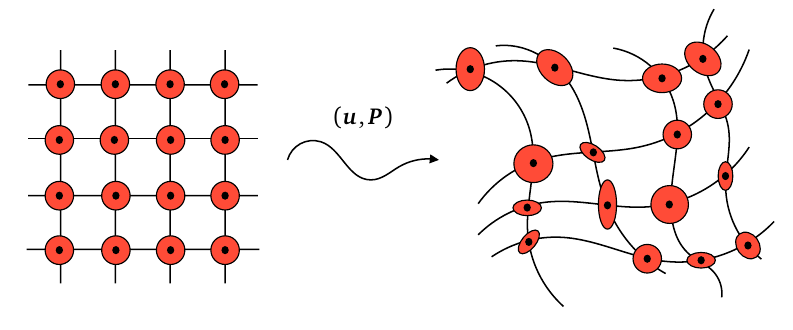} \end{tabular} 
				&
				\begin{tabular}{c}
					\includegraphics[scale=1.3]{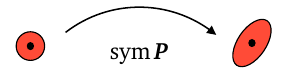} \end{tabular} \tabularnewline
				(a) & (b)\tabularnewline
			\end{tabular}
			\par\end{centering}
		\caption{
			(a) The kinematical assumptions underlying the micromorphic approach; the macroscopic displacement 
			$\pu-$distorts the lattice, while the microdistortion $\pP-$represents the affine deformation of the substructure attached in each macroscopic (lattice) point. (b) $\sym\pP$ \underline{controls} the size and shape of the affine mapping (affine microstructure), but \underline{not} the rotations of the ellipses. In principle, the microdistortion $\pP$ is independent of the macroscopic displacement gradient $\text{D}\pu$. Strain and curvature measures are needed to fully specify the response of the micromorphic model.}
			\label{fig:figure1}
	\end{figure}

	The classical linear Mindlin-Eringen micromorphic model \cite{mindlin1964micro,eringen1966mechanics,eringen1968mechanics} without
	body loads can be written as a minimization problem
	\begin{equation}
		I(\pu,\pP)=\int_{\Omega}W\left(\text{D}\boldsymbol{u},\boldsymbol{P},\text{D}\boldsymbol{P}\right)\text{dV}\qquad\longrightarrow\;\min\,\left(\pu,\pP\right)
	\end{equation}
	for the macroscopic displacement field $\pu:\Omega\subseteq\bR^{3}\fr\bR^{3}$
	and the non-symmetric microdistortion \linebreak ${\pP:\Omega\subseteq\bR^{3}\fr\bR^{3\times3}}$ (see Fig. \ref{fig:figure1}).
	Boundary conditions of Dirichlet-type
	\begin{equation}
		\left.\pu\right|_{\Gamma}=\widehat{\pu},\qquad\qquad\left.\pP\right|_{\Gamma}=\widehat{\pP}\qquad\textrm{at}\qquad\Gamma\subseteq\partial\Omega\label{eq:classical bc micromorphic}
	\end{equation}
	can typically be prescribed and lead to standard existence and uniqueness
	results in the space $\left(\pu,\pP\right)\in H^{1}\!\left(\Omega,\bR^{3}\right)\times H^{1}\!\left(\Omega,\bR^{3\times3}\right)$.
	Here, $W\!\left(\text{D}\boldsymbol{u},\boldsymbol{P},\text{D}\boldsymbol{P}\right)$
	denotes a non-negative quadratic form in the set of strain and curvature measures
	\begin{equation}
		\text{D}\boldsymbol{u}-\boldsymbol{P},\qquad\sym\pP,\qquad\text{D}\boldsymbol{P}.
	\end{equation}
	Depending on constitutive parameters, the energy $W$ satisfies typically
	estimates of the form

	\begin{subnumcases}{W\!\left(\text{D}\boldsymbol{u},\boldsymbol{P},\text{D}\boldsymbol{P}\right)\geqslant}
		c^{+}\left(\left\Vert \text{D}\boldsymbol{u}-\boldsymbol{P}\right\Vert ^{2}+\left\Vert \sym\boldsymbol{P}\right\Vert ^{2}+\left\Vert \text{D}\boldsymbol{P}\right\Vert ^{2}\right) \label{eq:energy1a}
		\\
		c^{+}\left(\left\Vert \sym\!\!\left(\text{D}\boldsymbol{u}-\boldsymbol{P}\right)\right\Vert ^{2}+\left\Vert \sym\boldsymbol{P}\right\Vert ^{2}+\left\Vert \text{D}\boldsymbol{P}\right\Vert ^{2}\right). \label{eq:energy1b}
	\end{subnumcases}
	In the case of availability of the traditional uniformly positive
	definite estimate (\ref{eq:energy1a}), the unique
	solution $(\pu,\pP)$ can be found in $H^{1}\!\left(\Omega,\bR^{3}\right)\times H^{1}\!\left(\Omega,\bR^{3\times3}\right)$
	provided merely $\left.\pu\right|_{\Gamma}=\widehat{\pu}$ is given.
	Indeed, since $\sym \pP\in L^{2}$ is automatically controlled, we may split $\text{D}\pu-\pP$ into symmetric and skew-symmetric parts and we  get $\text{D}\pu\in L^{2}$
	from the classical Korn's inequality \cite{neff2012maxwell,neff2002korn} and \eqref{eq:classical bc micromorphic}$_1$
	and this in turn controls $\pP\in L^{2}(\Omega)$ so that
	$\pu\in H^{1}\!\left(\Omega,\bR^{3}\right)$ and $\pP\in H^{1}\!\left(\Omega,\bR^{3\times3}\right)$.
	If, on the other side, only estimate (\ref{eq:energy1b})
	is available, then, $\pu\in H^{1}\!\left(\Omega,\bR^{3}\right)$ follows
	as before, but for $\pP$ we only control
	\begin{equation}
		\left\Vert \sym\pP\right\Vert ^{2}+\left\Vert \text{D}\pP\right\Vert ^{2}\label{eq:Stima su P}
	\end{equation}
	instead of $\left\Vert \pP\right\Vert ^{2}+\left\Vert \text{D}\pP\right\Vert ^{2}$ and without a Dirichlet-boundary condition for $\pP$ the microdistortion
	is not uniquely determined. If we add, as is classically done
	\begin{equation}
		\left.\pP\right|_{\Gamma}=\widehat{\pP}\label{eq:Condition su P}
	\end{equation}
	then (\ref{eq:Stima su P}) yields $\pP\in H^{1}\!\left(\Omega,\bR^{3\times3}\right)$,
	by invoking Poincaré's inequality for $\left\|\text{D}\pP \right\|^2 $ even without using the term $\left\|\sym\pP \right\|^2$. However, the interesting open question is whether the boundary
	condition (\ref{eq:Condition su P}) is really necessary in order to have
	$\pP\in H^{1}\!\left(\Omega,\bR^{3\times3}\right)$. The
	answer given in this paper is no. In fact, we will show subsequently, that the tangential boundary
	condition 
	\begin{equation}
		\left.\pP\times\pnu\,\right|_{\Gamma}= \widehat{\pP}\times\pnu\,|_{\Gamma}
	\end{equation}
	is already sufficient, where $\pnu$ denotes the outer unit normal on $\Gamma$. From a modelling perspective it is then suggestive
	to require the \underline{consistent coupling condition}
	\begin{equation}
		\left.\text{D}\boldsymbol{u}\times\pnu\,\right|_{\Gamma}=\left.\pP\times\pnu\,\right|_{\Gamma},\label{eq:new Bc}
	\end{equation}
	already known from, and naturally appearing in, the relaxed micromorphic model \cite{neff2014unifying,neff2016real,madeo2014band,madeo2016complete,madeo2016reflection,barbagallo2016transparent,aivaliotis2019microstructure}.
	Observe that the left-hand side of (\ref{eq:new Bc}) is known, once
	$\left.\pu\right|_{\Gamma}$ is given.
	
	In this paper we will fully investigate the new boundary condition (\ref{eq:new Bc}) for the classical Mindlin-Eringen micromorphic model.
	We will show existence and uniqueness under this weaker requirement, we will
	derive the missing normal Neumann-type condition on the Dirichlet boundary $\Gamma$ for
	the third order moment stress tensor $\gbm\simeq\text{D}\pP$ and we will compare
	analytical solutions for (\ref{eq:Condition su P}) versus (\ref{eq:new Bc}).
	This allows us to better understand the different Dirichlet conditions
	on $\pP$ (which both lead to existence and uniqueness).
	
	Furthermore, we consider the smooth transition of the micromorphic model
	to a second gradient continuum thanks to the  consistent coupling condition \eqref{eq:new Bc}. We
	will also briefly comment on the related notion of redundancy \cite{romano2016micromorphic}
	of elastic energies in generalized continua together with a discussion of invariance requirements beyond objectivity and isotropy.
    \subsection{The Cosserat model}\label{Cosserat model}
\begin{figure}[H]
		\begin{centering}
			\includegraphics[scale=0.8]{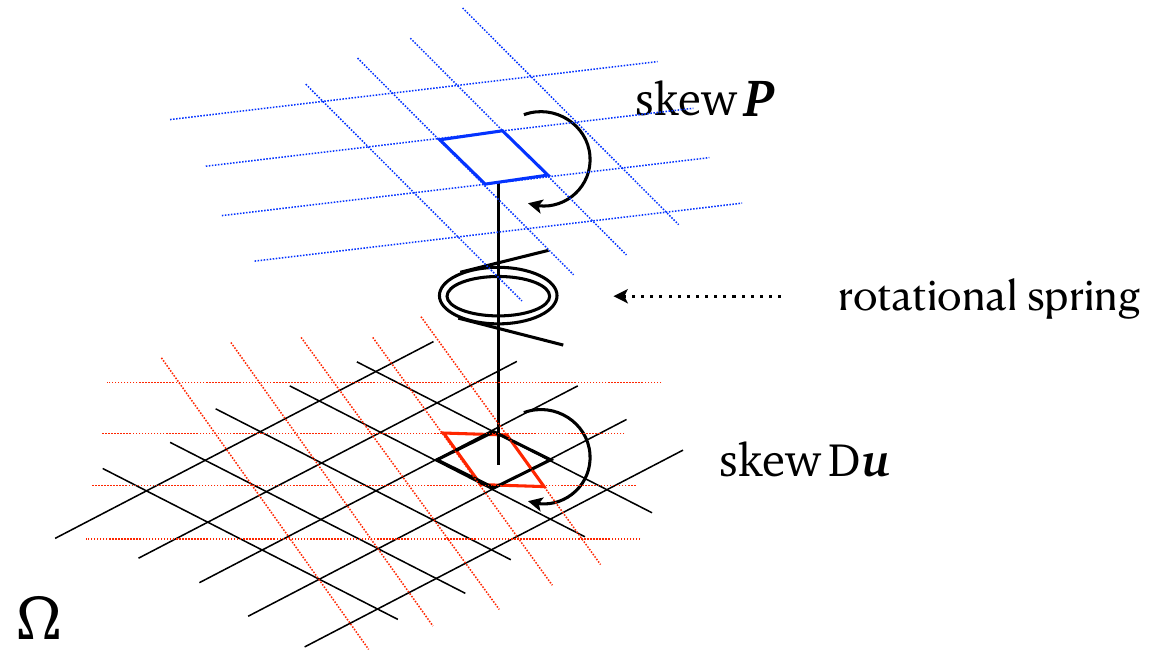}
			\par\end{centering}
		\caption{The Cosserat couple modulus $\muc$ stipulates an elastic spring between the infinitesimal continuum rotation $\skew\text{D}\pu$ and the infinitesimal microrotation $\skew\pP$. This rotational spring with spring constant $\muc$ is not needed to establish well-posedness of the model \cite{Lp_trace}. To the contrary, the Cosserat model is only operative with $\muc>0$.}
		\label{fig:mucRotSpring}
	\end{figure}

    In the isotropic linear Cosserat model the kinematics consists of the displacement $\pu$ and the infinitesimal microrotation $\boA\in\so$ (see Fig.~\ref{fig:mucRotSpring}). In a variational setting the problem reads 
    \begin{equation}
    I\left(\pu,\boA\right)=\int_{\Omega}W(\pD\pu,\boA,\Curl\boA)\,\text{dV}\longrightarrow\min\quad\left(\pu,\boA\right).\nonumber
    \end{equation}
    Boundary conditions of Dirichlet-type
    \begin{equation}
    \left.\pu\right|_{\Gamma}=\widehat{\pu},\qquad\qquad\left.\boA\right|_{\Gamma}=\widehat{\boA}\qquad\textrm{at}\qquad\Gamma\subseteq\partial\Omega
    \end{equation}
    can typically be prescribed and lead to standard existence and uniqueness
    results in the space $\left(\pu,\boA\right)\in H^{1}(\Omega,\bR^3)\times H^{1}(\Omega,\so)$. The energy can be written as 
    \begin{align}
    W(\pD\pu,\boA,\Curl\boA) 
    & 
    =\mu\left\Vert \sym\pD\pu\right\Vert ^{2}+\muc\!\left\Vert \skew\pD\pu-\boA\right\Vert ^{2}+\frac{\lambda}{2}\tr\!^{2}(\pD\pu)\nonumber
    \\
    & 
    \qquad+\frac{\mu\,L_{\text{c}}^{2}}{2}\left(\alpha_{1}\!\left\Vert \dev\sym\Curl\boA\right\Vert ^{2}+\alpha_{2}\!\left\Vert \skew\Curl\boA\right\Vert ^{2}+\frac{\alpha_{3}}{3}\tr\!^{2}(\Curl\boA)\right).
    \end{align}
    The equilibrium equations, in the absence of body forces, are therefore the following
    \begin{align}
    \Div\bosigma 
    & 
    =0,
    \\
    2\muc\,\skew\!(\pD\pu-\boA)-\skew\Curl\bom 
    &
     =0, \nonumber
    \end{align}
    with
    \begin{align}
    \bosigma
    & 
    =2\mu\,\sym\pD\pu+\lambda\tr\!(\pD\pu)\id+2\muc\,\skew\!(\pD\pu-\boA),
    \\
    \bom
    &
     =\mu L_{\text{c}}^{2}\left(\alpha_{1}\dev\sym\Curl\boA+\alpha_{2}\,\skew\Curl\boA+\frac{\alpha_{3}}{3}\tr\!(\Curl\boA)\id\right). \nonumber
    \end{align}
    Due to the special structure of the linear Cosserat model (rigid substructure and Cosserat couple modulus $\muc>0$), existence and uniqueness can also be shown with only $\left.\pu\right|_{\Gamma}$ prescribed at the boundary, while $\boA$ is left everywhere free, i.e., 
    \begin{equation}
    \left.\pu\right|_{\Gamma}=\widehat{u},\qquad\qquad\left.\bom\times\pnu\right|_{\partial\Omega}=0.
    \end{equation}
    This is the case because from Korn's inequality, $\pu$ is controlled in $H^1$, such that $\skew\pD\pu\in L^2$ and then $\muc>0$ allows to conclude that $A\in L^2$ and Nye's formula \eqref{Nye} shows that $\boA\in H^1(\Omega,\so)$ and the problem is strictly convex in the space $ H^1(\Omega,\bR^3)\times H^1(\Omega,\so)$.
	\section{The novel consistent coupling boundary condition on the Dirichlet part}
	The expression of the strain energy for the classical linear isotropic micromorphic
	continuum without mixed terms (like $\langle\text{sym}\boldsymbol{P},\text{sym}\left(\text{D}\boldsymbol{u}-\boldsymbol{P}\right)\rangle$,
	etc.) and simplified curvature energy \cite{rizzi2021bending,rizzi2021torsion} can be written as
	\begin{align}
		W(\pD\pu,\pP,\pD\pP)= 
		& 
		\,\mue\left\lVert \sym\!(\pD\pu-\pP)\right\rVert^{2}+\muc\left\lVert \skew\!(\pD\pu-\pP)\right\rVert^{2}
		+\dfrac{\le}{2}\text{tr}^{2}(\pD\pu-\pP)\nonumber 
		\\*
		\qquad 
		& 
		+\mum\left\lVert \sym\pP\right\rVert ^{2}+\dfrac{\lam}{2}\text{tr}^{2}(\pP)\label{eq:energy_MM}
		\\*
		\qquad 
		& 
		+\frac{\mu\,L_{\text{c}}^{2}}{2}\left(a_{1}\left\lVert \pD(\text{dev}\,\sym\pP)\right\rVert ^{2}
		+a_{2}\left\lVert \pD(\skew\pP)\right\rVert ^{2}
		+\frac{2}{9}\,a_{3}\left\lVert \pD(\text{tr}(\pP)\id)\right\rVert^{2}\right).\nonumber 
	\end{align}
	Thus, the equilibrium equations without body forces are   
	\begin{align}
	\text{Div}\,\bosigma & =\boldsymbol{0}\,,\label{eq:equiMic_MM}
	\\*
	\quad\bosigma-2\mu_{\text{micro}}\,\text{sym}\,\pP-\lambda_{\text{micro}}\text{tr}(\pP)\id+\,\DIV \gbm &=\boldsymbol{0}\,,\nonumber
	\end{align}
    where 
    \begin{align}
    \bosigma
    &  
    \coloneqq2\mue\,\sym\!(\pD\pu-\pP)
    	+2\muc\,\skew\!(\pD\pu-\pP)
    	+\le\text{tr}(\pD\pu-\pP)\id \label{eq:bos}
    \end{align}
    is the second order force stress tensor and
    \begin{align}
    \gbm
    &\coloneqq
    \mu L_{\text{c}}^{2}\,\left(a_{1}\pD(\text{dev}\,\sym\pP)+a_2\,\pD(\skew\pP)+\frac{2}{9}\,a_3\,\pD(\text{tr}(\pP)\id)\right)
    \, .
    \label{eq:m_third_MM}
    \end{align}
	is the third order moment stress tensor, while the classical system of boundary conditions is  (see Appendix \ref{Derivation new bc})
	\begin{equation}
		\begin{aligned}\begin{array}{cc}
				(\ref{eq:CBD}a) & \left.\pu\right|_{\Gamma}=\widehat{\pu}\\
				\\
				(\ref{eq:CBD}b) & \left.\pP\right|_{\Gamma}=\widehat{\pP}
			\end{array} &  &  &  & \qquad\qquad &  &  &  & \left.\begin{array}{r}
				\left.\boldsymbol{\sigma}\!\cdot\!\boldsymbol{\nu}\,\right|_{\partial\Omega\setminus\overline{\Gamma}}=\boldsymbol{0}\\
				\\
				\left.\boldsymbol{\mathfrak{m}}\!\cdot\!\boldsymbol{\nu}\,\right|_{\partial\Omega\setminus\overline{\Gamma}}=\boldsymbol{0}
			\end{array}\right\} .\quad(\ref{eq:CBD}c)\end{aligned}
		\label{eq:CBD}
	\end{equation}
    \begin{align*}
    \textrm{``strong local interation''} &\sim \pD\pP
    \\
    \textrm{``local rotational interation''} &\sim \Curl\pP
    \\
    \textrm{``weak local rotational interation''} &\sim \sym\Curl\pP
    \\
    \textrm{``very weak local rotational interation''} &\sim \dev\sym\Curl\pP &  (\textrm{conformal case \cite{lewintan2020korn}}) 
    \end{align*}
	\begin{rema}
		Let us explicitly remark that, assuming the needed regularity of the
		boundary to guarantee the existence of a local frame $\left\{ \pnu,\boldsymbol{\tau}_{1},\boldsymbol{\tau}_{2}\right\} $
		, where $\boldsymbol{\tau}_{1},\boldsymbol{\tau}_{2}$ are two orthonormal
		tangent vector fields, we have that 
		\begin{equation}
			\left.\pP\,\right|_{\Gamma}=0\qquad\Longleftrightarrow\qquad\left.\pP\!\cdot\!\pnu\,\right|_{\Gamma}=\left.\pP\!\cdot\!\boldsymbol{\tau}_{1}\,\right|_{\Gamma}=\left.\pP\!\cdot\!\boldsymbol{\tau}_{2}\,\right|_{\Gamma}=0
		\end{equation}
		or, equivalently, considering the decomposition of the identity operator
		as
		\begin{equation}
			\id=\pnu\otimes\pnu+\boldsymbol{\tau}_{1}\otimes\boldsymbol{\tau}_{1}+\boldsymbol{\tau}_{2}\otimes\boldsymbol{\tau}_{2},
		\end{equation}
		we have that
		\begin{align}
			\left.\pP\,\right|_{\Gamma}=0\qquad\Longleftrightarrow\qquad 0 
			&
			=\left.\pP\!\cdot\!\id\,\right|_{\Gamma} =\left.\pP\!\cdot\!\left(\pnu\otimes\pnu+\boldsymbol{\tau}_{1}\otimes\boldsymbol{\tau}_{1}+\boldsymbol{\tau}_{2}\otimes\boldsymbol{\tau}_{2}\right)\right|_{\Gamma}
			\\
			& =\left.\pP\!\cdot\!\left(\pnu\otimes\pnu\right)\right|_{\Gamma}+\left.\pP\!\cdot\!\left(\boldsymbol{\tau}_{1}\otimes\boldsymbol{\tau}_{1}\right)\right|_{\Gamma}+\left.\pP\!\cdot\!\left(\boldsymbol{\tau}_{2}\otimes\boldsymbol{\tau}_{2}\right)\right|_{\Gamma}\nonumber 
			\\
			& =\left.\left(\pP\!\cdot\!\pnu\right)\otimes\pnu\right|_{\Gamma}+\left.\left(\pP\!\cdot\!\boldsymbol{\tau}_{1}\right)\otimes\boldsymbol{\tau}_{1}\right|_{\Gamma}+\left.\left(\pP\!\cdot\!\boldsymbol{\tau}_{2}\right)\otimes\boldsymbol{\tau}_{2}\right|_{\Gamma}\nonumber 
			\\
			& =\left.\pP\!\cdot\!\begin{pmatrix}\pnu \left. \right|   \boldsymbol{\tau}_{1}  \left. \right|  \boldsymbol{\tau}_{2}\end{pmatrix}\right|_{\Gamma}.\nonumber 
		\end{align}
		Moreover, the following relation holds
		\begin{equation*}
			\left.\pP\!\cdot\!\boldsymbol{\tau}_{1}\,\right|_{\Gamma}=\left.\pP\!\cdot\!\boldsymbol{\tau}_{2}\,\right|_{\Gamma}=0\qquad\Longleftrightarrow\qquad\left.\pP\times\pnu\,\right|_{\Gamma}=0. \vspace{-7mm}
		\end{equation*}
	\end{rema}
	
	Existence and uniqueness is easily established in the space $(\pu,\pP)\in H^{1}(\Omega,\bR^{3})\times H^{1}(\Omega,\bR^{3\times3})$
	including the case $\muc\geqslant0$, since if $\muc=0$ the full
	gradient of $\pP$ allows, via Poincaré's inequality and using
	(\ref{eq:CBD}b) to control $\pP\in L^{2}(\Omega)$ even
	though $\skew\pP$ is not controlled locally, and this is true even
	if $\mum=0$.\footnote{Quite strange!}
	
	Let us consider a Lipschitz domain $\Omega$
	and let $\Gamma$ be an open subset of the boundary $\partial\Omega$.
	We denote by 
	\begin{equation}
		\begin{aligned}H_{0,\Gamma}^{1}\left(\Omega,\bR^{n\times n}\right) & \defi\overline{C_{0,\Gamma}^{\infty}\left(\overline{\Omega},\bR^{n\times n}\right)}^{\,\left\Vert \cdot\right\Vert _{1,2,\Omega}} & \textrm{where}\\
			C_{0,\Gamma}^{\infty}\left(\overline{\Omega},\bR^{n\times n}\right) & \defi\left\{ \pP\in C^{\infty}(\overline{\Omega},\bR^{n\times n})\;\left.\right|\;\mathsf{dist}\left(\supp\pP,\Gamma\right)>0\right\} 
		\end{aligned}
		\label{eq:Null-trace space}
	\end{equation}
	the classical Sobolev space of functions with zero-trace over $\Gamma$.
	
	The classical boundary conditions for the micromorphic model (\ref{eq:energy_MM})
	presuppose the knowledge of $\left.\pP\right|_{\Gamma}$, which can be
	difficult to motivate from a modelling perspective in some cases.
	Typically, at a ``glued'' part of the body
	$\Gamma$ we know the displacement $\left.\pu\right|_{\Gamma}$. This implies the knowledge
	of the tangential part of $\text{D}\pu$, namely
	\begin{equation}
		\left.\text{D}\pu\times\pnu\,\right|_{\Gamma}.\label{eq:Tangential gradient of u}
	\end{equation}
	Accordingly, it is possible to impose the \underline{consistent coupling
		boundary condition} for $\pP$
	\begin{equation}
		\left.\pP\times\pnu\,\right|_{\Gamma}=\left.\text{D}\pu\times\pnu\,\right|_{\Gamma}.\label{eq:Consistent coupling boundary condition}
	\end{equation}
	Clearly, if $\pP\in H^{1}\!\left(\Omega,\bR^{3\times3}\right)$, we have
	the possibility to evaluate $\left.\pP\times\pnu\right|_{\Gamma}$
	and/or $\left.\pP\!\cdot\!\pnu\right|_{\Gamma}$.

\begin{figure}[H]
		\begin{centering}
			\includegraphics[scale=0.6]{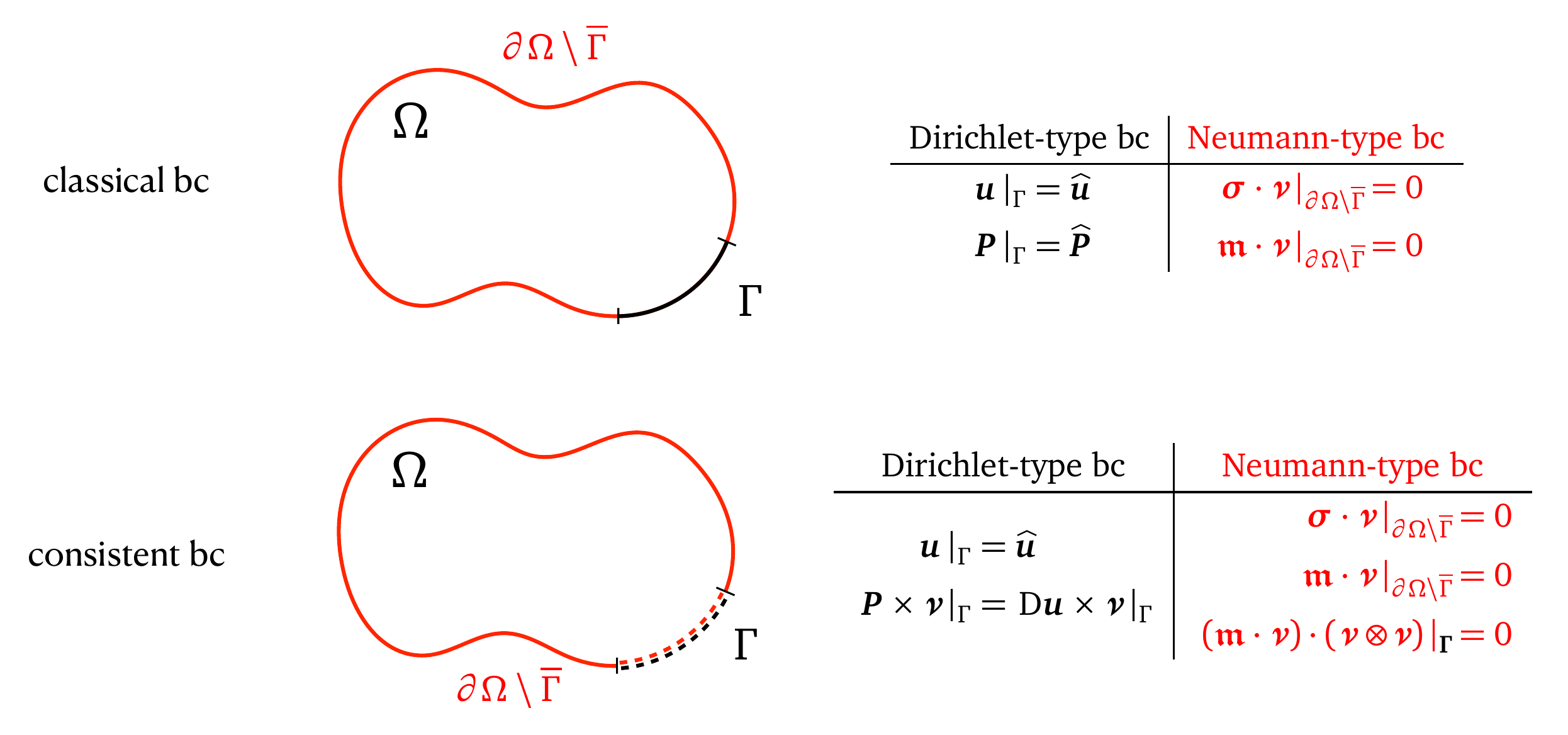}
			\par\end{centering}
		\caption{Comparing the classical Dirichlet type boundary condition with the new condition of consistent coupling. The issue is that at the standard Dirichlet part $\Gamma$ there is also a Neumann type boundary condition appearing for consistent coupling.}
		\label{fig:CCBCvsDBCvsNBC}
	\end{figure}

	 Thus, we can study (\ref{eq:energy_MM}) with a weaker
	set of boundary conditions (see Fig.~\ref{fig:CCBCvsDBCvsNBC}), i.e., 
	\begin{align}
		\left.\pu\right|_{\Gamma} & =\widehat{\pu} & {\color{red}\left.\bosigma\!\cdot\!\pnu\,\right|_{\partial\Omega\setminus\overline{\Gamma}}}\, & {\color{red}=0}\nonumber 
		\\
		\left.\pP\times\pnu\right|_{\Gamma} & =\left.\text{D}\pu\times\pnu\right|_{\Gamma} & {\color{red}\left.\gbm\!\cdot\!\pnu\,\right|_{\partial\Omega\setminus\overline{\Gamma}}} \,& {\color{red}=0}\label{eq:Consistent boundary conditions}
		\\
		&  & {\color{red}\left.\left(\gbm\!\cdot\!\pnu\right)\!\cdot\!\left(\pnu\otimes\pnu\right)\,\right|_{\color{black}{\boldsymbol{\Gamma}}}} \,& {\color{red}=0\quad{\color{black}\left(\textrm{new Neumann-type bc}\right)}{\color{black}.}}\nonumber 
	\end{align}
	The new appearing boundary condition $\left(\boldsymbol{\mathfrak{m}}\!\cdot\!\boldsymbol{\nu}\right)\!\cdot\!\left(\boldsymbol{\nu}\otimes\boldsymbol{\nu}\right)$ at the Dirichlet portion of the boundary $\Gamma$ will be derived in the Appendix.
	
\begin{figure}[H]
		\centering
		\begin{subfigure}{0.49\textwidth}
			\centering
			\includegraphics[width=\textwidth]{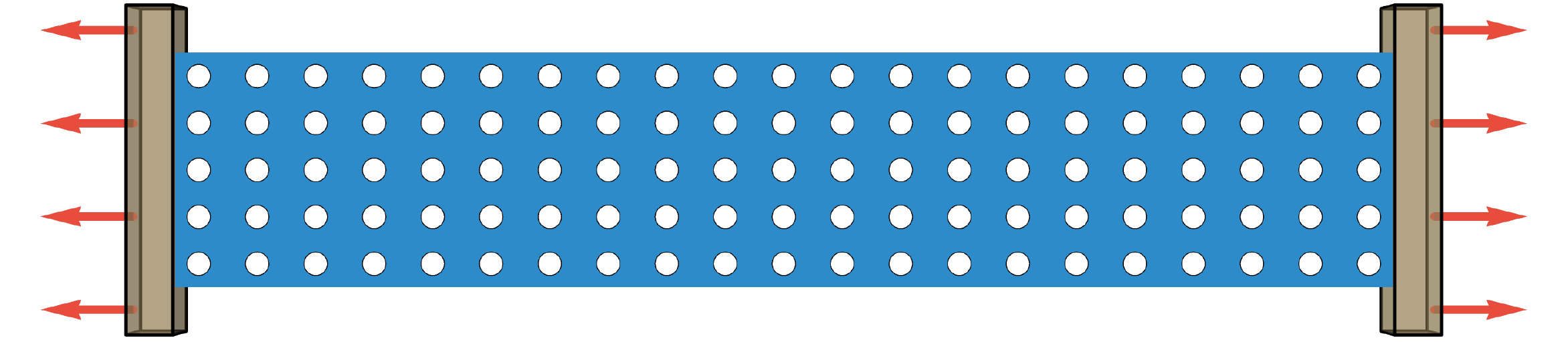}
		\end{subfigure}
		\hfill
		\begin{subfigure}{0.49\textwidth}
			\centering
			\includegraphics[width=\textwidth]{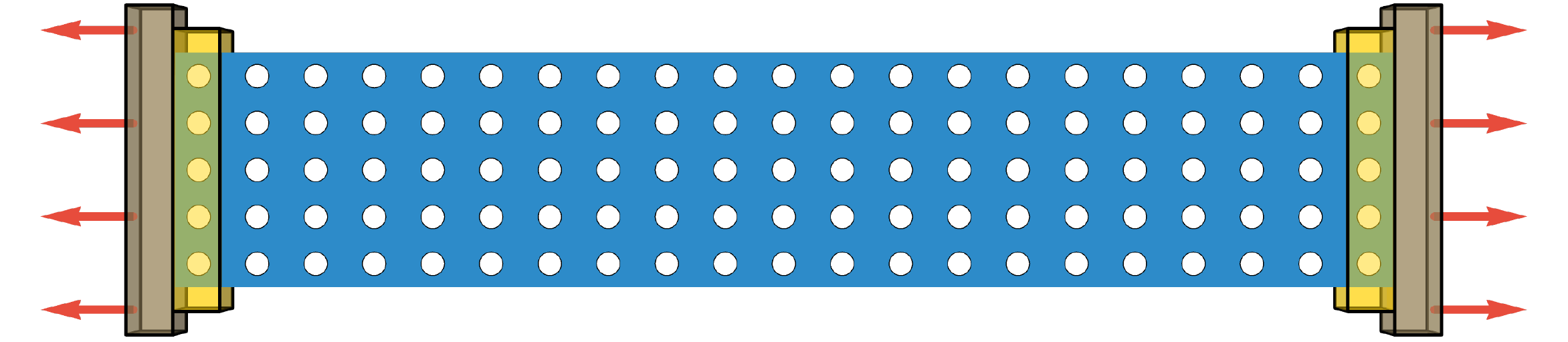}
		\end{subfigure}
		\caption{Sketch of a 2D-extensional test in which (\textit{left}) the consistent coupling boundary conditions should be applied due to a ``glued-like'' surface constraint and one in which (\textit{right}) the full Dirichlet boundary conditions should be applied due to the presence of a fully clamped constraint.
		Also here, in the context of a second gradient model, it is particularly clear that the choice of having the component of the gradient normal to the surface equal to zero in the clamped case (\textit{right}) is due to the fact that again, a stripe of material cannot move in the normal direction to the surface, while in the ``glued-like'' case (\textit{left}), the component of the gradient normal to the surface is free.
		}
		\label{fig:beam_bc}
	\end{figure}

\begin{figure}[H]
	\centering
	\begin{subfigure}{0.45\textwidth}
		\centering
		\includegraphics[width=\textwidth]{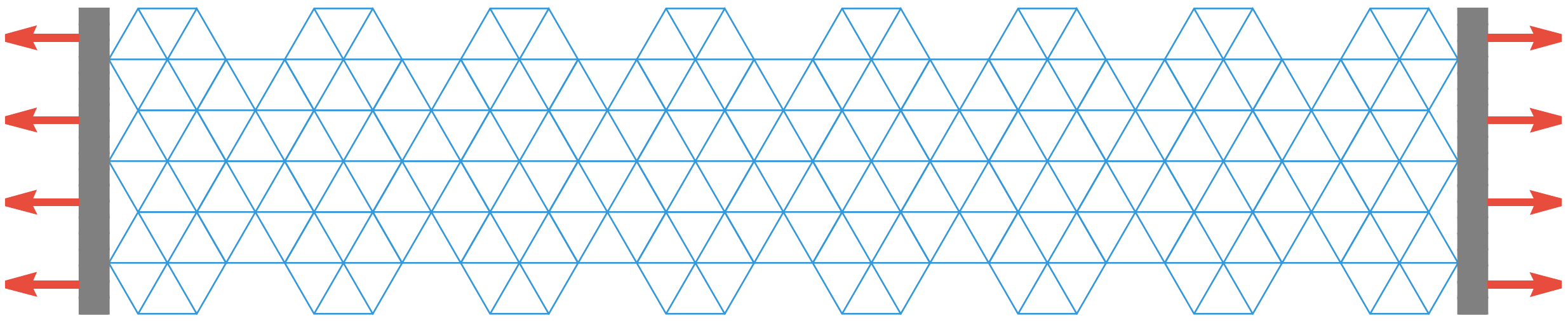}
	\end{subfigure}
	\hspace{0.8cm}
	\begin{subfigure}{0.45\textwidth}
		\centering
		\includegraphics[width=\textwidth]{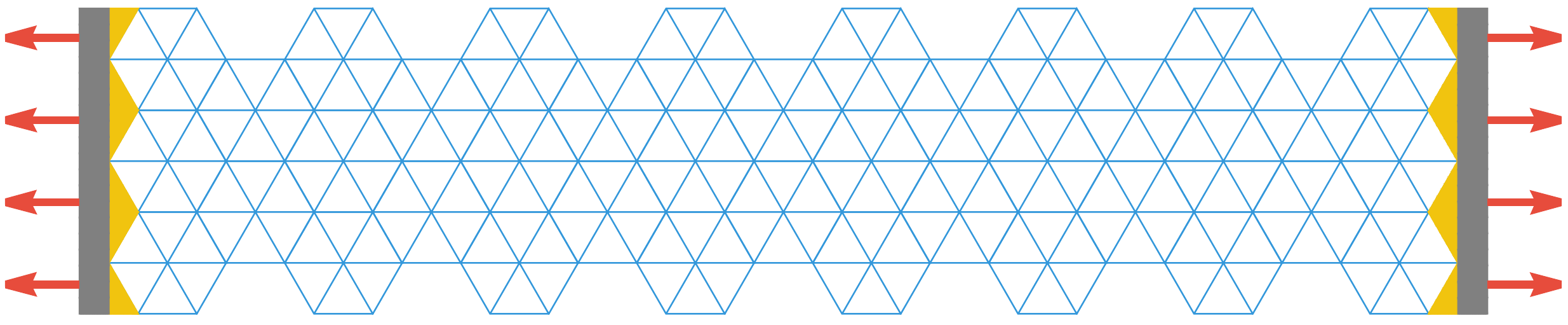}
	\end{subfigure}
	\caption{
		Sketch of a 2D-extensional test on a lattice structure in which (\textit{left}) the consistent coupling boundary conditions should be applied due to a ``hinged-like'' constraint and one in which (\textit{right}) the full Dirichlet boundary conditions should be applied due to the presence of a fully clamped constraint.
		In the context of  a second gradient model, it is particularly clear that the justification for the choice of having the component of the gradient normal to the surface equal to zero in the clamped case (\textit{right}) is due to the fact that two consecutive points cannot move in the normal direction to the surface, while in the ``hinged-like'' case (\textit{left}), the component of the gradient normal to the surface is free.
	}
	\label{fig:lattice_bc}
	\end{figure}

	In the two following sketches (Fig.~\ref{fig:beam_bc}-\ref{fig:lattice_bc}) it is shown how the choice of the boundary conditions (consistent coupling or full Dirichlet) can be related to the constraints that may be applied to a real structure (solid or lattice).
	

	\section{Existence and uniqueness of minimizers for consistent coupling boundary conditions for the classical micromorphic model}
	
	The existence and uniqueness proof for the classical micromorphic model with classical boundary conditions  is based on showing that the minimization
	problem
	\[
	\int_{\Omega}W\left(\text{D}\boldsymbol{u},\boldsymbol{P},\text{D}\boldsymbol{P}\right)\text{dV}\longrightarrow\min,\qquad\left(\pu,\pP\right)\in H_{0,\Gamma}^{1}(\Omega,\bR^{3})\times H_{0,\Gamma}^{1}(\Omega,\bR^{3\times3}),
	\]
	with homogeneous boundary conditions on $\Gamma$ for $\bou$ and $\boP$, is strictly convex
	and coercive in the space $H_{0,\Gamma}^{1}(\Omega,\bR^{3})\times H_{0,\Gamma}^{1}(\Omega,\bR^{3\times3})$. In order to study the same problem for the new consistent boundary condition we need to introduce a new space of functions replacing $H_{0,\Gamma}^{1}(\Omega,\bR^{3\times3})$. 
	Let us consider the Hilbert-space
	\begin{equation}
		\mathscr{H}^{\,\sharp}(\Omega)\defi\left\{ \pP\in H^{1}(\Omega,\bR^{3\times3})\;\left. \right| \;\left.\pP\times\pnu\right|_{\Gamma}=0\right\} \label{eq: New space}
	\end{equation}
	equipped with the norm
	\[
	\left\Vert \pP\right\Vert _{\sharp}^{2}\defi\left\Vert \sym\pP\right\Vert _{\bR^{3\times3}}^{2}+\left\Vert \text{D}\pP\right\Vert _{\bR^{3\times3\times3}}^{2}.
	\]
	\begin{rema}\label{rem_norm}
		Since 
		\[
		\left.\skew\pP\times\pnu\right|_{\Gamma}=0\quad\Longleftrightarrow\quad\left.\skew\pP\right|_{\Gamma}=0,
		\]
		it is clear that
		$\left\Vert \cdot\right\Vert _{\sharp}^{2}$
		is a norm on $\mathscr{H}^{\,\sharp}$.
	\end{rema}
	\begin{rema}
		If, instead of considering the tangent part of $\pP$ on the boundary we account for the normal part, i.e., we require $\left.\pP\!\cdot\!\pnu \right|_\Gamma=0 $, the introduced norm $\left\| \cdot \right\|_\sharp $ is not anymore a norm on the corresponding vector space (it will be simply a seminorm). In other words, setting
		\begin{equation}
			\mathscr{H}^{\,\flat}(\Omega)\defi\left\{ \pP\in H^{1}(\Omega,\bR^{3\times3})\;\left.\right|\;\left.\pP\!\cdot\!\pnu\,\right|_{\Gamma}=0\right\} 
		\end{equation}
		the resulting couple $\left(\mathscr{H}^{\,\flat}(\Omega),\left\Vert \cdot\right\Vert _{\sharp}\right)$ is not a normed space because the positive definiteness of $\left\Vert \cdot\right\Vert _{\sharp}$ fails. Indeed,
		\begin{equation}
			\left.\begin{array}{r}
				\sym\pP=0\\
				\\
				\text{D}\pP=0
			\end{array}\right\} \quad\Rightarrow\quad \pP(x)=\boA\in\so\;\;\textrm{constant}\qquad\textrm{and}\qquad\left.\boA\!\cdot\!\pnu\,\right|_{\Gamma}=0\quad\nRightarrow\quad \boA=0.
		\end{equation}
		The reason is related to the property established in \eqref{rule0}: it suffices to consider $\boA=\Anti\pnu\neq 0$ to obtain
		\begin{equation*}
			\boA\!\cdot\!\pnu=\left(\Anti\pnu\right)\!\cdot\!\pnu=\pnu\times\pnu=0,\qquad \forall \,\pnu\neq0\;\textrm{unit vector}. \vspace{-7mm}
		\end{equation*}
	\end{rema}
	\begin{rema}
		We can also impose the consistent coupling boundary condition in the linearized Cosserat model ( Section \ref{Cosserat model}). In this case, only the skew-symmetric part $\boA=\skew\pP$ of $\pP$ intervenes and the homogeneous consistent coupling boundary condition on $\Gamma$ reads accordingly
		\[
			\left.\pu\,\right|_{\Gamma}=0,\qquad\qquad\left.\skew\text{D} \pu\times\pnu\,\right|_{\Gamma}=\left.\boA\times\pnu\,\right|_{\Gamma},
		\]
		giving
		\[
		\left.\pu\,\right|_{\Gamma}=0\quad\Longrightarrow\quad\left.\skew\text{D}  \pu\times\pnu\,\right|_{\Gamma}=0\quad\Longrightarrow\quad\left.\boA\times\pnu\,\right|_{\Gamma}=0\quad\Longleftrightarrow\quad \boA=0.
		\]	
		This implies that in the Cosserat model with rigid microstructure, a fully prescribed condition on $\pu$ at $\Gamma $ implies full clamping of $\boA$ at $\Gamma $ as well if consistent coupling is required.
	\end{rema}
	\begin{rema}
		The generalized Cauchy-type force stress tensor $\bosigma$ becomes symmetric if and only if $\muc=0$ (in other words, the Cosserat couple modulus is responsible for the loss of symmetry of the force stress tensor $\bosigma$).
	\end{rema}
 	\noindent Returning to the micromorphic case, the new minimization problem hence reads
	\begin{equation}
		\int_{\Omega}W(\text{D}\boldsymbol{u},\pP,\text{D}\pP)\,\text{dV}\longrightarrow\min,\qquad\left(\pu,\pP\right)\in H_{0,\Gamma}^{1}(\Omega,\bR^{3})\times \mathscr{H}^{\,\sharp}(\Omega).
	\end{equation}
	Denoting by $\left\Vert \cdot\right\Vert _{1,2,\Omega}$ the classical $H^1$-Sobolev norm, we would like to prove the norm-equivalence of spaces
	\begin{equation}\label{Equivalence of spaces}
		\left(\mathscr{H}^{\,\sharp}(\Omega),\left\Vert \cdot\right\Vert _{\sharp}\right)\simeq\left(\mathscr{H}^{\,\sharp}(\Omega),\left\Vert \cdot\right\Vert _{1,2,\Omega}\right)
	\end{equation}
	and with an abuse of notation we will denote $\left(\mathscr{H}^{\,\sharp}(\Omega),\left\Vert \cdot\right\Vert _{\sharp}\right)$ simply by $\mathscr{H}^{\,\sharp}(\Omega)$.

	\subsection{Failure of a simple direct existence argument}
	
	One is tempted to think that showing the equivalence of spaces given in \eqref{Equivalence of spaces}
	is an easy exercise. However, we will exhibit where the simple argument
	fails. To proceed, let us observe first that
	\begin{align}
		\left\Vert \text{D}\pP\right\Vert _{\bR^{3\times3\times3}}^{2} & \defi\sum_{i=1}^{3}\left\Vert \partial_{i}\pP\right\Vert _{\bR^{3\times3}}^{2}=\sum_{i=1}^{3}\left\Vert \partial_{i}\left(\sym\pP+\skew\pP\right)\right\Vert ^{2} =\sum_{i=1}^{3}\left\Vert \partial_{i}\sym\pP+\partial_{i}\skew\pP\right\Vert ^{2}\label{eq:Splitting on the gradient}\\
		& =\sum_{i=1}^{3}\left\Vert \partial_{i}\sym\pP\right\Vert ^{2}+\left\Vert \partial_{i}\skew\pP\right\Vert ^{2}+2\underbrace{\left\langle \partial_{i}\sym\pP,\partial_{i}\skew\pP\right\rangle }_{=\,0}=\sum_{i=1}^{3}\left\Vert \partial_{i}\sym\pP\right\Vert ^{2}+\left\Vert \partial_{i}\skew\pP\right\Vert ^{2}\nonumber\\
		&=\left\Vert \text{D}\,\sym\pP\right\Vert ^{2}+\left\Vert \text{D}\,\skew\pP\right\Vert ^{2}.\nonumber
	\end{align}
	Therefore $\left\Vert \text{D}\pP\right\Vert^{2}$ splits properly into $\left\Vert \text{D}\,\sym\pP\right\Vert^{2}+\left\Vert \text{D}\,\skew\pP\right\Vert^{2}$. Moreover,
	\begin{equation}
		\begin{aligned}
		\left\Vert \sym\pP\right\Vert ^{2}+\left\Vert \text{D}\pP\right\Vert ^{2}=0\quad\Longrightarrow\quad & \pP=\boldsymbol{A}\in\so\quad a.e.
		\\
		\Longrightarrow\quad & \!\!\left.\pP\times\pnu\right|_{\Gamma}=\left.\boldsymbol{A}\times\pnu\right|_{\Gamma}=0 & \Longrightarrow\quad\boldsymbol{A}=0,
		\end{aligned}
		\label{eq:norm validity}
	\end{equation}
	showing Remark \ref{rem_norm}. Since $\sym\pP$ is already controlled in the norm $\left\| \cdot\right\|_\sharp $ it remains to control $\skew\pP$. The idea would be to use the control of $\left\| \text{D}\,\skew\pP\right\|^2 $ from \eqref{eq:Splitting on the gradient} in this respect.
	This may lead us to the question whether the classical Poincaré inequality
	for $\skew\pP$, which reads 
	\begin{equation}
		\int_{\Omega}\left\Vert \text{D}\,\skew\pP\right\Vert ^{2}dx\geqslant C_{\textrm{Poincaré}}^{+}\int_{\Omega}\left\Vert \skew\pP\right\Vert ^{2}dx,\qquad\left.\skew\pP\times\pnu\right|_{\Gamma}=\left.\skew\pP\right|_{\Gamma}=0,\label{eq:classical Poincar=0000E9 inequality on P}
	\end{equation}
	could be extended to the case when full boundary conditions on $\pP$ (but not on $\skew\pP$)
	are prescribed: is it true that 
	\begin{equation}
		\int_{\Omega}\left\Vert \text{D}\,\skew\pP\right\Vert ^{2}dx\geqslant C_{\textrm{Poincaré}}^{+}\int_{\Omega}\left\Vert \skew\pP\right\Vert ^{2}dx,\qquad\left.\pP\times\pnu\right|_{\Gamma}=0\,?\label{eq:Extended Poincar=0000E9 inequality on skewP}
	\end{equation}
	The answer is again no. Indeed, it is sufficient to consider the constant
	field 
	\begin{equation}
		\pP=\begin{pmatrix}a & 0 & 0\\
			b & 0 & 0\\
			c & 0 & 0
		\end{pmatrix}\qquad\textrm{with }a,b,c\in\bR\setminus\left\{ 0\right\} .
	\end{equation}
	In this case
	\begin{equation}
		\skew\pP=\frac{1}{2}\begin{pmatrix}0 & -b & -c\\
			b & 0 & 0\\
			c & 0 & 0
		\end{pmatrix},\qquad {\rm D} \,\skew\pP=0\quad\textrm{and}\quad\left.\pP\times\pnu\right|_{\Gamma}=0\label{eq:cc}
	\end{equation}
	but
	\begin{equation}
		\left.\skew\pP\times\pnu\right|_{\Gamma}=\frac{1}{2}\begin{pmatrix}0 & -c & -b\\
			0 & 0 & 0\\
			0 & 0 & 0
		\end{pmatrix}\neq0
		\qquad
		{\rm with}
		\qquad
		 \pnu=\boe_1=(1,0,0)\, .
		 \label{eq:ccc}
	\end{equation}
	Therefore, using (\ref{eq:Splitting on the gradient}) together with
	Poincaré's inequality cannot lead to a valid argument. We will next
	resort to a highly nontrivial result on coercive inequalities for incompatible tensor fields  \cite{neff2015poincare,lewintan2021p,lewintan2020korn,lewintan2019ne}.
	
	\begin{rema}
		The previous argument fails because  
		\begin{equation}
		\left.\pP\times\pnu\right|_{\Gamma}=\left.\sym\pP\times\pnu\right|_{\Gamma}+\left.\skew\pP\times\pnu\right|_{\Gamma}=0\quad\nRightarrow\quad\left.\skew\pP\times\pnu\right|_{\Gamma}=0,\label{eq:impli3}
		\end{equation}
		i.e., the cross product with $\times\pnu$ on the boundary does not preserve the orthogonal split into symmetric and skew-symmetric parts.
		Indeed, if we choose again $\pnu=\boe_1=(1,0,0)$ we obtain
		\begin{equation}
		\sym\pP\times\pnu=\begin{pmatrix}0 & \frac{1}{2}P_{\left(13\right)} & -\frac{1}{2}P_{\left(12\right)}\vphantom{\Big| }\\
		0 & \frac{1}{2}P_{\left(23\right)} & -P_{22}\vphantom{\Big| }\\
		0 & P_{33} & -\frac{1}{2}P_{\left(23\right)}\vphantom{\Big| }
		\end{pmatrix}\qquad\textrm{and}\qquad\skew\pP\times\pnu=\frac{1}{2}
		\begin{pmatrix}
		0 & P_{\left[13\right]} & -P_{\left[12\right]}\vphantom{\Big| }\\
		0 & P_{\left[23\right]} & 0\vphantom{\Big| }\\
		0 & 0 & P_{\left[23\right]} \vphantom{\Big| }
		\end{pmatrix}\label{eq:skew and sym}
		\end{equation}
		and if we solve with respect to the components of $\pP$ the equation
		\begin{equation}
		\left.\sym\pP\times\pnu\right|_{\Gamma}+\left.\skew\pP\times\pnu\right|_{\Gamma}=0\quad\Longleftrightarrow\quad\left.\sym\pP\times\pnu\right|_{\Gamma}=-\left.\skew\pP\times\pnu\right|_{\Gamma},
		\end{equation}
		we obtain
		\begin{equation}
		P_{12}=P_{13}=P_{22}=P_{23}=P_{32}=P_{33}=0,
		\end{equation}
		i.e., every 
		\[
		\pP=\begin{pmatrix}P_{11} & 0 & 0\\
		P_{21} & 0 & 0\\
		P_{31} & 0 & 0
		\end{pmatrix}\qquad\textrm{solves}\qquad\left.\sym\pP\times\pnu\right|_{\Gamma}=-\left.\skew\pP\times\pnu\right|_{\Gamma}
		\]
		but 
		\begin{equation}
		\skew\begin{pmatrix}P_{11} & 0 & 0\\
		P_{21} & 0 & 0\\
		P_{31} & 0 & 0
		\end{pmatrix}=\frac{1}{2}\begin{pmatrix}0 & -P_{21} & -P_{31}\\
		P_{21} & 0 & 0\\
		P_{31} & 0 & 0
		\end{pmatrix}\label{eq:bb}
		\end{equation}
		and
		\begin{equation}
		\frac{1}{2}\begin{pmatrix}0 & -P_{21} & -P_{31}\\
		P_{21} & 0 & 0\\
		P_{31} & 0 & 0
		\end{pmatrix}\times\pnu=\frac{1}{2}\begin{pmatrix}0 & -P_{31} & -P_{21}\\
		0 & 0 & 0\\
		0 & 0 & 0
		\end{pmatrix}\neq0.\label{eq:bbb}
		\end{equation}
		Via the map $\Anti:\bR^{3}\fr\so$, we can rewrite the problem as
		follows:
		
		\medskip
		
		$
		\qquad\qquad\left.\sym\pP\times\pnu\right|_{\Gamma}+\left.\skew\pP\times\pnu\right|_{\Gamma}=0
		\qquad\overset{\eqref{new id}}{\Longleftrightarrow}\qquad\left.\sym\pP\!\cdot\!\Anti\pnu\right|_{\Gamma}+\left.\skew\pP\!\cdot\!\Anti\pnu\right|_{\Gamma}=0. \nonumber
		$
	\end{rema}
	
	\subsection{Existence argument for consistent coupling via the incompatible Korn's inequality}
	Let us remember that in $H(\Curl;\Omega,\bR^{3\times3})$ the tangential trace $\pP\times\pnu$ is naturally defined.
	In \cite{neff2015poincare,lewintan2021p,lewintan2020korn,lewintan2019ne,neff2011canonical} the following generalized incompatible Korn's inequality has been established.
	\begin{thm}\label{Korn}
		There exists a constant $c^{+}>0\;\textrm{such that for all}\;\pP\in H\!\left(\Curl;\Omega,\bR^{3\times3}\right),\;\textrm{with}\;\left.\pP\times\pnu\right|_{\Gamma}=0$
		\begin{equation}
			\left\Vert \pP\right\Vert _{L^{2}\left(\Omega\right)}^{2}\leqslant c^{+}\left(\left\Vert \sym\pP\right\Vert _{L^{2}\left(\Omega\right)}^{2}+\left\Vert \Curl\pP\right\Vert _{L^{2}\left(\Omega\right)}^{2}\right).\label{eq:New Korn}
		\end{equation}
	\end{thm}
	We can use this estimate for our case next. Indeed, adding $\left\Vert \text{D}\pP\right\Vert _{L^{2}\left(\Omega\right)}^{2}$
	to both sides in inequality (\ref{eq:New Korn}) we have
	\begin{align}
		\left\Vert \pP\right\Vert _{L^{2}\left(\Omega\right)}^{2}+\left\Vert \text{D}\pP\right\Vert _{L^{2}\left(\Omega\right)}^{2} 
		& \overset{\left.\pP\times\pnu\right|_{\Gamma}=\,0}{\leqslant}c^{+}\left(\left\Vert \sym\pP\right\Vert _{L^{2}\left(\Omega\right)}^{2}+\left\Vert \Curl\pP\right\Vert _{L^{2}\left(\Omega\right)}^{2}\right)+\left\Vert \text{D}\pP\right\Vert _{L^{2}\left(\Omega\right)}^{2}\nonumber 
		\\
		& \underset{{\phantom{\left.\pP\times\pnu\right|_{\Gamma}=\,0}}}{\leqslant}\max\left\{ 1,c^{+}\right\} \left(\left\Vert \sym\pP\right\Vert _{L^{2}\left(\Omega\right)}^{2}+\left\Vert \Curl\pP\right\Vert _{L^{2}\left(\Omega\right)}^{2}+\left\Vert \text{D}\pP\right\Vert _{L^{2}\left(\Omega\right)}^{2}\right)\label{eq:Korn gen}
		\\
		& \overset{{\phantom{\left.\pP\times\pnu\right|_{\Gamma}=\,0}}}{\overset{\eqref{Curl Bound}}{\leqslant}}\max\left\{ 1,c^{+}\right\} \left(\left\Vert \sym\pP\right\Vert _{L^{2}\left(\Omega\right)}^{2}+3\left\Vert \text{D}\pP\right\Vert _{L^{2}\left(\Omega\right)}^{2}\right)\nonumber 
		\\
		& \overset{{\phantom{\left.\pP\times\pnu\right|_{\Gamma}=\,0}}}{\leqslant}\underbrace{3\max\left\{ 1,c^{+}\right\} }_{c_{1}}\left(\left\Vert \sym\pP\right\Vert _{L^{2}\left(\Omega\right)}^{2}+\left\Vert \text{D}\pP\right\Vert _{L^{2}\left(\Omega\right)}^{2}\right).\nonumber 
	\end{align}
	This shows  that
	\begin{equation}
		\left(\mathscr{H}^{\,\sharp}(\Omega),\left\Vert \cdot\right\Vert _{\sharp}\right)\simeq\left(\mathscr{H}^{\,\sharp}(\Omega),\left\Vert \cdot\right\Vert _{1,2,\Omega}\right)
	\end{equation}
	and leads to existence and uniqueness of the classical micromorphic model with consistent coupling condition following the direct methods of the calculus of variations for a strictly convex problem in $H^1(\Omega,\bR^3)\times\mathscr{H}^{\,\sharp}(\Omega)$.
	\begin{rema}
	    The consistent coupling boundary condition can be weakened imposing only either 
	    $$
	    \sym\!(\!\left.\pP\times\pnu\right|_{\Gamma}) =\sym\!(\!\left.\text{D}\pu\times\pnu\right|_{\Gamma})
	    $$
	    or 
	    $$
	    \dev\sym\!(\!\left.\pP\times\pnu\right|_{\Gamma}) =\dev\sym\!(\!\left.\text{D}\pu\times\pnu\right|_{\Gamma})
	    $$
	    on $\Gamma$ (see for details Appendix \ref{GC}). 
	\end{rema}
	\section{Comparison of analytical solutions for the different boundary conditions}
\begin{figure}[H]
		\begin{centering}
			\begin{tikzpicture}
				\node (M1) [circle] at (2,-1) {};
				\filldraw[gray,pattern=north east lines] (0,0) -| (0.2,-2) -|  (0,0);
				\draw [ snake=coil, red, segment amplitude=5pt, segment length=5pt] (0.2,-1) -- (M1);
				\draw [ snake=coil, blue, segment amplitude=7pt, segment length=6pt] (1.82,-1) -- (3.8,-1);
				\draw [thick,decoration={brace,amplitude=6pt,mirror,raise=0.3cm},decorate] (0.2,-1) -- (1.7,-1) 
				node [pos=0.5,anchor=north,yshift=-0.55cm] {$\me$}; 
				\draw [thick,decoration={brace,amplitude=6pt,mirror,raise=0.3cm},decorate] (1.82,-1) -- (3.8,-1) 
				node [pos=0.5,anchor=north,yshift=-0.55cm] {$\mum$}; 
				\draw [thick,decoration={brace,amplitude=8pt,mirror,raise=0.3cm},decorate] (0.2,-1.8) -- (3.8,-1.8) 
				node [pos=0.5,anchor=north,yshift=-0.55cm] {$\muM$}; 
				\node (M2) [circle,xshift=6cm] at (2,-1) {};
				\filldraw[gray,pattern=north east lines,xshift=6cm] (0,0) -| (0.2,-2) -|  (0,0);
				\draw [ snake=coil, red, segment amplitude=5pt, segment length=5pt,xshift=6cm] (0.2,-1) -- (M2);
				\draw [ snake=coil, blue, segment amplitude=7pt, segment length=6pt,xshift=6cm] (1.82,-1) -- (3.8,-1);
				\draw [thick,decoration={brace,amplitude=6pt,mirror,raise=0.3cm},decorate,xshift=6cm] (0.2,-1) -- (1.7,-1) 
				node [pos=0.5,anchor=north,yshift=-0.55cm] {{$\kappa_{e}$}}; 
				\draw [thick,decoration={brace,amplitude=6pt,mirror,raise=0.3cm},decorate,xshift=6cm] (1.82,-1) -- (3.8,-1) 
				node [pos=0.5,anchor=north,yshift=-0.55cm] {{$\kappa_{\textrm{micro}}$}}; 
				\draw [thick,decoration={brace,amplitude=8pt,mirror,raise=0.3cm},decorate,xshift=6cm] (0.2,-1.8) -- (3.8,-1.8) 
				node [pos=0.5,anchor=north,yshift=-0.55cm] {$\kappa_{\textrm{macro}}$}; 
				\node[align=left] at (8.5,-4) {\small $\begin{aligned}
						\kappa_{\rm macro}&=\kappa_{\rm micro}\,(\kappa_{\rm micro}+\kappa_{\rm e})^{-1}\kappa_{\rm e} \\ 
						\kappa_{\rm e}&=\kappa_{\rm micro}\,(\kappa_{\rm micro}-\kappa_{\rm macro})^{-1}\kappa_{\rm macro}\\
						\kappa_{\rm micro}&=\kappa_{\rm e}\,(\kappa_{\rm e}-\kappa_{\rm macro})^{-1}\kappa_{\rm macro}
					\end{aligned}$};
				\node[draw=red,align=left,inner sep=1.5ex,thick] at (2.2,-4) {\small $\begin{aligned}
						\muM&=\mum\,(\mum+\mu_{\text{e}})^{-1}\mu_{\text{e}} \\ 
						\mu_{\text{e}}&=\mu_{\textrm{micro}}\,(\mu_{\textrm{micro}}-\mu_{\textrm{macro}})^{-1}\mu_{\textrm{macro}}\\
						\mu_{\textrm{micro}}&=\mu_{\text{e}}\,(\mu_{\text{e}}-\mu_{\textrm{macro}})^{-1}\mu_{\textrm{macro}}
					\end{aligned}$} ;
				\draw[->, thick,gray] (3.8,-1) -- (4.6,-1);
				\draw[->, thick,gray] (9.8,-1) -- (10.6,-1);
				\node[] at (4.2,-0.7) {$F$};
				\node[] at (10.2,-0.7) {$F$};
			\end{tikzpicture}
			\par\end{centering}
		\caption{For all micromorphic models without mixed terms in the lower order energy part in \eqref{eq:energy_MM} the appearing material parameters are related by a Reuss-like homogenization formula \cite{neff2020identification}. The scale-independent response is governed by two springs in series. If $\mum=\muM$ (or $\kappa_{\rm micro}=\kappa_{\rm macro}$) then $\me=+\infty$ (or $\kappa_{\rm e}=+\infty$). The macroscopic stiffness is determined by classical homogenization \cite{neff2020identification}. The Cosserat couple modulus $\muc\geqslant0$ is not appearing in these formulas. The macroscopic shear modulus $\muM$ and the macroscopic bulk modulus $\kappa_{\rm macro}$ are uniquely determined by classical homogenization \cite{neff2020identification}.}
		\label{fig:springs}
	\end{figure}
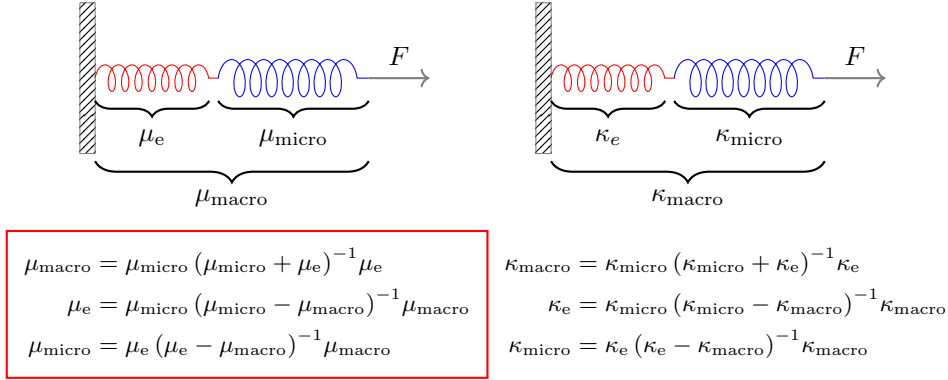

	Having established the well-posedness of the new formulation with consistent coupling, we are now comparing analytical solutions for the different possible Dirichlet boundary conditions (see Fig.~\ref{fig:stiff_shear_MM}-\ref{fig:stiff_extens_MM}-\ref{fig:stiff_shear_RM}-\ref{fig:stiff_extens_RM}-\ref{fig:stiff_shear_SG}-\ref{fig:stiff_extens_SG}). These solutions were derived in \cite{rizzi2021bending,rizzi2021torsion,rizzi2021extension}.
	For the limit $L_{\rm c}\to 0$, all the models converge to a classical Cauchy material with stiffness $\mathbb{C}_{\rm macro}$ (see Fig.~\ref{fig:springs} for more details)\footnote{In 3D, the relation between the bulk modulus and the Lamé parameters is $\kappa_i=\lambda_{i}+\frac{2}{3}\mu_i$, where $i=\{{\rm e,micro,macro}\}$.}.
	\subsection{Classical Mindlin-Eringen micromorphic-model}
	The accounted energy density is that in \eqref{eq:energy_MM} with corresponding equilibrium equations in  \eqref{eq:equiMic_MM}.
	\subsubsection{Simple shear}
	\begin{figure}[H]
		\centering
		\begin{subfigure}{0.49\textwidth}
			\centering
			\includegraphics[width=\textwidth]{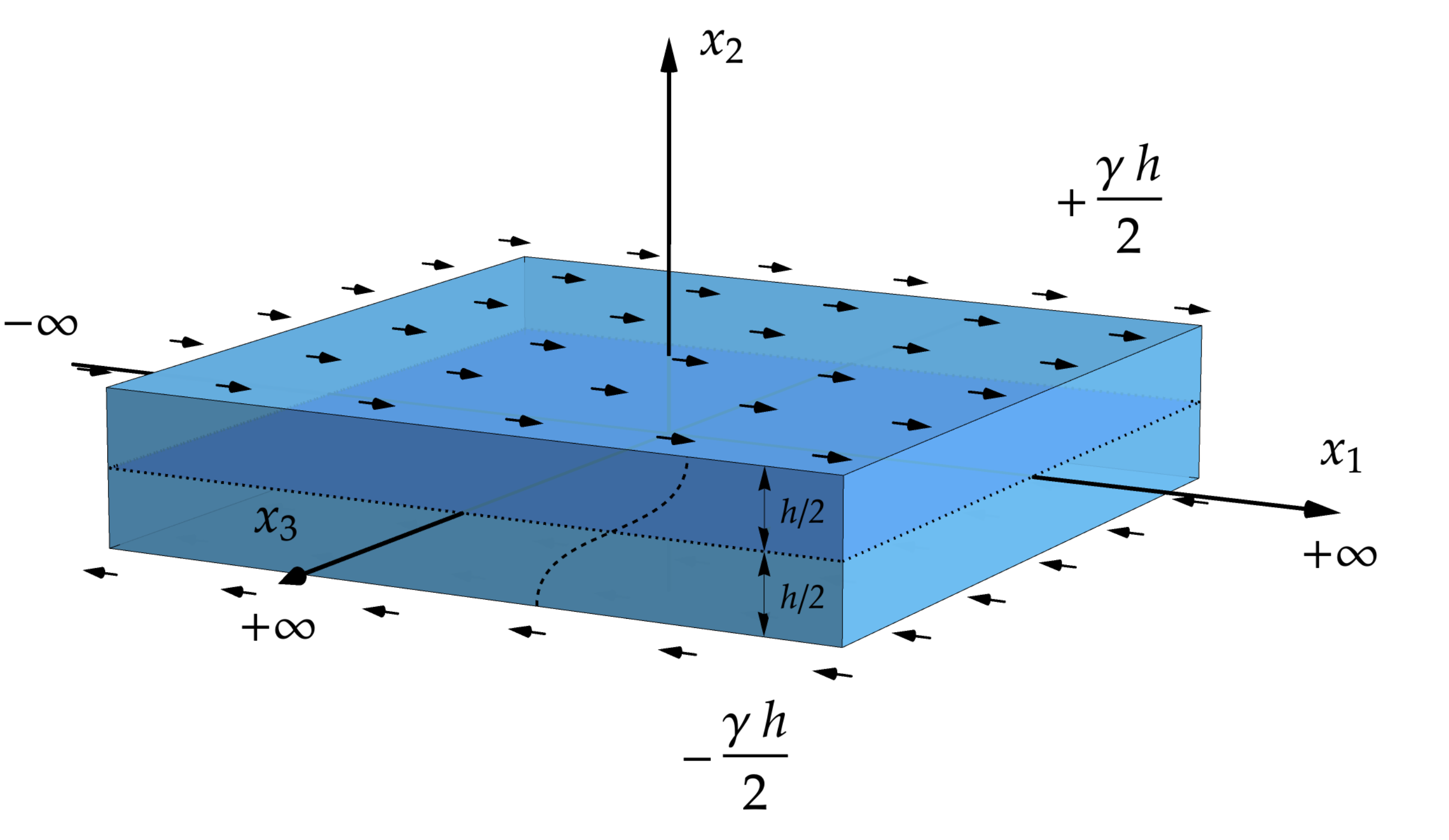}
		\end{subfigure}
		\hfill
		\begin{subfigure}{0.49\textwidth}
			\centering
			\includegraphics[width=\textwidth]{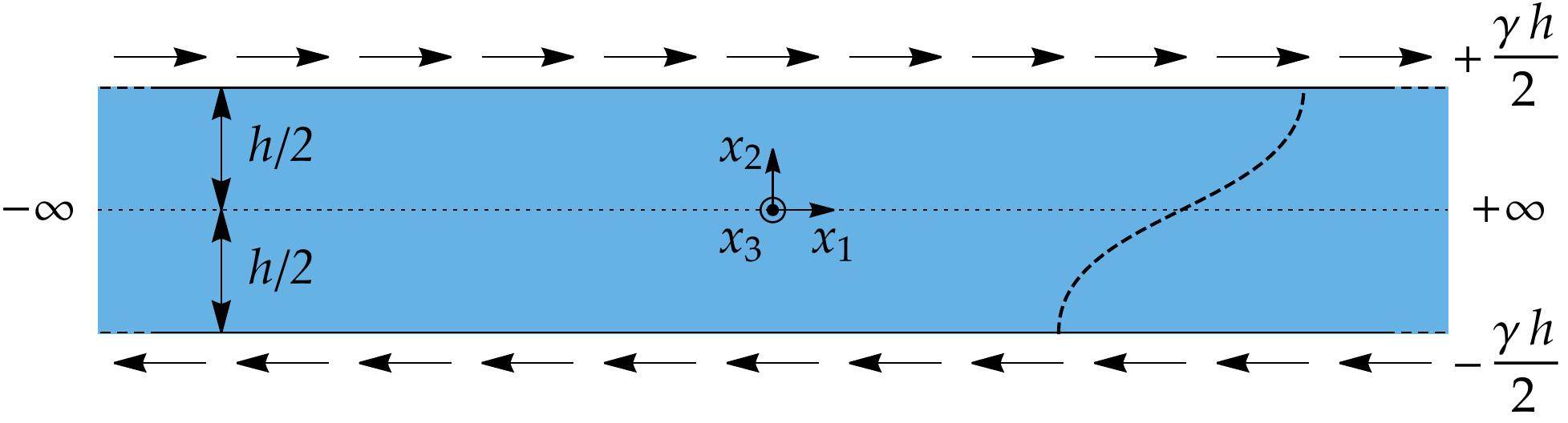}
		\end{subfigure}
		\caption{Sketch of an infinite stripe in the $x_1-$ and  $x_3-$direction of height $h$ subjected to shear boundary conditions.}
		\label{fig:shear}
	\end{figure}
	According to the reference system shown in Fig.~\ref{fig:shear}, the ansatz for the displacement field and the classical micromorphic model is
	\begin{align}
		\label{eq:ansatz_shear_MM}
		\boldsymbol{u}(x_2) &=
		\left(
		\begin{array}{c}
			u_{1}(x_{2}) \\
			0 \\
			0 
		\end{array}
		\right)
		\, ,
		\qquad\qquad
		\boldsymbol{P}(x_2) =
		\left(
		\begin{array}{ccc}
			0 & P_{12}(x_2) & 0 \\
			P_{21}(x_2) & 0 & 0 \\
			0 & 0 & 0  \\
		\end{array}
		\right)
		\, ,
		\\*
		\text{D}\boldsymbol{u}(x_2) &=
		\left(
		\begin{array}{ccc}
			0 & u_{1,2}(x_2) & 0 \\
			0 & 0 & 0 \\
			0 & 0 & 0  \\
		\end{array}
		\right)
		\, .
		\notag
	\end{align}
	For the classical micromorphic model, it is possible to choose between two sets of boundary conditions for the shear problem.
	The first  possible choice is the following full Dirichlet boundary conditions set on $\bou$ and $\boP$,
	\begin{gather}
		\left. u_{1}\right|_{x_{2} = \pm h/2} = \pm \frac{\boldsymbol{\gamma} \, h}{2} \, ,
		\qquad\qquad\qquad
		\left.\pP\right|_{x_{2} = \pm h/2}  = 0 \, ,
		\label{eq:BC_shear_MM}
	\end{gather}
	while the second possible choice is the following mixed boundary conditions set
	\begin{gather}
		\left.u_{1}\right|_{x_{2} = \pm h/2} = \pm \frac{\boldsymbol{\gamma} \, h}{2} \, ,
		\qquad\qquad
		\begin{cases}
			\left.
			\text{D}\boldsymbol{u}(x_2) \times \boldsymbol{\nu}
			\right|_{x_{2} = \pm h/2}
			=
			\left.
			\boldsymbol{P}(x_2) \times \boldsymbol{\nu}
			\right|_{x_{2} = \pm h/2}
			\, ,
			\\*
			\left.
			\left(\boldsymbol{\mathfrak{m}}(x_2) \!\cdot\!\boldsymbol{\nu}\right) \!\cdot\! \left( \boldsymbol{\nu}\otimes\boldsymbol{\nu} \right)
			\right|_{x_{2} = \pm h/2}
			= 0
			\, ,
		\end{cases}
		\label{eq:BC_mix_shear_MM}
	\end{gather}
	where
	\begin{align}
		\text{D}\boldsymbol{u}\times \boldsymbol{\nu}
		=
		\boldsymbol{P}\times \boldsymbol{\nu} \, ,
		\qquad\qquad\qquad
		\left(
		\begin{array}{ccc}
			0 & 0 & 0 \\
			0 & 0 & 0 \\
			0 & 0 & 0 \\
		\end{array}
		\right)
		=
		\left(
		\begin{array}{ccc}
			0 & 0 & 0 \\
			0 & 0 & P_{21} \\
			0 & 0 & 0 \\
		\end{array}
		\right)
		\, .
	\end{align}
	The boundary conditions on $\boldsymbol{P}$ in eqs.(\ref{eq:BC_mix_shear_MM})
	prescribe \textbf{consistent coupling boundary conditions} for $\boldsymbol{P}$, which fix its behaviour accordingly to the behaviour of $\text{D}\boldsymbol{u}$ along the tangential direction on the boundary (${P_{11}=P_{33}=0}$), and an additional \textbf{double-stress-free condition} along the normal at the boundary.
	\begin{figure}[H]
		\centering
		\includegraphics[height=5.5cm]{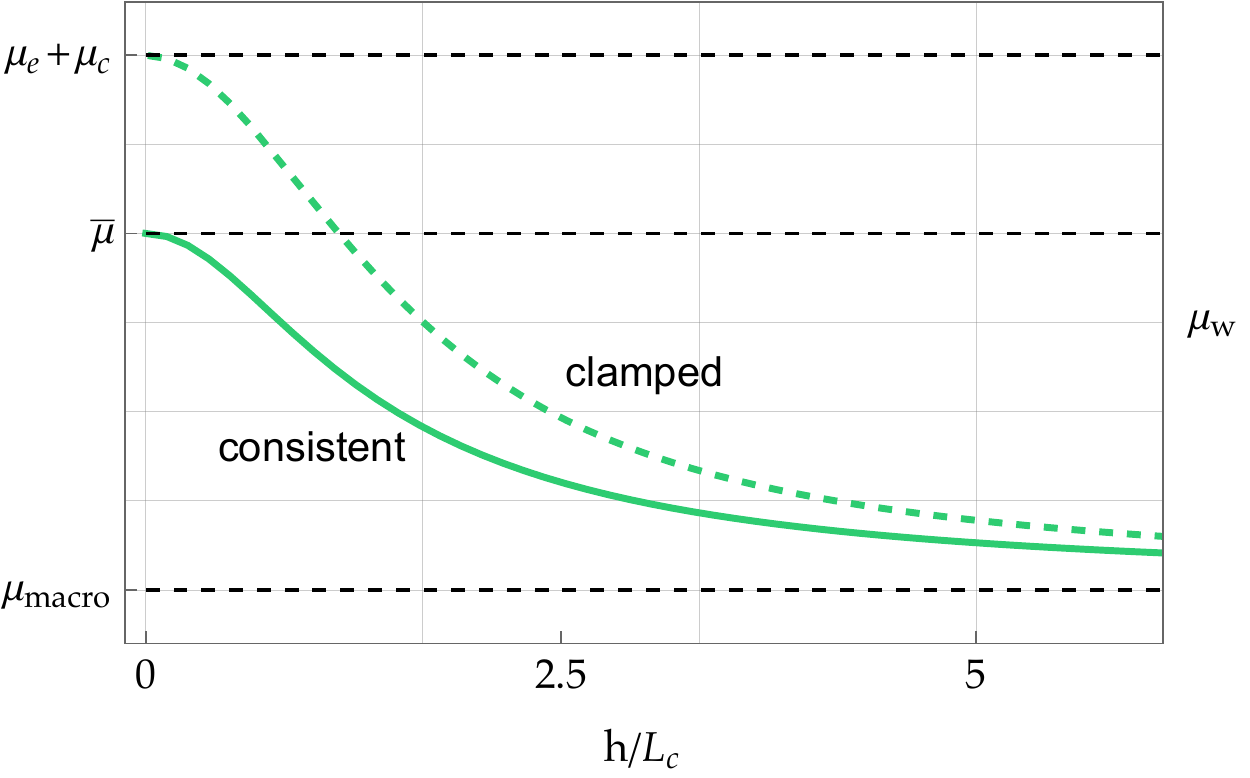}
		\caption{ Comparison of shear stiffness for different Dirichlet-boundary conditions.
			Here $\overline{\mu}=\dfrac{\left(\mu_{\textrm{e}}+\mu_{\textrm{c}}\right)\mu_{\textrm{micro}}}{\mu_{\textrm{e}}+\mu_{\textrm{c}}+\mu_{\textrm{micro}}}$.
		}
		\label{fig:stiff_shear_MM}
	\end{figure}
	\subsubsection{Uniaxial extension}
	
	\begin{figure}[H]
		\centering
		\begin{subfigure}{0.49\textwidth}
			\centering
			\includegraphics[width=\textwidth]{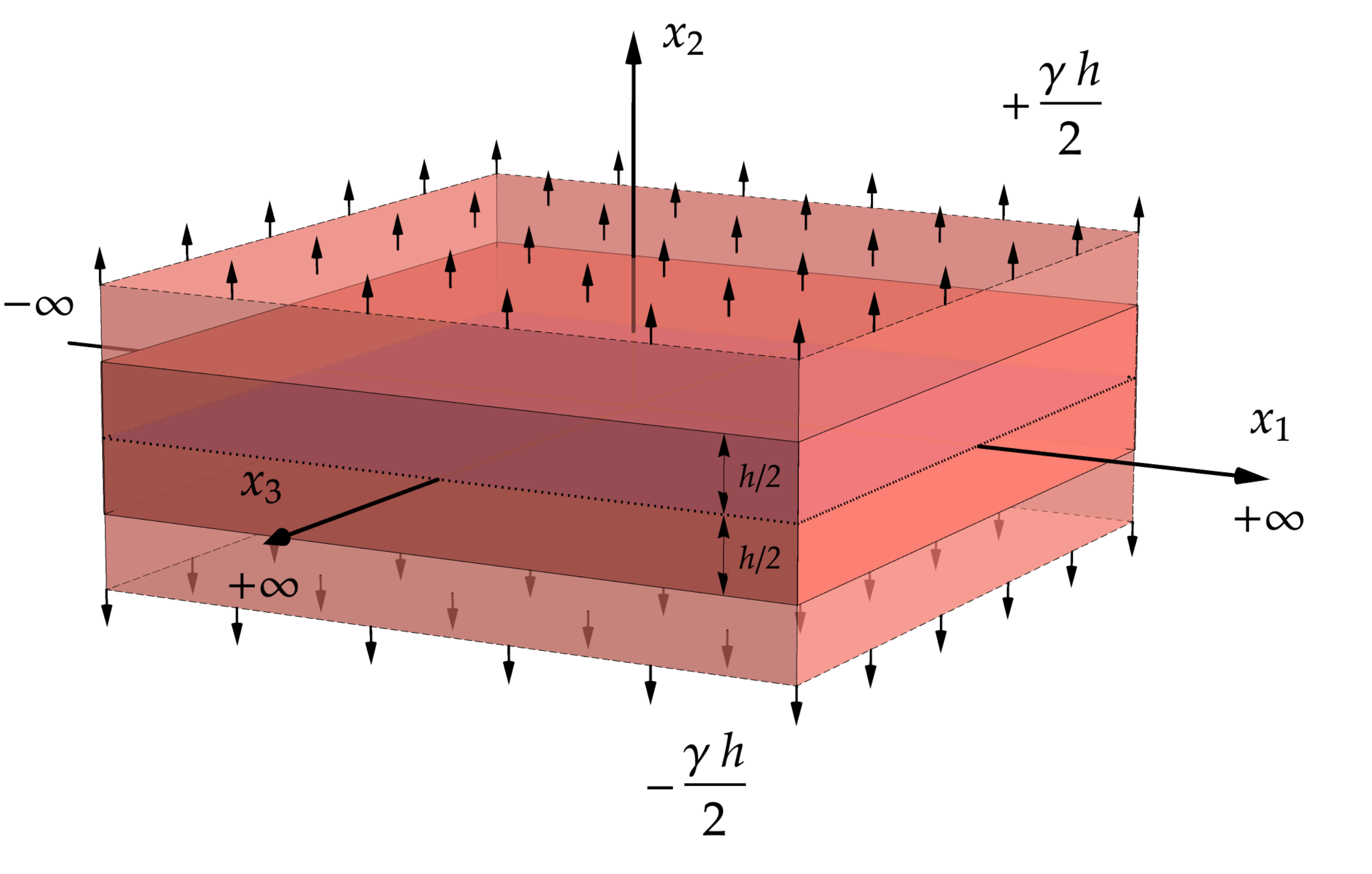}
		\end{subfigure}
		\hfill
		\begin{subfigure}{0.49\textwidth}
			\centering
			\includegraphics[width=\textwidth]{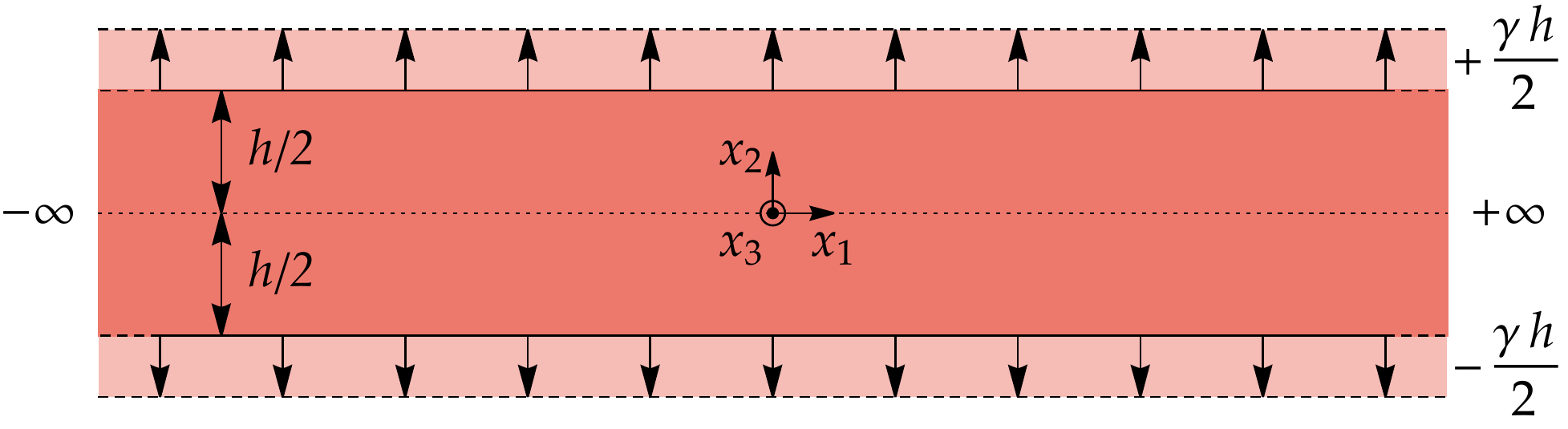}
		\end{subfigure}
		\caption{Sketch of an infinite stripe in the $x_1-$ and  $x_3-$direction of height $h$ subjected to uniaxial extension boundary conditions.}
		\label{fig:extens}
	\end{figure}
	According to the reference system shown in Fig.~\ref{fig:extens}, the ansatz for the displacement field and the classical micromorphic model is
	\begin{align}
		\label{eq:ansatz_extens_MM}
		\pu(x_2) &=
		\left(
		\begin{array}{c}
			0 \\
			u_{2}(x_{2}) \\
			0 
		\end{array}
		\right)
		\, ,
		\qquad\qquad
		\boldsymbol{P}(x_2) =
		\left(
		\begin{array}{ccc}
			P_{11}(x_2) & 0 & 0 \\
			0 & P_{22}(x_2) & 0 \\
			0 & 0 & P_{33}(x_2)  \\
		\end{array}
		\right)
		\, ,
		\\*
		\text{D}\boldsymbol{u}(x_2) &=
		\left(
		\begin{array}{ccc}
			0 & 0 & 0 \\
			0 & u_{2,2}(x_2) & 0 \\
			0 & 0 & 0  \\
		\end{array}
		\right)
		\, .
		\notag
	\end{align}
	For the classical micromorphic and micro-strain models, it is possible to choose between two sets of boundary conditions for the uniaxial extension problem.
	The first possible choice is the following full Dirichlet boundary conditions set on $\bou$ and $\boP$
	\begin{gather}
		\left. u_{2}\right|_{x_{2} = \pm h/2} = \pm \frac{\boldsymbol{\gamma} \, h}{2} \, ,
		\qquad\qquad\qquad
		\left. \pP\right|_{x_{2} = \pm h/2}  = 0 \, ,
		\label{eq:BC_extens_MM}
	\end{gather}
	while the second possible choice is the following mixed boundary conditions set
	\begin{gather}
		\left. u_{1}\right|_{x_{2} = \pm h/2} = \pm \frac{\boldsymbol{\gamma} \, h}{2} \, ,
		\qquad\qquad
		\begin{cases}
			\left.
			\boldsymbol{P}(x_2) \times \boldsymbol{\nu}
			\right|_{x_{2} = \pm h/2}
			=
			\left.
			\text{D}\boldsymbol{u}(x_2) \times \boldsymbol{\nu}
			\right|_{x_{2} = \pm h/2}
			\, ,
			\\*
			\left.
			\left(\boldsymbol{\mathfrak{m}}(x_2) \!\cdot\!\boldsymbol{\nu}\right) \!\cdot\! \left( \boldsymbol{\nu}\otimes\boldsymbol{\nu} \right)
			\right|_{x_{2} = \pm h/2}
			= 0
			\, ,
		\end{cases}
		\label{eq:BC_mix_extens_MM}
	\end{gather}
	where 
	\begin{align}
		\text{D}\boldsymbol{u}\times \boldsymbol{\nu}
		=
		\boldsymbol{P}\times \boldsymbol{\nu} \, ,
		\qquad\qquad\qquad
		\left(
		\begin{array}{ccc}
			0 & 0 & 0 \\
			0 & 0 & 0 \\
			0 & 0 & 0 \\
		\end{array}
		\right)
		=
		\left(
		\begin{array}{ccc}
			0 & 0 & P_{11} \\
			0 & 0 & 0 \\
			-P_{33} & 0 & 0 \\
		\end{array}
		\right)
		\, .
	\end{align}
	The boundary conditions on $\boldsymbol{P}$ in eqs.(\ref{eq:BC_mix_extens_MM})
	prescribe \textbf{consistent coupling boundary conditions} for $\boldsymbol{P}$, which fix its behaviour accordingly to the behaviour of $\text{D}\boldsymbol{u}$ along the tangential direction on the boundary (${P_{11}=P_{33}=0}$), and an additional \textbf{double-stress-free condition} along the normal at the boundary.
	\begin{figure}[H]
		\centering
		\includegraphics[height=5.5cm]{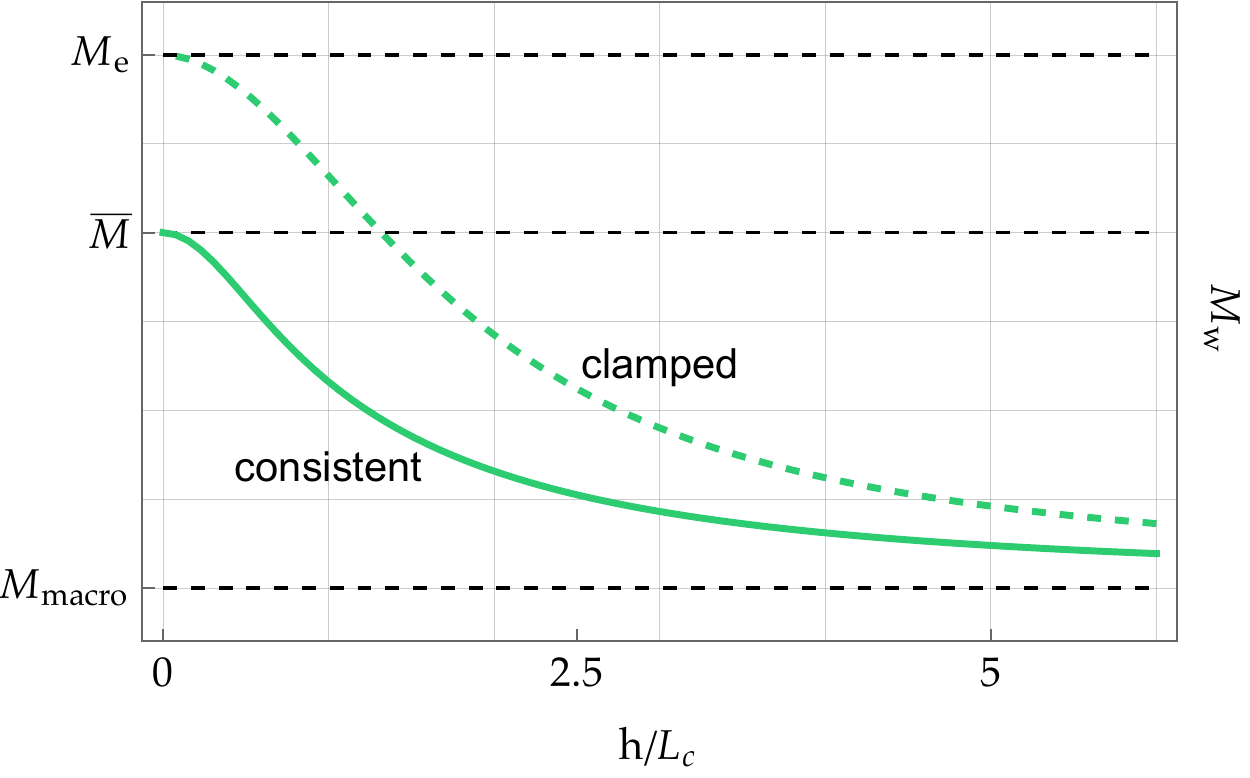}
		\caption{ Comparison of longitudinal stiffness for different boundary conditions.
			Here $\overline{M}=\dfrac{M_{\text{e}}M_{\text{micro}}}{M_{\text{e}}+M_{\text{micro}}}$.
		}
		\label{fig:stiff_extens_MM}
	\end{figure}
	\subsection{Relaxed micromorphic model}
	The  general expression of the strain energy for the isotropic relaxed micromorphic continuum is \cite{rizzi2021bending,rizzi2021torsion,rizzi2021extension} 
	\begin{align}
		W \left(\boldsymbol{\text{D}u}, \boldsymbol{P},\text{Curl}\,\boldsymbol{P}\right)
		=
		&
		\, \mu_{\text{e}} \left\lVert \text{sym} \left(\boldsymbol{\text{D}u} - \boldsymbol{P} \right) \right\rVert^{2}
		+ \mu_{\text{c}} \left\lVert \text{skew} \left(\boldsymbol{\text{D}u} - \boldsymbol{P} \right) \right\rVert^{2}
		+ \frac{\lambda_{\text{e}}}{2} \text{tr}^2 \left(\boldsymbol{\text{D}u} - \boldsymbol{P} \right) 
		\notag
		\\*
		&
		+ \mu_{\text{micro}} \left\lVert \text{sym}\,\boldsymbol{P} \right\rVert^{2}
		+ \frac{\lambda_{\text{micro}}}{2} \text{tr}^2 \left(\boldsymbol{P} \right)
		\label{eq:energy_RM}
		\\*
		&
		+ \frac{\mu \,L_{\text{c}}^2 }{2} \,
		\left(
		a_1 \, \left\lVert \text{dev sym} \, \text{Curl} \, \boldsymbol{P}\right\rVert^2 +
		a_2 \, \left\lVert \text{skew} \,  \text{Curl} \, \boldsymbol{P}\right\rVert^2 +
		\frac{a_3}{3} \, \text{tr}^2 \left(\text{Curl} \, \boldsymbol{P}\right)
		\right)
		\notag
		\, ,
		\notag
	\end{align}
	where ($\mu_{\text{e}}$,$\lambda_{\text{e}}$), ($\mu_{\text{micro}}$,$\lambda_{\text{micro}}$), $\mu_{\text{c}}$, $L_{\text{c}} > 0$, and ($a_1$,$a_2$,$a_3$) are the parameters related to the meso-scale, the parameters related to the micro-scale, the Cosserat couple modulus, the characteristic length, and the three general isotropic curvature parameters, respectively.
	This energy expression represents the  most general isotropic form possible for the relaxed micromorphic model.
	In the absence of body forces, the equilibrium equations  are then
		\begin{equation}
		\text{Div}\,\bosigma=\boldsymbol{0}\,,
		\qquad\qquad
		\bosigma
		- 2 \mu_{\text{micro}}\,\text{sym}\,\pP
		- \lambda_{\text{micro}} \text{tr} (\pP) \id
		- \Curl\bom
		=\boldsymbol{0}\,,
		\end{equation}
		where
		\begin{align}
		\boldsymbol{\sigma}
		& 
		\coloneqq 2\mue\,\sym\!(\pD\pu - \pP) 
		+ 2\mu_{\text{c}}\,\skew\!(\pD\pu - \pP)
		+ \le \text{tr}(\pD\pu - \pP) \id\,
		\end{align}
		is the second order force stress tensor and
		\begin{align}
		\bom
		&
		\coloneqq 
		\mu L_{\text{c}}^{2}
		\Big(
		a_1 \, \text{dev sym} \, \Curl\pP +
		a_2 \, \skew \Curl\pP +
		a_3 \, \text{tr} (\Curl\pP)\id
		\Big)
		\,.
		\notag
		\end{align}
	    is the second order moment tensor in the relaxed micromorphic model.
	\subsubsection{Simple shear}
	According to the reference system shown in Fig.~\ref{fig:shear}, the ansatz for the displacement field and the relaxed micromorphic model is \cite{rizzi2021shear}
	\begin{align}
		\label{eq:non_zero_compo}
		\boldsymbol{u}(x_2) &=
		\left(
		\begin{array}{c}
			u_{1}(x_{2}) \\
			0 \\
			0 
		\end{array}
		\right)
		\, ,
		\qquad\qquad
		\boldsymbol{P}(x_2) =
		\left(
		\begin{array}{ccc}
			0 & P_{12}(x_2) & 0 \\
			P_{21}(x_2) & 0 & 0 \\
			0 & 0 & 0  \\
		\end{array}
		\right)
		\, ,
		\\*
		\text{D}\boldsymbol{u}(x_2) &=
		\left(
		\begin{array}{ccc}
			0 & u_{1,2}(x_2) & 0 \\
			0 & 0 & 0 \\
			0 & 0 & 0  \\
		\end{array}
		\right)
		\, .
		\notag
	\end{align}
	The boundary conditions for the simple shear in the relaxed micromorphic model are the following
	\begin{equation}
		\left. u_{1}\right|_{x_{2} = \pm h/2} = \pm \frac{\boldsymbol{ \gamma} \, h}{2} \, ,
		\qquad\qquad
		\left. \pP\right|_{x_{2} = \pm h/2} \times \boldsymbol{\nu}
		=
		\left. \text{D}\boldsymbol{u}\right|_{x_{2} = \pm h/2} \times \boldsymbol{\nu}
		\, ,
		\label{eq:BC_RM_shear}
	\end{equation}
	where $\boldsymbol{\nu}$ is the normal unit vector on the upper and lower surface.
	\begin{figure}[H]
		\centering
		\includegraphics[height=5.5cm]{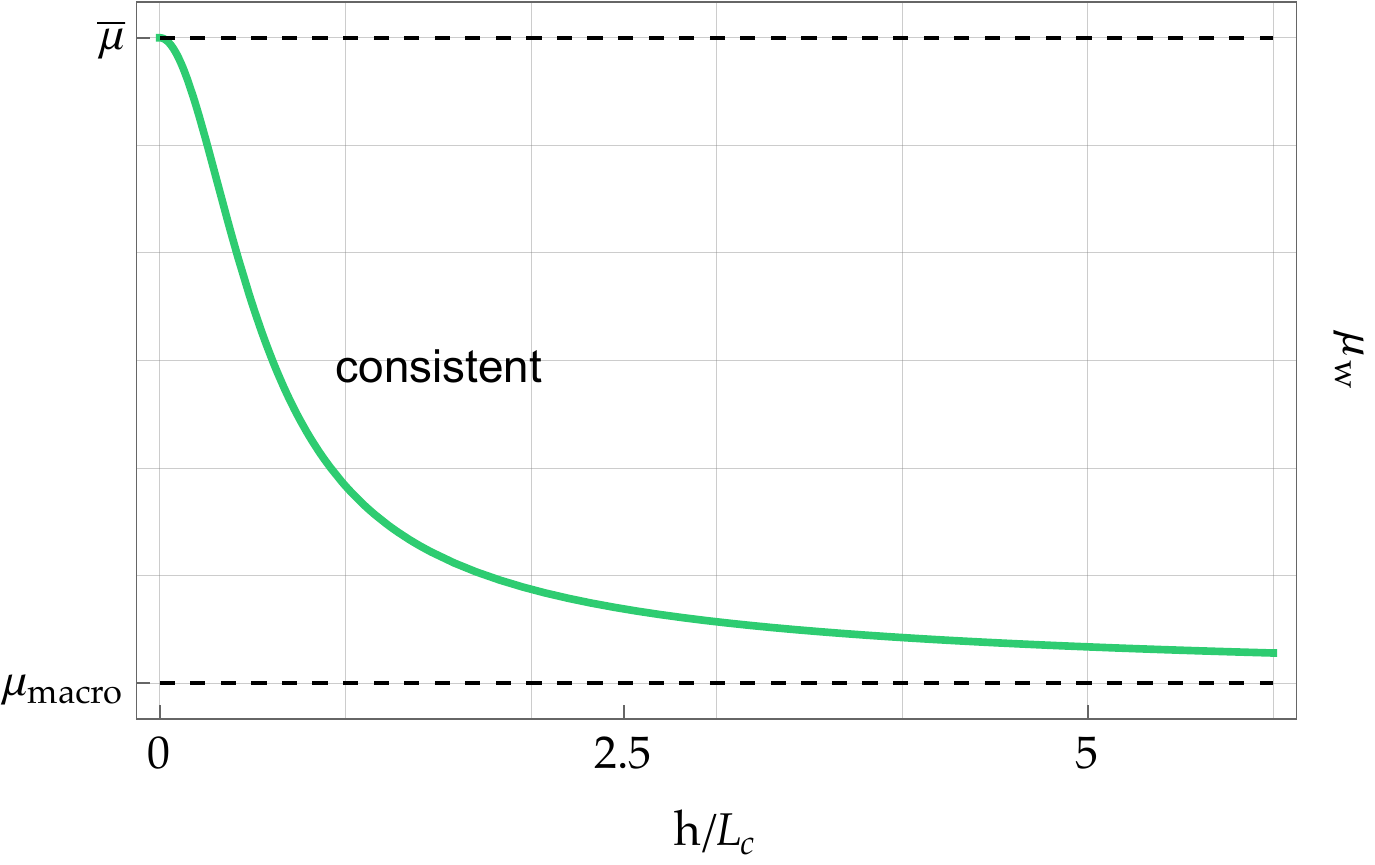}
		\caption{
			Longitudinal stiffness for consistent coupling in the relaxed micromorphic model.
			The stiffness for $L_{\text{c}}\to\infty$ is bounded and its value $\overline{\mu}=\dfrac{\left(\mu_{\textrm{e}}+\mu_{\textrm{c}}\right)\mu_{\textrm{micro}}}{\mu_{\textrm{e}}+\mu_{\textrm{c}}+\mu_{\textrm{micro}}}$ coincides with the one obtained with the consistent coupling conditions in the full micromorphic model.}
		\label{fig:stiff_shear_RM}
	\end{figure}
	
	It is underlined that only the two limits for $0\leftarrow L_{\text{c}}\rightarrow\infty$ coincide between the relaxed micromorphic model and the classical micromorphic model with consistent coupling boundary conditions, not the full solution.
	\subsubsection{Uniaxial extension}
	According to the reference system shown in Fig.~\ref{fig:extens}, the ansatz for the displacement field in the relaxed micromorphic model is
	\begin{align}
		\label{eq:ansatz_RM}
		\boldsymbol{u}(x_2) &=
		\left(
		\begin{array}{c}
			0 \\
			u_{2}(x_{2}) \\
			0 
		\end{array}
		\right)
		\, ,
		\qquad\qquad
		\boldsymbol{P}(x_2) =
		\left(
		\begin{array}{ccc}
			P_{11}(x_2) & 0 & 0 \\
			0 & P_{22}(x_2) & 0 \\
			0 & 0 & P_{33}(x_2)  \\
		\end{array}
		\right)
		\, ,
		\\*
		\text{D}\boldsymbol{u}(x_2) &=
		\left(
		\begin{array}{ccc}
			0 & 0 & 0 \\
			0 & u_{2,2}(x_2) & 0 \\
			0 & 0 & 0  \\
		\end{array}
		\right)
		\, .
		\notag
	\end{align}
	The boundary conditions for the uniaxial extension in the relaxed micromorphic model are
	\begin{align}
		\left. u_{2}\right|_{x_{2} = \pm h/2} = \pm \frac{\boldsymbol{\gamma} \, h}{2}
		\, , 
		\qquad\qquad
		\left. \pP\right|_{x_{2} = \pm h/2} \times \boldsymbol{\nu}
		=
		\left. \text{D}\pu\right| _{x_{2} = \pm h/2} \times \boldsymbol{\nu}
		\, .
		\label{eq:BC_RM_extension}
	\end{align}
	where $\boldsymbol{\nu}$ is the normal unit vector to the upper and lower surface.
	\begin{figure}[H]
		\centering
		\includegraphics[height=5.5cm]{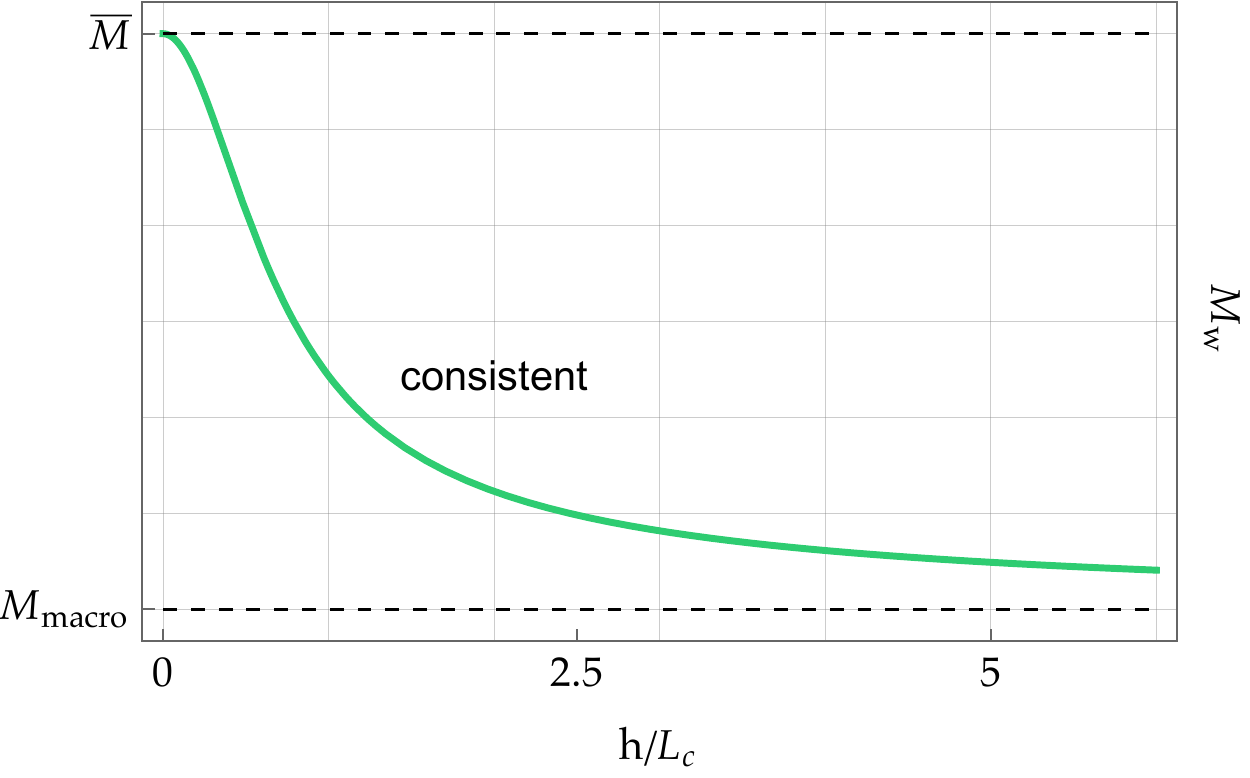}
		\caption{ Longitudinal stiffness for consistent coupling in the relaxed micromorphic model.
			The stiffness for $L_{\text{c}}\to\infty$ is bounded and its value $\overline{M}=\dfrac{M_{\textrm{e}}\, M_{\textrm{micro}}}{M_{\textrm{e}}+M_{\textrm{micro}}}$ coincides with the one obtained with the consistent coupling conditions in the full micromorphic model.}
		\label{fig:stiff_extens_RM}
	\end{figure}
	
	Again, it is underlined that only the two limits for $0\leftarrow L_{\text{c}}\rightarrow\infty$ coincide between the relaxed micromorphic model and the classical micromorphic model with consistent coupling boundary conditions, not the full solution.
	\subsection{Second gradient model}
	The strain energy density for the isotropic second gradient model with simplified curvature \cite{mindlin1964micro,dellisola2009generalized,rizzi2021shear,rizzi2021torsion,altenbach2019higher,rizzi2021extension,altenbach2019higher} is
	\begin{align}
		W \left(\boldsymbol{\text{D}u}, \boldsymbol{\text{D}^2 u}\right)
		= &
		\, \mu_{\text{macro}} \left\lVert \text{sym}\,\boldsymbol{\text{D}u} \right\rVert^{2}
		+ \frac{\lambda_{\text{macro}}}{2} \text{tr}^2 \left(\boldsymbol{\text{D}u} \right)
		\label{eq:energy_Strain_Grad}
		\\*
		&
		+ \frac{\mu \, L_{\text{c}}^2}{2}
		\left(
		a_1 \, \left\lVert \text{D} \Big(\text{dev} \, \text{sym} \, \boldsymbol{\text{D} u}\Big) \right\rVert^2
		+ a_2 \, \left\lVert \text{D} \Big(\text{skew} \, \boldsymbol{\text{D} u}\Big) \right\rVert^2
		+  \frac{2}{9} \, a_3 \, \left\lVert \text{D}
		\Big( 
		\text{tr} \left(\boldsymbol{\text{D} u}\right) \, \boldsymbol{\mathbbm{1}}
		\Big) \right\rVert^2
		\right)
		\, ,
		\notag
	\end{align}
	while the equilibrium equations without body forces are the following:
	\begin{align}
		\text{Div}
		\bigg[
		2 \mu_{\text{macro}} \,\text{sym}\,\boldsymbol{\text{D}u}
		+ \lambda_{\text{macro}} \text{tr} \left(\boldsymbol{\text{D}u}\right) \boldsymbol{\mathbbm{1}}
		\hspace{8cm}
		\label{eq:equi_Strain_Grad}\nonumber
		\\*
		- 
		\text{DIV}
		\bigg[
		\underbrace{
		\mu L_{\text{c}}^{2}
		\left( 
			a_1 \, \text{D}\left(\text{dev} \, \text{sym} \left(\boldsymbol{\text{D}u}\right)\right)
			+ a_2 \, \text{D}\left(\text{skew}  \left(\boldsymbol{\text{D}u}\right)\right)
			+ \frac{2}{9} \, a_3 \, \text{D}\left(\text{tr} \left(\boldsymbol{\text{D}u}\right)\boldsymbol{\mathbbm{1}}\right)
		    \right) }_{\mathlarger{=:\gbm\in\bR^{3\times3\times3}}}
		\bigg]
		\bigg]
		=
		\\*
		\text{Div}\bigg[
		2 \mu_{\text{macro}} \,\text{sym}\,\boldsymbol{\text{D}u}
		+ \lambda_{\text{macro}} \text{tr} \left(\boldsymbol{\text{D}u}\right) \boldsymbol{\mathbbm{1}}
		\hspace{7.2cm}\nonumber
		\\*
		\quad- \mu L_{\text{c}}^{2} \,
		\left(
		a_1 \, \text{dev} \, \text{sym} \, \boldsymbol{\Delta} \left(\boldsymbol{\text{D}u}\right)
		+ a_2 \, \text{skew} \, \boldsymbol{\Delta} \left(\boldsymbol{\text{D}u}\right)
		+ \frac{2}{9} \, a_3 \, \text{tr} \left(\boldsymbol{\Delta} \left(\boldsymbol{\text{D}u}\right)\right)\boldsymbol{\mathbbm{1}}
		\right)
		\bigg]
		= \boldsymbol{0} \, .
		\notag
	\end{align}
	\subsubsection{Simple shear}
	According to the reference system shown in Fig.~\ref{fig:shear}, the ansatz for the displacement field is \cite{rizzi2019identificationII,rizzi2021shear}
	\begin{align}
		\label{eq:ansatz_shear_SG}
		\boldsymbol{u}(x_2) =
		\left(
		\begin{array}{c}
			u_{1}(x_{2}) \\
			0 \\
			0 
		\end{array}
		\right)
		\, ,
		\qquad\qquad\qquad
		\text{D}\boldsymbol{u}(x_2) =
		\left(
		\begin{array}{ccc}
			0 & u_{1,2}(x_2) & 0 \\
			0 & 0 & 0 \\
			0 & 0 & 0  \\
		\end{array}
		\right)
		\, .
	\end{align}
	For the second gradient model, it is possible to choose between two sets of boundary conditions for the shear problem.
	The first  possible choice is the following full Dirichlet boundary conditions set,
	\begin{gather}
		\left. u_{1}\right|_{x_{2} = \pm h/2} = \pm \frac{\boldsymbol{\gamma} \, h}{2} \, ,
		\qquad\qquad\qquad
		\left. u'_{1}\right|_{x_{2} = \pm h/2} = 0
		\quad
		\Leftrightarrow
		\quad
		\text{D}\bou = 0
		\, ,
		\label{eq:BC_shear_SG}
	\end{gather}
	while the second possible choice is the following mixed boundary conditions set
	\begin{gather}
		\left. u_{1}\right|_{x_{2} = \pm h/2} = \pm \frac{\boldsymbol{\gamma} \, h}{2} \, ,
		\qquad\qquad
		(\left. \boldsymbol{\mathfrak{m}}\right|_{x_{2} = \pm h/2} \!\cdot\! \boldsymbol{\nu}) \!\cdot\!\pnu
		= 0
		\, ,
		\label{eq:BC_mix_shear_SG}
	\end{gather}
	\begin{figure}[H]
		\centering
		\includegraphics[height=5.5cm]{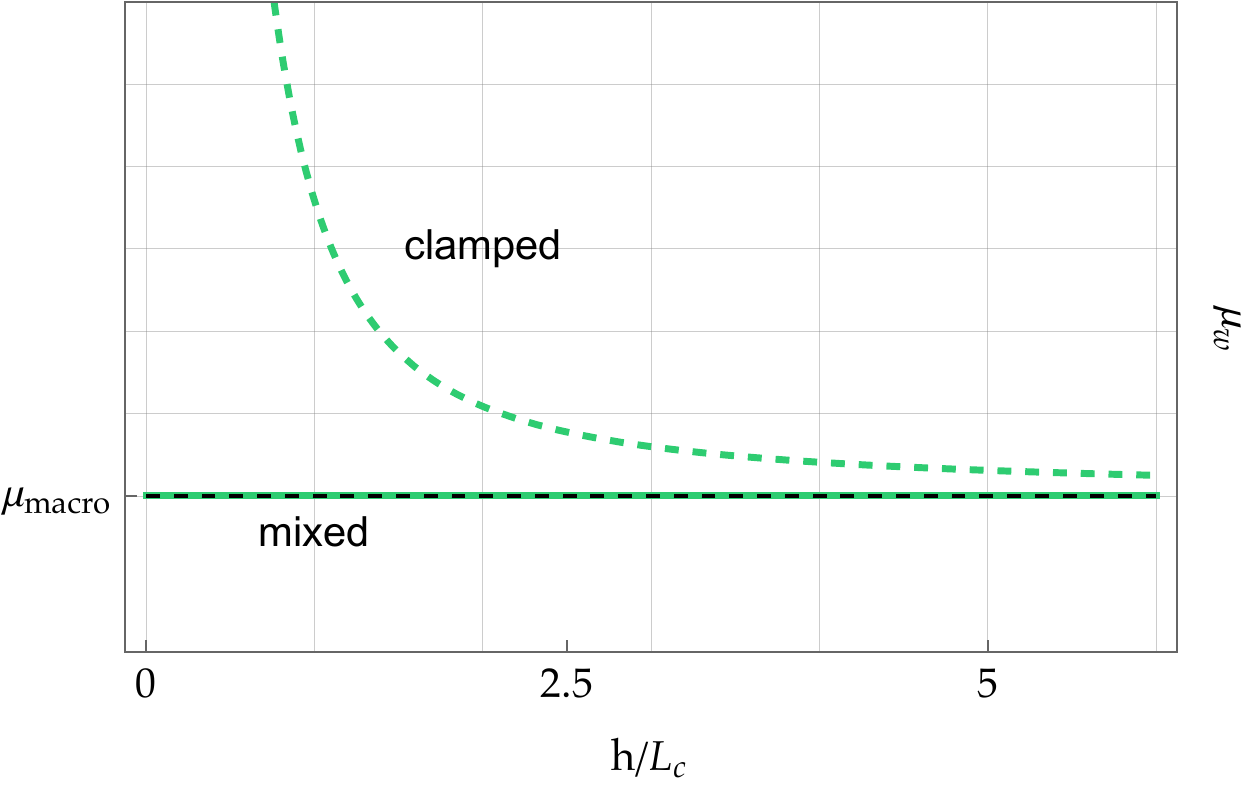}
		\caption{
			Comparison of shear stiffness for different boundary conditions.
			The fully clamped solution shows a stiffness singularity while the mixed boundary condition yields the homogenous solution with stiffness $\mu_{\text{macro}}$.
		}
		\label{fig:stiff_shear_SG}
	\end{figure}
	\subsubsection{Uniaxial extension}
	According to the reference system shown in Fig.~\ref{fig:extens}, the ansatz for the displacement field and the classical micromorphic model is \cite{rizzi2019identificationII}
	\begin{align}
		\label{eq:ansatz_extens_SG}
		\boldsymbol{u}(x_2) =
		\left(
		\begin{array}{c}
			0 \\
			u_{2}(x_{2}) \\
			0 
		\end{array}
		\right)
		\, ,
		\qquad\qquad\qquad
		\text{D}\boldsymbol{u}(x_2) =
		\left(
		\begin{array}{ccc}
			0 & 0 & 0 \\
			0 & u_{2,2}(x_2) & 0 \\
			0 & 0 & 0  \\
		\end{array}
		\right)
		\, .
	\end{align}
	For the second gradient model, it is possible to choose between two sets of boundary conditions for the uniaxial extension problem.
	The first possible choice is the following full Dirichlet boundary conditions set,
	\begin{gather}
		\left. u_{2}\right|_{x_{2} = \pm h/2} = \pm \frac{\boldsymbol{\gamma} \, h}{2} \, ,
		\qquad\qquad\qquad
		\left. u'_{2}\right|_{x_{2} = \pm h/2} = 0
		\quad
		\Leftrightarrow
		\quad
		\text{D}\bou = 0
		\, ,
		\label{eq:BC_extens_SG}
	\end{gather}
	while the second possible choice is the following mixed boundary conditions set
	\begin{gather}
		\left. u_{2}\right|_{x_{2} = \pm h/2} = \pm \frac{\boldsymbol{\gamma} \, h}{2} \, ,
		\qquad\qquad
		(\left. \boldsymbol{\mathfrak{m}}\right|_{x_{2} = \pm h/2} \!\cdot\! \boldsymbol{\nu}) \!\cdot\!\pnu
		= 0
		\, ,
		\label{eq:BC_mix_extens_SG}
	\end{gather}
	\begin{figure}[H]
		\centering
		\includegraphics[height=5.5cm]{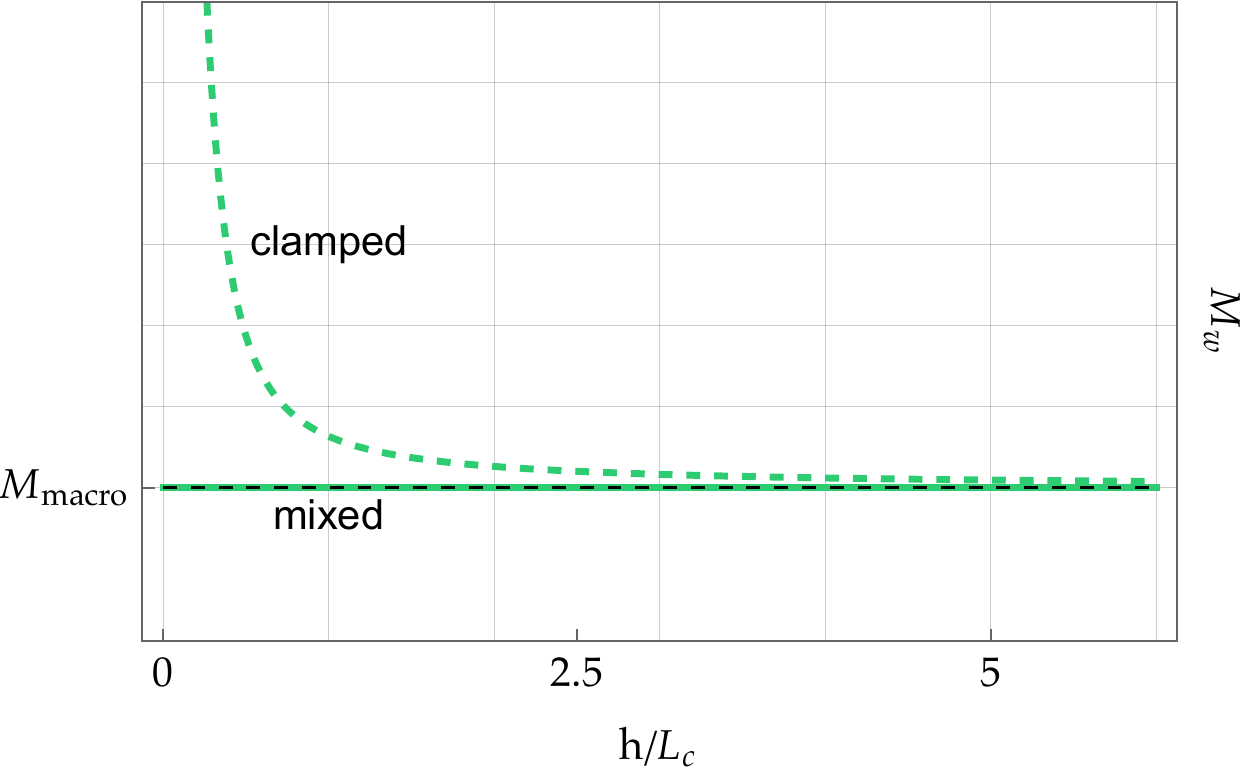}
		\caption{
			Comparison of longitudinal stiffness for different boundary conditions.
			The fully clamped solution shows a stiffness singularity while the mixed boundary condition yields the homogenous solution with stiffness $M_{\text{macro}}$.
		}
		\label{fig:stiff_extens_SG}
	\end{figure}
	
	\section{Non-redundant formulation}
	
	In this section we want to show that the classical micromorphic model with zero Cosserat couple modulus $\muc=0$ but full curvature expression is not
	redundant according to the following definition introduced by Romano et al. \cite{romano2016micromorphic}:
	\begin{defn}\label{def: Romano def}
		A formal differential operator is said to meet a \textit{constraint qualification} if no other implicit description can be constructed by extracting a strictly more economical condition from the given operator. An operator for classical micromorphic continua which does not fulfill the constraint qualification requirement is termed \textit{redundant}.
	\end{defn}
	
	\noindent According to this definition, the classical micromorphic model with positive Cosserat couple modulus $\muc> 0$ is redundant because if a subset of the terms in the definition of the elastic energy is zero, in this case  $\D\!\pu-\pP = 0$ and $\sym \pP = 0$, then this implies that all the remaining terms are zero too, in this case $\D\!\pP=0$ (see \eqref{eqt:numerical_scheme1}). 
	
	
	
	In order to provide a more precise definition, we specialize definition \ref{def: Romano def} to the relation between the first order strain (involving only $\pD\pu$ and $\pP$), and the curvature (involving only $\pD\pP$).
	
	\begin{defn}
		\label{def: ours}
		Let us consider a suitable linear micromorphic model and an elastic energy density $W$ involving the first order derivatives of the displacement $\pu$ and the micro-distortion field $\pP$, which splits into two parts: $W(\pD\pu,\pP,\pD\pP)=W_1(\pD\pu,\pP)+W_2(\pD\pP)$.
		The model is said \textit{non redundant} if none of the two following possibilities occurs: 
		\begin{align}
                W_{1}(\pD\pu,\pP)=0\quad\Longrightarrow\quad W(\pD\pu,\pP,\pD\pP)&=0, &\textrm{or}
                \\
                W_{2}(\pD\pP)=0\quad\Longrightarrow\quad W(\pD\pu,\pP,\pD\pP)&=0. 
        \end{align}
	\end{defn}
	
	\noindent Also according to this second definition we have that the classical micromorphic model with positive Cosserat couple modulus $\muc> 0$ is redundant and the relaxed micromorphic model with $\muc=0$ is non redundant. Furthermore, we avoid all the indeterminacies related to definition \ref{def: Romano def}. 
	
	In the next section, we will apply the new definition to a family of generalized continua to establish which of them are redundant or non redundant.
	\subsection{A view on rigid body deformations and zero energy modes in generalized continua}
	\label{sec:rigid_body_motion}
\begin{figure}[H]
   	\begin{centering}
   		\begin{tabular}{ lc }
   			a) classical micromorphic
   			& 
   			\begin{tabular}{c}
   				\includegraphics[scale=1.3]{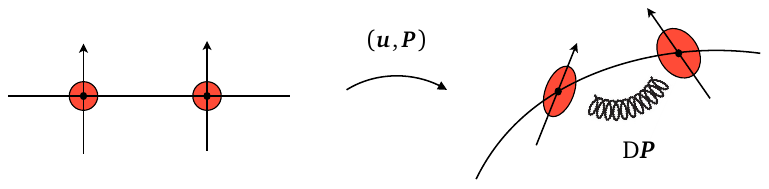} 
   			\end{tabular} 
   			 \tabularnewline
   	        b) relaxed micromorphic
   	        & 
   	        \begin{tabular}{c}
   	    	    \includegraphics[scale=1.3]{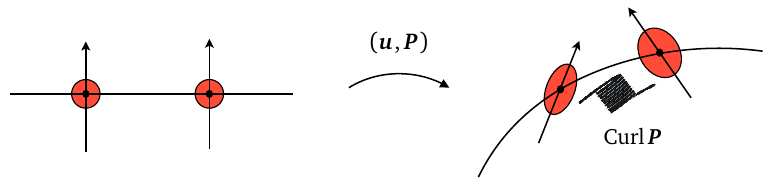} 
   	        \end{tabular} 
   	        \tabularnewline
   		    c) Cosserat
   		    & 
   		    \begin{tabular}{c}
   		   	    \includegraphics[scale=1.32]{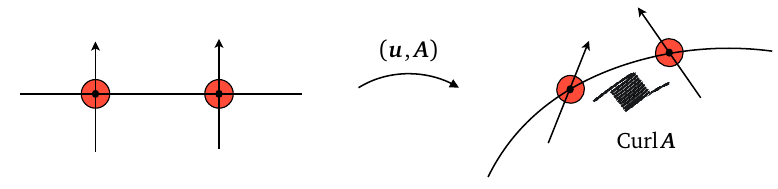} 
   		    \end{tabular} 
   		    \tabularnewline
   		\end{tabular}
   		\par\end{centering}
    	\caption{
    		a,b)  the difference in local rotations inherent in $\pP$ are controlled by the curvature term (classical micromorphic via $\text{D}\pP$, relaxed micromorphic via $\Curl\pP$). c) In the Cosserat model the material particles can only rotate against each other and they are connected by a ``rotational spring'', naturally given by $\Curl\boA$ see \eqref{Nye}.
    	}
    	\label{fig:Pcurvature}
    \end{figure}

	In classical isotropic linear elasticity the notion of infinitesimal rigid transformations and zero elastic energy modes is equivalent (see Fig.~\ref{fig:Pcurvature}). Indeed, if 
	\begin{align}
		W_{\textrm{lin}}\left(\text{D}\pu\right) & =\mu\,\left\Vert \sym\pD\pu\right\Vert ^{2}+\frac{\lambda}{2}\,\tr^{2}(\pD\pu)=\mu\,\left\Vert \dev\sym\pD\pu\right\Vert ^{2}+\frac{\kappa}{2}\,\tr^{2}(\pD\pu),\qquad\mu,\kappa>0
	\end{align}
	then
	\begin{equation*}
		W_{\textrm{lin}}\left(\text{D}\pu\right)=0\quad\Longleftrightarrow\quad\sym\text{D}\pu=0\quad\Longleftrightarrow\quad\pu=\!\underbrace{\overline{\boA}\!\cdot\!\boox+\overline{\bob}}_{\textrm{\ensuremath{\substack{\\
						\textrm{infinitesimal}\\
						\textrm{rigid movement}
					}
		}}},\qquad \overline{\boA}\in\so,\overline{\bob}\in\bR^3.
	\end{equation*}
	Indeed,
	\begin{equation}
		\textrm{dist}^{2}\left(\text{D}\boldsymbol{u},\so\right)\defi\inf_{\boA\in\so}\left\Vert \text{D}\boldsymbol{u}-\boA\right\Vert _{\bR^{3\times3}}^{2}=\left\Vert \sym\text{D}\boldsymbol{u}\right\Vert _{\bR^{3\times3}}^{2}
	\end{equation}
	\begin{equation*}
		\textrm{dist}^{2}\left(\text{D}\boldsymbol{u},\so\right)=0\quad\Longrightarrow\quad\sym\text{D}\boldsymbol{u}=0\quad\Longrightarrow\quad\pu=\overline{\boA}\!\cdot\!\boox+\overline{\bob}.
	\end{equation*}
	However, in generalized continua, the relation between the two notions is more involved and depends critically on specific further constitutive assumptions. Let us focus on the pair $(\pu,\pP)$ of displacement $\pu$ and microdistortion $\pP$ in the micromorphic model and the generic isotropic energy expression
	\begin{subequations}
		\begin{numcases}{W\left(\text{D}\boldsymbol{u},\boldsymbol{P},\text{D}\boldsymbol{P}\right)=}
			\left\Vert \text{D}\boldsymbol{u}-\boldsymbol{P}\right\Vert ^{2}+\left\Vert \sym\boldsymbol{P}\right\Vert ^{2}+\left\Vert \text{D}\boldsymbol{P}\right\Vert ^{2},
			& 
			$\muc>0$ \label{eq:energy21a}
			\\
			\left\Vert \sym\!(\text{D}\pu-\pP)\right\Vert ^{2}+\left\Vert \sym\pP\right\Vert ^{2}+\left\Vert \text{D}\pP\right\Vert ^{2},
			&
			$\muc=0.$ \label{eq:energy21b}
		\end{numcases}
	\end{subequations}
	The case \eqref{eq:energy21b} corresponds to a zero Cosserat couple modulus $\muc\equiv0$. Obviously
	
	\begin{subequations}
		\begin{numcases}{	\hspace{-1cm}\underbrace{W\left(\text{D}\boldsymbol{u},\boldsymbol{P},\text{D}\boldsymbol{P}\right)=0}_{\textrm{zero energy mode}}\quad\Longleftrightarrow\quad}
			\textrm{case \eqref{eq:energy21a}}\quad\pu=\overline{\boA}\!\cdot\!\boox+\overline{\bob},\quad\text{D}\pu=\pP=\overline{\boA}
			&
			$\muc>0$ \label{eq:energy31a}
			\\
			\textrm{case \eqref{eq:energy21b}}\quad\pu=\overline{\boA}\!\cdot\!\boox+\overline{\bob},\quad\text{D}\pu=\overline{\boA},\;\pP=\overline{\overline{\boA}}, \label{eq:energy31b}
			&
			$\muc=0$ 
		\end{numcases}
	\end{subequations}
	where $\overline{\boA},\overline{\overline{\boA}}\in\so$ are two (in general different) constant infinitesimal rotations. 
	Let us see next how to involve Dirichlet boundary conditions to fix the rigid body modes at zero energy.
	\begin{description}
		\item[Linear elasticity:] 
		\begin{equation}
			\begin{tikzcd}[
				column sep=10mm,
				row sep=-3mm,
				every arrow/.append style={line width=0.7pt},
				]
				\sym\text{D}\pu=0
				\arrow[r,Rightarrow]
				&
				\pu=\overline{\boA}\!\cdot\!\boox+\overline{\bob}
				&
				\textrm{and\quad  b.c.}
				&
				\left.\pu\right|_{\Gamma}=0
				\arrow[r,Rightarrow]
				&
				\pu\equiv0
				\, .
			\end{tikzcd}
		\end{equation}
	\end{description}
	Thus, prescribing $\left.\pu\right|_{\Gamma}=0$ fixes the constant infinitesimal rotation $\overline{\boA}$ and the constant translation $\overline{\bob}$.

    \begin{description}	
		\item[Cosserat:] 
		\begin{equation}
		\hspace{-7mm}
		\begin{aligned}
		\begin{tikzcd}[
		column sep=10mm,
		row sep=3mm,
		every arrow/.append style={line width=0.7pt},
		/tikz/column 1/.append style={anchor=base east},
		/tikz/column 2/.append style={anchor=base east},
		/tikz/column 3/.append style={anchor=base west},
		every label/.append style = {font = \normalsize },
		every arrow/.append style={line width=0.7pt},
		]
		\vphantom{pD\pu-\boA=0}
		\arrow[dd, start anchor=north east, end anchor=south east, no head, xshift=0em, decorate, decoration={brace},"{\Longrightarrow}" right=7pt,phantom]
		&
		\sym\pD\pu=0\;\overset{\textrm{Korn}}{\Longrightarrow}\;\pu=\overline{\boA}\!\cdot\!\boox+\overline{\bob}
		\arrow[dd, start anchor=north east, end anchor=south east, no head, xshift=0em, decorate, decoration={brace},"{\underset{\vphantom{\text{b.c.}}}{\overset{\text{b.c.}}{\Longrightarrow}}}" right=7pt]
		&
		\\
		\pD\pu-\boA=0
		&
		\boA=\skew\pD\pu=\overline{\boA} 
		&
		\left.\pu\right|_{\Gamma}=0\;\Rightarrow\;\skew\pD\pu=  0=\overline{\boA}=\boA.
		\\
		\vphantom{pD\pu-\boA=0}
		&
		\overline{\boA}\in\so,\;\overline{\bob}\in\bR^3
		& 
		\\
		\pD\boA=0
		    \arrow[rr,Rightarrow, "\textrm{not needed}"'description,end anchor=west,start anchor={[xshift=17.5mm]},end anchor={[xshift=-8mm]},phantom]
		&
		&
		\vphantom{pD\pu-\boA=0}
		\end{tikzcd}
		\end{aligned}
		\end{equation}
	\end{description}
	As alluded to above, in the Cosserat model the mere prescription of $\left.\pu\right|_{\Gamma}$ suffices for existence and uniqueness. Moreover, the Cosserat model is redundant according to definition \ref{def: ours}.

	\begin{description}
		\item[Micromorphic family:]  
		
		\bigskip
		
		\item[(a)] classical micromorphic model with positive Cosserat couple modulus $\mu_{\text{c}}>0$ and only $\bou$ clamped at the 
		           
		           \vspace{-1.5mm}\hspace{-3mm} boundary $\Gamma$
		
		\begin{equation}
			\hspace{-3mm}
			\begin{aligned}
			\begin{tikzcd}[
			column sep=10mm,
			row sep=-3mm,
			every arrow/.append style={line width=0.7pt},
			/tikz/column 1/.append style={anchor=base east},
			/tikz/column 2/.append style={anchor=base east},
			/tikz/column 3/.append style={anchor=base west},
			every label/.append style = {font = \normalsize },
			every arrow/.append style={line width=0.7pt},
			]
			\text{D}\pu-\pP=0
			\arrow[dd, start anchor=north east, end anchor=south east, no head, xshift=0em, decorate, decoration={brace}, "\Longrightarrow" right=3pt]
			&
			\sym\text{D}\pu=0\;\overset{\textrm{Korn}}{\Longrightarrow}\;\pu=\overline{\boA}\!\cdot\!\boox+\overline{\bob}
			\arrow[dd, start anchor=north east, end anchor=south east, no head, xshift=0em, decorate, decoration={brace},"{\underset{\vphantom{\text{rigidity}}}{\overset{\text{rigidity \cite{neff2008curl}}}{\Longrightarrow}}}" right=3pt]
			&
			\\
			&
			&
			\qquad\text{D}\pu=\overline{\boA}=\pP \phantom{\underline{\bob}}\label{eqt:numerical_scheme1}
			\arrow[dddd, Rightarrow,"\;+\,\textrm{bc}",end anchor={[xshift=-4.45mm]}]
			\\
			\sym\pP=0 
			&
			\phantom{\overline{\bob}} \overline{\boA}\in\so,\;\overline{\bob}\in\bR^3  
			&
			\\
			&
			\phantom{\text{D}\pu=\overline{\boA}=\pP}
			&
			\\[2mm]
			\underbrace{\phantom{\pP=}\text{D}\pP=0}_{\mathclap{\substack{\\ \text{zero elastic energy,}\\ \text{classical micromorphic,}\\\mu_{\textrm{c}}>0}}}
			\arrow[rr,Rightarrow, "\textrm{not needed}"'description,end anchor=west,start anchor={[xshift=8mm]},end anchor={[xshift=-8mm]},phantom]
			&
			&
			\phantom{\pu-\pP=0}
			\\
			&
			&
			\hspace{-5mm}\left.\pu\right|_{\Gamma}=0\;\Rightarrow\;\text{D}\pu=  0=\overline{\boA}=\pP
			\\
			&
			&
			\phantom{\hspace{-5mm}\left.\pu\right|_{\Gamma}=0\bigg| }\Rightarrow\;(\pu,\pP)=0
			\end{tikzcd}
			\end{aligned}
			\end{equation}
		\end{description}
		Thus, prescribing $\left. \pu\right|_\Gamma=0 $ fixes the infinitesimal rotation $\overline{\boA}$, the constant translation $\overline{\bob}$ and the constant 
		microdistortion $\pP$.
		This means formulation (a) is \underline{redundant} according to definition \ref{def: ours} since the curvature measure $\text{D}\pP$ is not needed in order to arrive at uniform infinitesimal rotations.
		Here, it can be seen that it is sufficient to only clamp $\bou$ at the boundary $\Gamma$.
		
		\begin{description}
		\item[(b)]  classical micromorphic model with $\mu_{\text{c}}=0$ and consistent coupling boundary conditions
		\begin{equation}
		    \hspace{-10mm}
			\begin{tikzcd}[
				column sep=10mm,
				row sep=-3mm,
				every arrow/.append style={line width=0.7pt},
				/tikz/column 1/.append style={anchor=base east},
				/tikz/column 2/.append style={anchor=base east},
				/tikz/column 3/.append style={anchor=base west},
				every label/.append style = {font = \normalsize },
				every arrow/.append style={line width=0.7pt},
				ampersand replacement=\&
				]
				\sym\!(\text{D}\pu-\pP)=0
				\arrow[dd, start anchor=north east, end anchor=south east, no head, xshift=0em, decorate, decoration={brace}, "\Longrightarrow" right=3pt]
				\&
				\sym\text{D}\pu=0\;\overset{\textrm{Korn}}{\Longrightarrow}\;\pu=\overline{\boldsymbol{A}}\!\cdot\!\boox+\overline{\boldsymbol{b}}
				\arrow[dd, start anchor=north east, end anchor=south east, no head, xshift=0em, decorate, decoration={brace},"\Longrightarrow" right=3pt]
				\&
				\text{D}\pu=\overline{\boldsymbol{A}} 
				\arrow[dddd, start anchor=north east, end anchor=south east, no head, xshift=0em, decorate, decoration={brace},"{\overset{\substack{\textrm{consist.} \\ \textrm{coupl.}}}{\Longrightarrow}\!\!\!
					\begin{blockarray}{l c l}
					     \begin{block}{ l c l}
					     \\
					     \left. \text{D}\pu\times\pnu\,\right|_\Gamma =\left. \pP\times\pnu\,\right|_\Gamma  \;\Rightarrow\\ 
					     \;\;\big.  \overline{\boA}\times\pnu\,\big|_\Gamma =\big. \overline{\overline{\boA}}\times\pnu\,\big|_\Gamma \;\,\Rightarrow\\
					     \end{block}
					\;\;\quad\qquad\overline{\boA}=\overline{\overline{\boA}} \\
					\end{blockarray}}" right=3mm]
				\\
				\&
				\&
				\vphantom{\text{D}\pu=\overline{\boldsymbol{A}} }
				\\
				\sym\pP=0 
				\&
				\phantom{\overline{\boldsymbol{b}}} \overline{\boldsymbol{A}}\in\so,\;\overline{\boldsymbol{b}}\in\bR^3  
				\&
				\\
				\&
				\phantom{\text{D}\pu=\overline{\boldsymbol{A}}=\pP}
				\&
				\\[2mm]
				\underbrace{\phantom{\pu-\pP=0}\text{D}\pP=0}_{\substack{\\ \text{zero elastic energy,}\\ \text{classical micromorphic,}\\\mu_{\textrm{c}}\equiv0}}
				\arrow[rr,Rightarrow,"\textrm{curvature needed, but too much, see } \eqref{qualification}"'{font = \footnotesize},end anchor=west,start anchor={[xshift=8mm]},end anchor={[xshift=-8mm]}]
				\&
				\&
				\;\;\pP=\overline{\overline{\boldsymbol{A}}}
				\\
				\&
				\hphantom{.}
				\arrow[r,Rightarrow,"{\small \!\!+\,\textrm{b.c.}}"below=1mm,start anchor={[xshift=1.5mm]},end anchor={[xshift=-2mm]}]
				\&
				\left.\pu\right|_{\Gamma}=0\;\Longrightarrow\; \text{D}\pu=\pP=0
			\end{tikzcd} 
		\end{equation}
		\end{description}
		We observe that the consistent coupling boundary condition replaces a positive Cosserat couple modulus $\muc>0$ in order to arrive at constant infinitesimal rotations.
		This case is non-redundant in the new definition \ref{def: ours}: the curvature measure $\text{D}\pP$ is needed in this argument.
		
		\begin{description}
		\item[(c)]  relaxed micromorphic model with $\mu_{\text{c}}=0$ and consistent coupling boundary conditions
          \begin{equation}\label{qualification}
			\hspace{-10mm}
			\begin{tikzcd}[
				column sep=10mm,
				row sep=-3mm,
				every arrow/.append style={line width=0.7pt},
				/tikz/column 1/.append style={anchor=base east},
				/tikz/column 2/.append style={anchor=base east},
				/tikz/column 3/.append style={anchor=base west},
				every label/.append style = {font = \normalsize },
				every arrow/.append style={line width=0.7pt},
				ampersand replacement=\&
				]
				\sym\!(\text{D}\pu-\pP)=0
				\arrow[dd, start anchor=north east, end anchor=south east, no head, xshift=0em, decorate, decoration={brace}, "\Longrightarrow" right=3pt]
				\&
				\sym\text{D}\pu=0\;\overset{\textrm{Korn}}{\Longrightarrow}\;\pu=\overline{\boldsymbol{A}}\!\cdot\!\boox+\overline{\boldsymbol{b}}
				\arrow[dddd, start anchor=north east, end anchor=south east, no head, xshift=0em, decorate, decoration={brace},"{
					\begin{blockarray}{l c l}
						\begin{block}{ l\} c l}
							\text{D}\pu=\overline{\boA} \;\;
							&
							\kern0.1cm \multirow{2}{*}{$\overset{\substack{\textrm{consist.} \\ \textrm{coupl.}}}{\Longrightarrow}$}  
							& 
							\kern-0.2cm \left. \text{D}\pu\times\pnu\,\right|_\Gamma =\left. \pP\times\pnu\,\right|_\Gamma  \;\Rightarrow
							\\ 
							\;\;\pP=\overline{\overline{\boA}}  
							&  
							& 
							\kern-0.6cm\;\;\;\quad \big.  \overline{\boA}\times\pnu\,|_\Gamma =\big.\overline{\overline{\boA}}\times\pnu\,\big|_\Gamma \;\Rightarrow\\
						\end{block}
						& 
						& 
						\quad\qquad\overline{\boA}=\overline{\overline{\boA}} \\
				\end{blockarray}}" right=5mm]
				\&
				\\
				\&
				\&
				\phantom{\text{D}\pu=\overline{\boldsymbol{A}}=\pP \phantom{\underline{\bob}}}
				\\
				\sym\pP=0 
				\&
				\phantom{\overline{\bob}} \overline{\boA},\overline{\overline{\boA}}\in\so,\;\overline{\bob}\in\bR^3
				\&
				\phantom{\overline{\bob}}  
				\\
				\&
				\phantom{\text{D}\pu=\overline{\boA}=\pP}
				\&
				\\[2mm]
				\underbrace{\phantom{\pu-\pP}\Curl\pP=0}_{\substack{\\ \textrm{zero elastic energy}\\ \! \!\!\text{relaxed micromorphic,}\\ \mu_{\textrm{c}}\,=\,0}}
				\arrow[r,Rightarrow,start anchor={[xshift=5mm]},end anchor={[xshift=15mm]}]
				\&
				\phantom{\text{D}\pu=\overline{\boA}=\pP\underline{\bob}} \pP=\overline{\overline{\boA}}
				\&
				\left.\pu\right|_{\Gamma}=0\;\;\,\Longrightarrow\;\;\text{D}\pu=\pP=0\\
				\&
				\hphantom{.}
				\arrow[ur,Rightarrow,"\!\!\!\!\!\!\!\!\!{\small +\,\textrm{b.c.}}"'sloped,start anchor={[yshift=2mm]},end anchor=south west]
				\&
				\phantom{\left.\pu\right|_{\Gamma}=0}\;\;\,\Longrightarrow\;\;(\pu,\pP)=0
			\end{tikzcd}
		\end{equation}
	    \end{description}
		Indeed, by the Nye's formula \cite{lewintan2020korn,nye1953some} for a skew-symmetric matrix field $\boA\in\so$ it holds:
		\begin{equation}\label{Nye}
		\D\!\axl \boA =\frac12\tr\!(\Curl \boA)\!\cdot\!\id-(\Curl \boA)^T\quad\Leftrightarrow \quad \Curl \boA=\tr\!(\D\!\axl\boA)\id-(\D\!\axl\boA)^{T},
		\end{equation}
		or in short $\D\! \boA = L(\Curl \boA)$ with a constant coefficient linear operator $L$. Thus, assuming for a matrix field $\pP$
		$$
		\sym \pP \equiv0 \quad\text{and}\quad \Curl \pP \equiv0,
		$$
		we must have $\pP=\boA(x)$, with a skew-symmetric matrix field, and by the Nye formula it follows $\D\! \boA\equiv0$ and hence $\pP=\overline{\overline{\boA}}\in\so$.
		This case is non-redundant according to definition \ref{def: ours}.
		\begin{description}
		\item[Second gradient model:] 
		$$
		\begin{tikzcd}[
			column sep=15mm,
			row sep=2mm,
			every arrow/.append style={line width=0.7pt},
			/tikz/column 1/.append style={anchor=base east},
			/tikz/column 3/.append style={anchor=base west},
			]
			\sym\text{D}\pu=0\;\;\,
			\arrow[r,Rightarrow,start anchor={[xshift=2mm]},end anchor={[xshift=-2mm]}]
			&
			\pu=\overline{\boldsymbol{A}}\!\cdot\!\boox+\overline{\boldsymbol{b}}
			\arrow[r,Rightarrow,start anchor={[xshift=2mm]},end anchor={[xshift=-2mm]}]
			&
			\text{D}\pu=\overline{\boldsymbol{A}}
			\\
			\underbrace{\text{D}\,(\text{D}\pu)=0}_{\textrm{zero elastic energy}}
			&
			\textrm{not needed,}
			\arrow[r,Rightarrow,"\!\!+\,\textrm{b.c.}"below=1mm,start anchor={[xshift=3.5mm]},end anchor={[xshift=-2mm]}]
			&
			\left.\pu\right|_{\Gamma}=0\;\Longrightarrow\;\pu\equiv0.
		\end{tikzcd}
		$$
		\end{description}
		Therefore the second gradient model is redundant, which is well known.
		Looking at the last subsection, this leads us to the question: what should be zero energy modes for  a generalized continuum model? We have seen that we arrive always at
		\begin{equation}
			\begin{aligned}\pu & =\overline{\boA}\!\cdot\!\boox+\overline{\boldsymbol{b}} & \textrm{classical for displacement}\\
				\pP & =\overline{\overline{\boldsymbol{A}}} & \overline{\boldsymbol{A}},\overline{\overline{\boldsymbol{A}}}\in\so.
			\end{aligned}
		\end{equation}
		In the latter, the missing crucial information is whether $\overline{\boldsymbol{A}}=\overline{\overline{\boldsymbol{A}}}$. As we have seen, this can be supplied by boundary conditions, e.g. the consistent coupling condition $\left.\text{D}\pu\times\pnu\,\right|_{\Gamma}=\left.\pP\times\pnu\,\right|_{\Gamma}$
		which then implies 
		\begin{equation}
			\label{eq:con_coup_cond_skew}
			\overline{\boldsymbol{A}}=\overline{\overline{\boldsymbol{A}}}.
		\end{equation}
	\section{Passage to a second gradient continuum}
	The second gradient continuum deserves a little detour since the effect of imposing higher order boundary conditions significantly modifies the solution. Indeed, this can be seen already for a didactic $1D-$ example.
	
\begin{figure}[H]
		\centering
		\begin{subfigure}{0.43\textwidth}
			\centering
			\includegraphics[width=\textwidth]{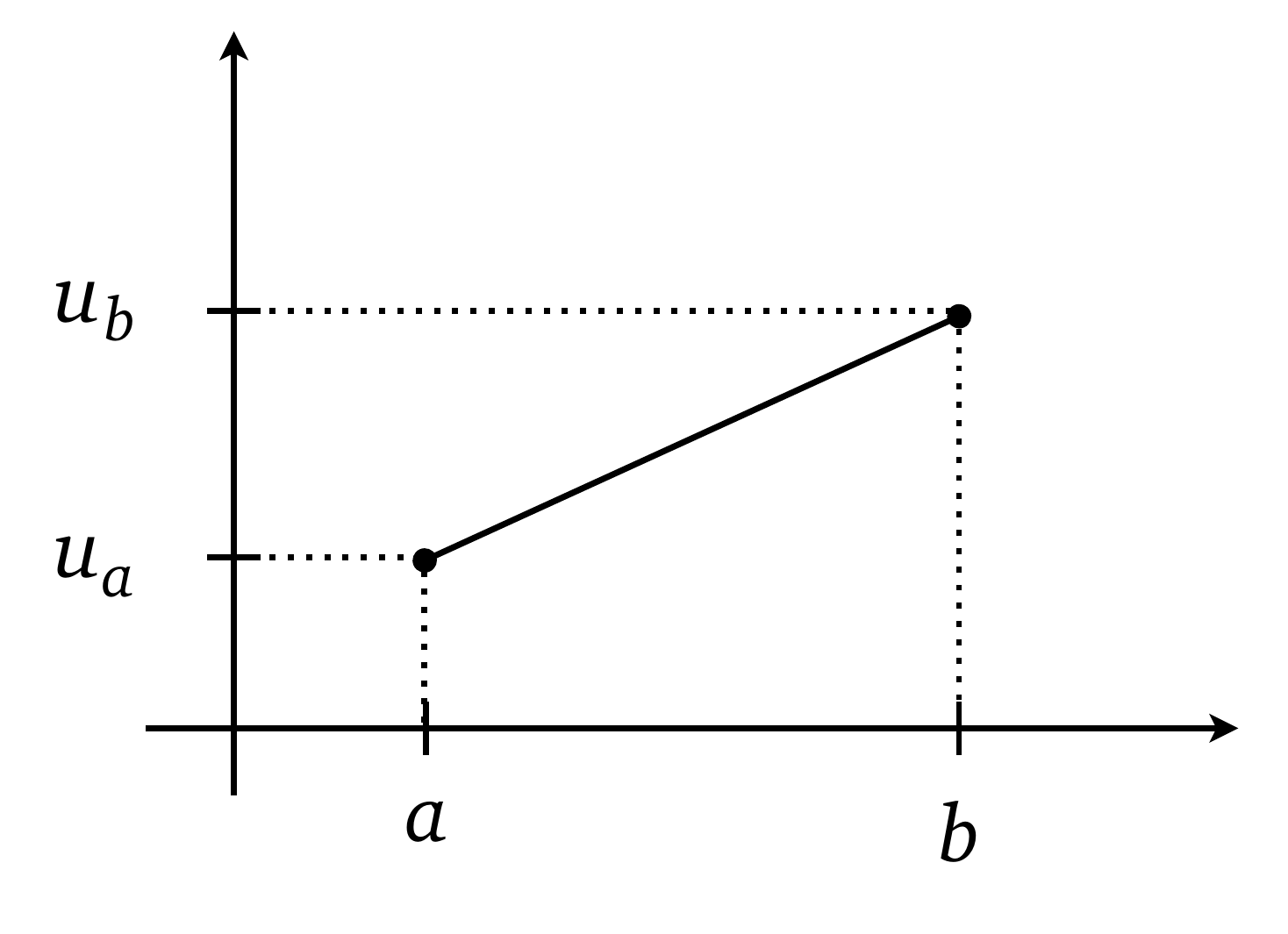}
		\end{subfigure}
		\hfill
		\begin{subfigure}{0.55\textwidth}
			\centering
			\includegraphics[width=\textwidth]{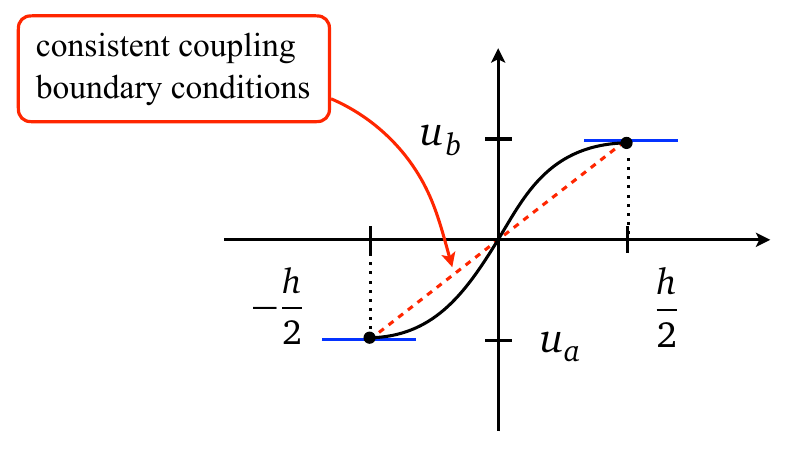}
		\end{subfigure}
		\caption{(Left) Unique homogeneous solution to a strictly convex first order problem \eqref{convex problem}. (Right)
			The ``normal derivative'' boundary conditions $u'\left(-\frac{h}{2}\right)=u'\left(\frac{h}{2}\right)=0$ trigger a non-homogeneous solutions when considering the second gradient formulation \eqref{non convex problem}. \label{fig.solution}
		}
		\label{fig:second}
	\end{figure}

	Let us look first at a classical first order problem
	\begin{equation}\label{convex problem}
		\int_{a}^b\left|u'\right|^{2}dx\longrightarrow\min\quad\textrm{subjected to the constraints}\quad u\left(a\right)=u_{a},\;u\left(b\right)=u_{b}.
	\end{equation}
	Since the integrand $\left| u' \right|^2$ is strictly convex in $u'$, the unique, homogenous solution is (see Fig. \ref{fig.solution} left)
	\begin{equation}\label{solution}
		u\left(x\right)=\frac{u_{b}-u_{a}}{b-a}\,\left(x-a\right)-u_{a}.
	\end{equation}
	Now consider the second gradient augmented formulation \eqref{non convex problem} together with describing additional ``normal'' derivatives at the boundary
	\begin{equation}\label{non convex problem}
		\int_{-\frac{h}{2}}^{\frac{h}{2}}\left|u'\right|^{2}+L_c^2\left|u''\right|^{2}dx\longrightarrow\min\quad\textrm{subjected to}\quad\begin{cases}
			u\left(-\frac{h}{2}\right)=u_{a},\;u\left(\frac{h}{2}\right)=u_{b},\\
			u'\left(-\frac{h}{2}\right)=0,\;u'\left(\frac{h}{2}\right)=0. 
		\end{cases} 
	\end{equation}
	The solution is now given by (considering $b=-a=h/2$ and $u_b=-u_a=(\boldsymbol{\gamma}\, h)/2$ without loosing any generality)
	\begin{equation}\label{solution second}
		u\left(x\right)=
		\frac{
			\frac{x}{h}
			\cosh \left(\frac{h}{2L_{\text{c}}}\right)
			-
			\sinh \left(\frac{x}{L_{\text{c}}}\right)
			\frac{L_{\text{c}}}{h}
		}{
			\frac{1}{2} \cosh \left(\frac{h}{2L_{\text{c}}}\right)
			-
			\sinh \left(\frac{h}{2L_{\text{c}}}\right)
			\frac{L_{\text{c}}}{h}
		}
		\frac{\boldsymbol{\gamma}  h}{2}
		\, .
	\end{equation}
	Here, the non-homogeneous ``normal''-derivative boundary conditions trigger a non-homogeneous solution that depends on $L_{\rm c}>0$ (see Fig.\,\ref{fig.solution} right). 
	However, if we study \eqref{non convex problem} only with $u\left(a\right)=u_{a},\;u\left(b\right)=u_{b}$ (the consistent coupling case) then we retrieve the homogeneous solution \eqref{solution} (see Fig.~\ref{fig:second}).
	
	
	\subsection{Second gradient continuum and micromorphic penalty formulation}
	The classical micromorphic model with positive Cosserat couple modulus $\muc>0$ is also often used to approximate, via a
	penalty approach, solutions of a corresponding second gradient continuum.
	Indeed, consider
	\begin{align}
		\int_{\Omega}\mue\left\lVert \sym\!(\pD\pu-\pP)\right\rVert^{2} 
		& 
		+\muc\left\lVert \skew\!(\pD\pu-\pP)\right\rVert^{2}
		+\dfrac{\le}{2}\tr\!^{2}(\pD\pu-\pP)\label{eq:second gradient}
		\\
		& 
		+\mum\left\lVert \sym\pP\right\rVert^{2}
		+\dfrac{\lam}{2}\tr\!^{2}(\pP)
		+\frac{\mu\,L_{\text{c}}^{2}}{2}\left\lVert \pD\pP\right\rVert^{2}\text{dV}\longrightarrow\min\quad\left(\pu,\pP\right)\nonumber 
	\end{align}
	and interpret $\mu_{\text{e}},\lambda_{\text{e}},\mu_{\text{c}}\fr\infty$
	as penalty parameters. In the penalty limit, we must have $\text{D}\pu=\pP$ and (\ref{eq:second gradient}) turns formally into
	\begin{equation}
		\int_{\Omega}\mu_{\text{micro}}\left\lVert \text{sym}\,\text{D}\pu\right\rVert ^{2}
		+\dfrac{\lam}{2}\tr\!^{2}(\pD\pu)
		+\frac{\mu\,L_{\text{c}}^{2}}{2}\left\lVert \pD(\pD\pu)\right\rVert^{2}\text{dV}\longrightarrow\min\quad\pu.\label{eq:limit problem}
	\end{equation}
	If we assume the consistent coupling boundary conditions in (\ref{eq:second gradient})
	\begin{equation}
		\left.\pu\right|_{\Gamma}=\widehat{\pu},\qquad\left.\text{D}\pu\times\pnu\right|_{\Gamma}=\left.\pP\times\pnu\right|_{\Gamma},\label{eq:consistent_bc}
	\end{equation}
	then the limit (\ref{eq:limit problem}) generates just the ``first order boundary conditions'' $\left.\pu\right|_{\Gamma}=\widehat{\pu}$
	and the corresponding third order stress tensor $\boldsymbol{\mathfrak{m}}$ (see eq.(\ref{eq:equi_Strain_Grad})) verifies Neumann
	boundary conditions throughout.
	
	If, on the other hand, in (\ref{eq:second gradient}) we assume
	\begin{equation}
		\left.\pu\right|_{\Gamma}=\widehat{\pu},\qquad\left.\pP\right|_{\Gamma}=\widehat{\pP}\label{eq:Classical bc}
	\end{equation}
	then this generates in the penalty limit a possible inconsistency: $\left.\text{D}\pu\times\pnu\right|_{\Gamma}$
	is already given through $\left.\pu\right|_{\Gamma}=\widehat{\pu}$, but also $\pP=\text{D}\pu$ and $\left.\pP\right|_{\Gamma}=\widehat{\pP}$ hold,
	which does not necessarily imply that $\left.\text{D}\pu\times\pnu\right|_{\Gamma}=\left.\widehat{\pP}\times\pnu\right|_{\Gamma}$.
	This may lead to a boundary layer. Note that both sets of boundary conditions lead to a unique solution.
	
	\subsection{Comparison of boundary conditions}
	In this section we present a summary of the possible boundary conditions in Table~\ref{tab:tab1}.
	\begin{table}[H]
		\rowcolors{1}{}{gray!15}
		{\scriptsize
			\global\tabulinesep =3mm 
			\begin{tabu}to 170mm { X[0.6,c] !{\color{myblue}\vrule width 1.2pt} X[c] | X[c]  !{\color{myblue}\vrule width 1.2pt}  X[c]  }   \arrayrulecolor{myblue}
				\tabucline [0.8mm]{}  
				&
				\multicolumn{2}{ c  !{\color{myblue}\vrule width 1.2pt}}{   { \Large$\textcolor{black}{\Gamma}$}    }
				&
				{ \Large$ \textcolor{black}{\partial\Omega\setminus\overline{\Gamma}}$}  
				\\
				\hline
				\hline
				&
				\textcolor{myblue}{\textbf{\textsf{Dirichlet}}}
				&
				\textcolor{myblue}{\textbf{\textsf{Neumann}}}
				&
				\textcolor{myblue}{\textbf{\textsf{Neumann}}}
				\\
				\tabucline [1.2pt]{}
				\textcolor{myblue}{\textbf{\textsf{linear elasticity}}}
				&
				$
				\left.\pu\right|_{\Gamma}=\widehat{\pu}
				$
				&
				$
				\diagup
				$
				&
				$
				\left.\bosigma\,\right|_{\partial\Omega\setminus\overline{\Gamma}}  =0
				$
				\\
				\tabucline [1.2pt]{}
				\textcolor{myblue}{\textbf{\textsf{Cosserat with fixed rotations}}}
				&
				$
				\left\{ \ensuremath{\begin{aligned}\left.\pu\right|_{\Gamma} & =\widehat{\pu}\\
					\left.\boA\right|_{\Gamma} & =\widehat{\boA}|_{\Gamma}
					\end{aligned}
				}\right.
				$
				&
				$
				\diagup
				$
				&
				$
				\left\{ \ensuremath{\begin{aligned}\left.\bosigma\!\cdot\!\pnu\,\right|_{\partial\Omega\setminus\overline{\Gamma}} & =0\\
					\left.\bom\times\pnu\,\right|_{\partial\Omega\setminus\overline{\Gamma}} & =0
					\end{aligned}
				}\right.
				$
				\\
				\tabucline [1.2pt]{}
				\textcolor{myblue}{\textbf{\textsf{Cosserat with consistent coupling}}}
				&
				$
				\left\{ \ensuremath{\begin{aligned}\left.\pu\right|_{\Gamma} & =\widehat{\pu}\\
					\left.\boA\times\pnu\,\right|_{\Gamma} & =\left.\skew\pD\pu\times\pnu\right|_{\Gamma}
					\end{aligned}
				}\right.
				$
				&
				$
				\diagup
				$
				&
				$
				\left\{ \ensuremath{\begin{aligned}\left.\bosigma\!\cdot\!\pnu\,\right|_{\partial\Omega\setminus\overline{\Gamma}} & =0\\
					\left.\bom\times\pnu\right|_{\partial\Omega\setminus\overline{\Gamma}}& =0
					\end{aligned}
				}\right.
				$
				\\
				\tabucline [1.2pt]{}
				\textcolor{myblue}{\textbf{\textsf{Cosserat with free rotations}}}
				&
				$
				\left.\pu\right|_{\Gamma}=\widehat{\pu}
				$
				&
				$
				\left.\bom\times\pnu\,\right|_{\Gamma} =0
				$
				&
				$
				\left\{ \ensuremath{\begin{aligned}\left.\bosigma\!\cdot\!\pnu\,\right|_{\partial\Omega\setminus\overline{\Gamma}} & =0\\
					\left.\bom\times\pnu\,\right|_{\partial\Omega\setminus\overline{\Gamma}} & =0
					\end{aligned}
				}\right.
				$
				\\
				\tabucline [1.2pt]{}
				\textcolor{myblue}{\textbf{\textsf{second gradient with free normal derivative}}}
				&
				$
				\left.\pu\right|_{\Gamma}=\widehat{\pu} \newline \left(\;\Rightarrow\;\left.\text{D}\pu\times\pnu\,\right|_{\Gamma}\;\textrm{is fixed}\;\right)
				$
				&
				$
				\begin{aligned}\left.\left(\gbm\!\cdot\!\pnu\right)\!\cdot\!\left(\pnu\otimes\pnu\,\right)\right|_{\Gamma} & =0\\
					\overset{\textrm{see } \eqref{conto second gradient}}{\Leftrightarrow}\quad\left.\left(\gbm\!\cdot\!\pnu\right)\!\cdot\!\pnu\,\right|_{\Gamma} & =0
				\end{aligned}
				$
				&
				$
				\left\{ \ensuremath{\begin{aligned}\left.\boldsymbol{\eta}\,\right|_{\partial\Omega\setminus\overline{\Gamma}} & =0\;\quad \textrm{see } \eqref{bc second gradient1}\\
						\left.\gbm\!\cdot\!\pnu\,\right|_{\partial\Omega\setminus\overline{\Gamma}} & =0
					\end{aligned}
				}\right.
				$
				\\
				\tabucline [0.2mm]{}
				\textcolor{myblue}{\textbf{\textsf{second gradient with fixed normal derivative}}}
				&
				$
				\left\{ \ensuremath{
					\begin{aligned}
					    \left.\pu\right|_{\Gamma} & =\widehat{\pu}
					    \\
						\left.\text{D}\pu\!\cdot\!\pnu\right|_{\Gamma}   & =\widehat{\bog}
					\end{aligned}
				}\right.
				$
				&
				$
				\diagup
				$
				&
				$
				\left\{ \ensuremath{
					\begin{aligned}
					    \left.\boldsymbol{\eta}\,\right|_{\partial\Omega\setminus\overline{\Gamma}} & =0
					    \\
						\left.\gbm\!\cdot\!\pnu\,\right|_{\partial\Omega\setminus\overline{\Gamma}} & =0
					\end{aligned}
				}\right.
				$
				\\
				\tabucline [1.2pt]{}   
				\textcolor{myblue}{\textbf{\textsf{micromorphic + classical bc}}}
				&
				$
				\left\{ \ensuremath{\begin{aligned}\left.\pu\right|_{\Gamma} & =\widehat{\pu}\\
						\left.\pP\right|_{\Gamma}  & =\widehat{\pP}
					\end{aligned}
				}\right.
				$
				&
				$
				\diagup
				$
				&
				$
				\left\{ \ensuremath{
					\begin{aligned}
					    \left.\bosigma\!\cdot\!\pnu\,\right|_{\partial\Omega\setminus\overline{\Gamma}} & =0\\
						\left.\gbm\!\cdot\!\pnu\,\right|_{\partial\Omega\setminus\overline{\Gamma}} & =0
					\end{aligned}
				}\right.
				$
				\\
				\tabucline [0.2mm]{}   
				\textcolor{myblue}{\textbf{\textsf{micromorphic + consistent bc}}}
				&
				$
				\left\{ \ensuremath{\begin{aligned}\left.\pu\right|_{\Gamma} & =\widehat{\pu}\\
						\left.\pP\times\pnu\,\right|_{\Gamma} & =\left.\text{D}\pu\times\pnu\right|_{\Gamma}
					\end{aligned}
				}\right.
				$
				&
				$
				\begin{aligned}\left.\left(\gbm\!\cdot\!\pnu\right)\!\cdot\!\left(\pnu\otimes\pnu\,\right)\right|_{\Gamma} & =0\\
					\overset{\textrm{see } \eqref{conto second gradient}}{\Leftrightarrow}\quad\left.\left(\gbm\!\cdot\!\pnu\right)\!\cdot\!\pnu\,\right|_{\Gamma} & =0
				\end{aligned}
				$
				&
				$
				\left\{ \ensuremath{\begin{aligned}\left.\bosigma\!\cdot\!\pnu\,\right|_{\partial\Omega\setminus\overline{\Gamma}} & =0\\
						\left.\gbm\!\cdot\!\pnu\,\right|_{\partial\Omega\setminus\overline{\Gamma}} & =0
					\end{aligned}
				}\right.
				$
				\\
				\tabucline [1.2pt]{}
				\textcolor{myblue}{\textbf{\textsf{relaxed micromorphic written in $\Curl\pP$}}}
				&
				$
				\left\{ \ensuremath{\begin{aligned}\left.\pu\right|_{\Gamma} & =\widehat{\pu}\\
						\left.\pP\times\pnu\,\right|_{\Gamma} & =\left.\text{D}\pu\times\pnu\right|_{\Gamma}
					\end{aligned}
				}\right.
				$
				&
				$
				\diagup
				$
				&
				$
					\left\{ \ensuremath{\begin{aligned}\left.\bosigma\!\cdot\!\pnu\,\right|_{\partial\Omega\setminus\overline{\Gamma}} & =0\\
					\left.\boldsymbol{m}\times\pnu\right|_{\partial\Omega\setminus\overline{\Gamma}}& =0
					\end{aligned}
				}\right.
				$
				\\
				\tabucline [1.2pt]{}
				\textcolor{myblue}{\textbf{\textsf{relaxed micromorphic written in $\text{D}\pP$}}}
				&
				$
				\left\{ \ensuremath{\begin{aligned}\left.\pu\right|_{\Gamma} & =\widehat{\pu}\\
					\left.\pP\times\pnu\,\right|_{\Gamma} & =\left.\text{D}\pu\times\pnu\right|_{\Gamma}
					\end{aligned}
				}\right.
				$
				&
				$
				\diagup
				$
				&
				$
				\left\{ \ensuremath{\begin{aligned}\left.\bosigma\!\cdot\!\pnu\,\right|_{\partial\Omega\setminus\overline{\Gamma}} & =0\\
					\left.\gbm\!\cdot\!\pnu\right|_{\partial\Omega\setminus\overline{\Gamma}}& =0
					\end{aligned}
				}\right.
				$
				\\
				\tabucline [0.8mm]{}  
			\end{tabu}
		}
		\caption{Comparison between different possibilities of imposing the boundary conditions in different models. We remember that the definition of the generalized traction in the second gradient model is given by
			$\boldsymbol{\eta}\defi\left.(\bosigma-\DIV\gbm)\!\cdot\!\pnu\right|_{\partial\Omega\setminus\overline{\Gamma}}\left.\,-\;\Div_{\partial\Omega,\top}\left(\gbm\!\cdot\!\pnu\right)\right|_{\partial\Omega\setminus\overline{\Gamma}}$. These expressions are derived in the appendix.}
			\label{tab:tab1}
	\end{table}	
	\section{Another invariance condition in non-linear and other higher order models - pure strain coupling}
	Let us again consider the generic energy expression for the classical linear Mindlin-Eringen micromorphic model for $\muc=0$, i.e.
	\begin{equation}\label{9.1}
		W(\text{D}\boldsymbol{u},\boldsymbol{P},\text{D}\boldsymbol{P})=\left\Vert \sym\!(\text{D}\boldsymbol{u}-\boldsymbol{P})\right\Vert ^{2}+\left\Vert \sym\boldsymbol{P}\right\Vert ^{2}+\left\Vert \text{D}\boldsymbol{P}\right\Vert ^{2}.
	\end{equation}
	It is easy to observe that the energy $W$ is invariant under the \underline{global transformation} 
	\begin{equation}\label{invariance1}
		\left(\pu,\pP\right)\longrightarrow\left(\overline{\boA}_{1}\!\cdot\!\pu+\overline{\bob},\pP+\overline{\boA}_{2}\right)
	\end{equation}
	for all constant infinitesimal rotations $\overline{\boA}_{1},\overline{\boA}_{2}\in\so$. Infinitesimal objectivity (frame indifference) is also included by considering arbitrary $\overline{\boA}_{1}=\overline{\boA}_{2} \in \so$.
	
	In fact, there is a non-linear counterpart for this model for the deformation $\bovarphi:\Omega\fr\bR^3$, $\bovarphi(\boox)=\boox+\pu(\boox)$, and the finite micro-distortion $\widehat{\pP}=\id+\pP$,
	\begin{equation}\label{9.2}
		\widehat{W}\big(\text{D}\bovarphi,\widehat{\pP},\text{D}\widehat{\pP}\big)=\big\Vert \text{D}\bovarphi^{T}\!\!\!\cdot\!\text{D}\bovarphi-\widehat{\pP}^{T}\!\!\!\cdot\!\widehat{\pP}\big\Vert ^{2}+\big\Vert \widehat{\pP}^{T}\!\!\!\cdot\!\widehat{\pP}-\id\big\Vert ^{2}+\big\Vert \text{D}\widehat{\pP}\big\Vert ^{2}
	\end{equation}
	which is now invariant under the \underline{global transformation}  
	\begin{equation}
		(\bovarphi,\widehat{\pP})\longrightarrow\big(\overline{\boQ}_{1}\!\cdot\!\bovarphi+\overline{\bob},\overline{\boQ}_{2}\!\cdot\!\widehat{\pP}\big)
	\end{equation}
	for all constant rotations $\overline{\boQ}_{1},\overline{\boQ}_{2}\in\SO$, mirroring \eqref{invariance1}.
	
\begin{figure}[H]
		\centering
		\includegraphics[width=\textwidth]{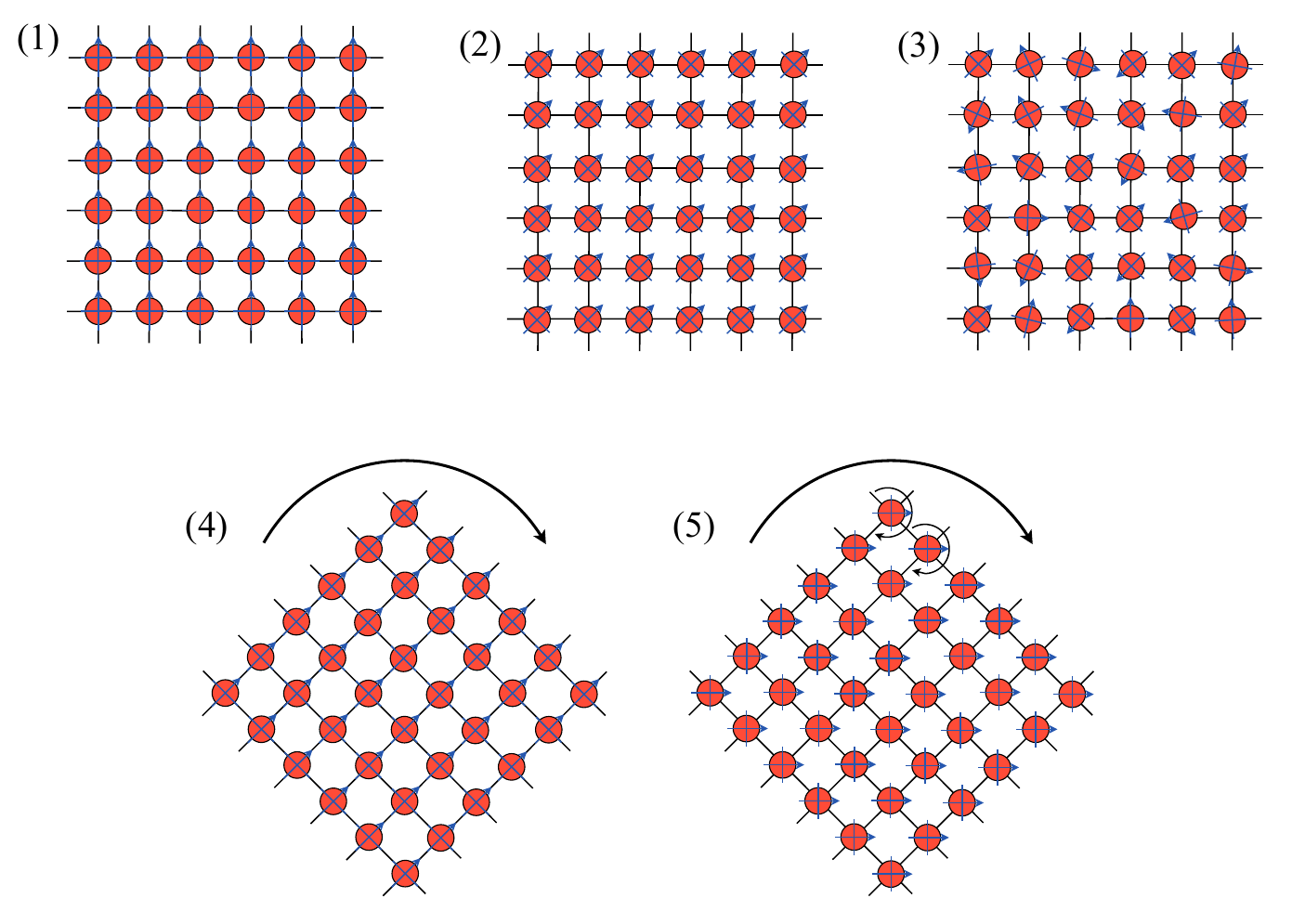}
		\caption{
			(1) shows the perfect alignement of the assumed microstructure with the underformed lattice. It represents a zero energy mode for the micromorphic model with $\muc\geqslant0$. In (2) the lattice remains undeformed and the microstructure is homogeneously rotated. This state represents a zero energy mode for the micromorphic model with zero Cosserat couple modulus $\muc=0$ (also for the relaxed micromorphic model). (3) shows an underformed lattice when the microstructure is inhomogeneously rotated. This represents only a zero energy mode, if $\muc=0$ and the curvature energy satisfies eq.\eqref{ccc}. Image (4) represents a rigidly rotated lattice while the microstructure is perfectly aligned and (5) represents a rigidly rotated lattice in which the microstructure is homogeneously rotated but not aligned with the lattice. Both (4) and (5) are zero energy modes for $\muc\geqslant0$ and $\mu_c=0$, respectively. The consistent coupling boundary condition will always align the microstructure rotation with the lattice without needing $\mu_c>0$. The invariance alluded to in image (3) should not be confounded with any concept of isotropy \cite{ebobisse2020fourth}.
		}
		\label{fig:net_rotation}
	\end{figure}

	Both formulations \eqref{9.1} and \eqref{9.2} can be modified, such that the global $\mathfrak{so}(3)-$invariance turns even into a corresponding \underline{local} $\so$-invariance (local SO(3)-invariance, respectively).
	This can be achieved by considering new higher order prototype elastic energies (see Fig.~\ref{fig:net_rotation})
	\begin{align}
		W_{\textrm{inc}}(\text{D}\pu,\boldsymbol{P},\text{D}^{2}\pP) 
		& =\left\Vert \sym\!(\text{D}\pu-\pP)\right\Vert ^{2}+\left\Vert \sym\pP\right\Vert ^{2}+\left\Vert \text{inc}(\sym\pP)\right\Vert ^{2},\label{ccc}
		\\
		\widehat{W}_{\textrm{inc}}(\text{D}\varphi,\widehat{\pP},\text{D}^{2}\widehat{\pP}) 
		& =\bigl\Vert\text{D}\bovarphi^{T}\!\!\cdot\!\text{D}\bovarphi-\widehat{\pP}^{T}\!\!\!\cdot\!\widehat{\pP}\bigr\Vert^{2}+\bigl\Vert\widehat{\pP}^{T}\!\!\!\cdot\!\widehat{\pP}-\id\bigr\Vert^{2}+\bigl\Vert\gR\,\bigl(\widehat{\pP}^{T}\!\!\!\cdot\!\widehat{\pP}\bigr)\bigr\Vert^{2}, \label{cccc}
	\end{align}
	where
	\begin{equation}
		\text{inc}(\sym\pP)\defi\Curl\!\!\left[\left(\Curl\sym\pP\right)^{T}\right]
	\end{equation}
	is the Kröner incompatibility tensor \cite{ebobisse2020fourth}  and 
	\begin{equation}\label{def-Riem-curvature}
		\gR^i_{ jkl}\bigl(\widehat{\pP}^{T}\!\!\!\cdot\!\widehat{\pP}\bigr)\coloneqq\frac{\partial\Gamma^i_{jl}}{\partial
			x_k}-\frac{\partial\Gamma^i_{jk}}{\partial x_l}
		+\Gamma^m_{jl}\,\Gamma^i_{mk}-\Gamma^m_{jk}\,\Gamma^i_{ml}\,,
	\end{equation}
	is the Riemann-Christoffel tensor \cite{ebobisse2020fourth} with 
	$$
	\Gamma^l_{ij}\coloneqq\frac12g^{kl}\!\left[\frac{\partial
		g_{ki}}{\partial x_j}+\frac{\partial g_{jk}}{\partial
		x_i}-\frac{\partial g_{ij}}{\partial x_k}\right]
	$$ being the
	Christoffel symbols of the second kind in the metric $(g_{ij})=\bigl(\widehat{\pP}^{T}\!\!\!\cdot\!\widehat{\pP}\bigr)_{ij}$
	whose inverse is $(g^{ij})=\bigl( \bigl(\widehat{\pP}^{T}\!\!\!\cdot\!\widehat{\boldsymbol{P}}\bigr)^{-1}\bigr) _{ij}$. Note that \cite{ebobisse2020fourth} 
	$$
	\gR\left(\left(\id+\pP\right)^{T}\!\!\!\cdot\!\left(\id+\pP\right)\right)=\text{inc}(\sym\pP)+\textrm{h.o.t}.\;.
	$$
	 It is clear that both formulations \eqref{ccc} and \eqref{cccc} are not equivalent to \eqref{9.1} and \eqref{9.2}, respectively. 
	Note that new invariance-condition \eqref{invariance1} cannot be satisfied in the classical micromorphic model with generic energy (\ref{9.1}), corresponding to $\muc>0$, which is, however, still infinitesimally objective, i.e. invariant under 
	\begin{equation}
		\left(\pu,\pP\right)\longrightarrow\left(\overline{\boA}\!\cdot\!\pu+\overline{\bob},\pP+\overline{\boA}\right)
	\end{equation}
	for all constant $\overline{\boA}\in\so$. In any case, the global indeterminacy in \eqref{invariance1} is removed (gauged) by applying the consistent coupling boundary condition 
	\begin{equation}
		\left.\text{D}\pu\times\pnu\,\right|_{\Gamma}=\left.\pP\times\pnu\,\right|_{\Gamma}
	\end{equation}
	which can therefore be also interpreted as fixing the gauge.

	\section{Discussion}
	
	For two different sets of boundary conditions the classical Mindlin-Eringen \cite{mindlin1964micro,eringen1969micromorphic} micromorphic model allows for a unique solution in the same product space $\left(\pu,\pP\right)\in H^{1}\!\left(\Omega,\bR^{3}\right)\times H^{1}\!\left(\Omega,\bR^{3\times3}\right)$.
	The difference comes from the requirement on the microdistortion $\pP$ at the part of the boundary $\Gamma\subset\partial\Omega$, where Dirichlet boundary conditions for the displacement $\pu$ are prescribed.
	If $\pP$ is fully clamped at $\Gamma$, then existence and uniqueness can be obtained even in the most degenerate case for zero Cosserat couple modulus $\muc=0$ and $\mum=\kam=0$.
	On the contrary, if the new \underline{consistent coupling boundary condition} is adopted, i.e., $\left.D\pu\times\pnu\right|_{\Gamma}=\left.\pP\times\pnu\right|_{\Gamma}$,	then existence and uniqueness follow in a highly non-trivial way by	using recently established incompatible Korn's inequalities \cite{neff2015poincare,lewintan2021p,lewintan2020korn,lewintan2019ne}.
	We have	shown that no other simple argument allows to bypass the use of these new inequalities. In the latter case, the suitable assumptions on the relevant material parameters are $\muc\geqslant0$ and $\mum,\kam>0.$ We remark that if $\muc>0$, $\mum,\kam>0$, then existence and uniqueness would even hold true with only Neumann boundary conditions on $\partial\Omega$ regarding $\pP$, i.e., $\pP$ is completely left free, inducing $\left.\gbm\!\cdot\!\pnu\right|_{\partial\Omega}=0$	as higher stress boundary condition.
	The consistent coupling boundary constraint has already been used in the relaxed micromorphic model, where only $\Curl\pP$ appears as curvature measure, allowing only to prescribe $\left.\pP\times\pnu\right|_{\Gamma}$ arbitrarily.
	Using this constraint also on $\pP$ in the micromorphic	model, in which $\text{D}\pP$ appears as curvature measure, escaped our
	attention for some time.
	We have presented again some analytical solutions for the respective boundary conditions which allow to see that the latter may trigger stiffer (\textit{clamped constraint}) or a softer (\textit{consistent coupling conditions}) response of the samples.
	As expected, the fully clamped boundary conditions predicts stiffer response.

	{\scriptsize
	\paragraph{{\scriptsize Acknowledgements.}}
	Angela Madeo acknowledges support from the European Commission through the funding of the ERC Consolidator Grant META-LEGO, N 101001759.
	Angela Madeo and Gianluca Rizzi acknowledge funding from the French Research Agency ANR, “METASMART” (ANR-17CE08-0006).
	Angela Madeo and Gianluca Rizzi acknowledge support from IDEXLYON in the framework of the “Programme Investissement d’Avenir” ANR-16-IDEX-0005.
	Patrizio Neff acknowledges support in the framework of the DFG-Priority Programme 2256 ``Variational Methods for Predicting Complex Phenomena in Engineering Structures and Materials'', funded by the Deutsche Forschungsgemeinschaft (DFG, German research foundation), Project-ID 422730790, and a collaboration of projects ‘Mathematical analysis of microstructure in supercompatible alloys’ (Project-ID 441211072) and ‘A variational scale-dependent transition scheme - from Cauchy elasticity to the relaxed micromorphic continuum’ (Project-ID 440935806). Peter Lewintan and Patrizio Neff were supported by the Deutsche Forschungsgemeinschaft (Project-ID 415894848).  Hassam Khan acknowledges the support of the German Academic Exchange Service (DAAD) and the Higher Education Commission of Pakistan (HEC).
}
\printbibliography
	{\footnotesize
	\section{Appendix}
	\subsection{The Lie-algebra $\so$ and the maps Anti and axl}
	
	Let us introduce the Lie-algebra
	\begin{equation}
	\gs\go(3)\defi\left\{ \begin{pmatrix}0 & -v_{3} & v_{2}\\
	v_{3} & 0 & -v_{1}\\
	-v_{2} & v_{1} & 0
	\end{pmatrix}\in\bR^{3\times3}\;\left. \right| \;\bov\in\bR^3\right\} 
	\end{equation}
	equipped with the matrix commutator bracket 
	$
	\left[A,B\right]=A\!\cdot\!\!B-B\!\cdot\!\!A.
	$
	We identify it with the Lie-algebra $(\bR^{3},\times)$ via the isomorphism
	\begin{equation}
	\axl:\gs\go(3)\longrightarrow\bR^{3},\qquad\begin{pmatrix}0 & -v_{3} & v_{2}\\
	v_{3} & 0 & -v_{1}\\
	-v_{2} & v_{1} & 0
	\end{pmatrix}\mapsto\begin{pmatrix}v_{1}\\
	v_{2}\\
	v_{3}
	\end{pmatrix}
	\end{equation}
	whose inverse is
	\begin{equation}
	\Anti:\bR^{3}\longrightarrow\gs\go(3),\qquad\begin{pmatrix}v_{1}\\
	v_{2}\\
	v_{3}
	\end{pmatrix}\mapsto\Anti\begin{pmatrix}v_{1}\\
	v_{2}\\
	v_{3}
	\end{pmatrix}\defi\begin{pmatrix}0 & -v_{3} & v_{2}\\
	v_{3} & 0 & -v_{1}\\
	-v_{2} & v_{1} & 0
	\end{pmatrix}.
	\end{equation}
	These isomorphisms are constructed in a such a way that 
	\begin{equation}\label{product convertion}
	\left(\Anti\boldsymbol{u}\right)\!\cdot\! \boldsymbol{v}=\boldsymbol{u}\times \boldsymbol{v}\qquad\forall(\boldsymbol{u},\boldsymbol{v})\in\bR^{3}\times\bR^{3}.
	\end{equation}
	Moreover,
	\begin{equation}
	(\Anti\pu)\!\cdot\!\pv=\pu\times\pv=-\pv\times\pu=-(\Anti\pv)\!\cdot\!\pu,\qquad\forall\pu,\pv\in\bR^{3}.
	\end{equation}
	We remember also the following useful relations: let $\boldsymbol{A}\in\so$ and $\boldsymbol{a}=\axl\boldsymbol{A}$, then inductively it can be proved that 
	\begin{equation}
	\boldsymbol{A}^{2n}=\left(-1\right)^{n-1}\left\Vert \boldsymbol{a}\right\Vert _{\bR^{3}}^{2n-2}\boldsymbol{A}^{2}\qquad\textrm{and}\qquad \boldsymbol{A}^{2n-1}=\left(-1\right)^{n-1}\left\Vert \boldsymbol{a}\right\Vert _{\bR^{3}}^{2n-2}\boldsymbol{A}\qquad\forall n\in\bN, 
	\end{equation}
	i.e.,
	\begin{equation}\label{eq:exponentAnti}
	\left(\Anti\boldsymbol{a}\right)^{2n}=\left(-1\right)^{n-1}\left\Vert \boldsymbol{a}\right\Vert _{\bR^{3}}^{2n-2}\left(\Anti\boldsymbol{a}\right)^{2}\qquad\textrm{and}\qquad\left(\Anti\boldsymbol{a}\right)^{2n-1}=\left(-1\right)^{n-1}\left\Vert \boldsymbol{a}\right\Vert _{\bR^{3}}^{2n-2}\Anti\boldsymbol{a}
	\end{equation}
	for all $n\in\bN$ and for all $\boldsymbol{a}\in\bR^{3}$.
	For $n=1$ we get
	\begin{align}\label{power}
	\boldsymbol{A}^{2}\!\cdot\! \boldsymbol{u}=\left(\Anti\boldsymbol{a}\right)^{2}\!\cdot\! \boldsymbol{u}=\boldsymbol{a}\times\left(\boldsymbol{a}\times \boldsymbol{u}\right)\overset{\eqref{cross product}}{=}\left\langle \boldsymbol{a},\boldsymbol{u}\right\rangle \boldsymbol{a}-\left\Vert \boldsymbol{a}\right\Vert ^{2}\boldsymbol{u}=\left(\boldsymbol{a}\otimes \boldsymbol{a}-\left\Vert \boldsymbol{a}\right\Vert ^{2}\id\right)\!\cdot\! \boldsymbol{u}
	\end{align}
	for all $(\boA,\bou)\in\so\times\bR^3$, hence for $n=2$ we obtain
	\begin{align}\label{lie prop}
	\boldsymbol{A}^{3}\!\cdot\!\boldsymbol{u} 
	& =\boldsymbol{A}\!\cdot\!\left(\boldsymbol{A}^{2}\!\cdot\!\boldsymbol{u}\right)\overset{\eqref{power}}{=}\boldsymbol{A}\!\cdot\!\left(\left\langle \boldsymbol{a},\boldsymbol{u}\right\rangle \boldsymbol{a}-\left\Vert \boldsymbol{a}\right\Vert ^{2}\boldsymbol{u}\right)
	=\left\langle \boldsymbol{a},\boldsymbol{u}\right\rangle \underbrace{\boldsymbol{A}\!\cdot\!\boldsymbol{a}}_{=\,0}-\left\Vert \boldsymbol{a}\right\Vert ^{2}\boldsymbol{A}\!\cdot\!\boldsymbol{u}
	=-\left\Vert \boldsymbol{a}\right\Vert ^{2}\boldsymbol{A}\!\cdot\!\boldsymbol{u},
	\end{align}
	because
	\begin{equation}\label{rule0}
	\boldsymbol{A}\!\cdot\!\boldsymbol{a}=\left(\Anti\boldsymbol{a}\right)\!\cdot\!\boldsymbol{a}\overset{\eqref{product convertion}}{=}\boldsymbol{a}\times\boldsymbol{a}=0\qquad \forall \boldsymbol{a}\in\bR^3.
	\end{equation}
	Moreover, considering a second order tensor $\bom$ and a vector $\pnu$, we have
	\begin{align}
	\left(\bom\times\pnu\right)_{ij}&=\left(\left(\bom\right)_{i}\times\pnu\right){}_{j}=\left(-\pnu\times\left(\bom\right)_{i}\right){}_{j}=\left(-(\Anti\pnu)\!\cdot\!\left(\bom\right)_{i}\right){}_{j}=-\left(\Anti\pnu\right)_{jk}\left(\left(\bom\right)_{i}\right)_{k}\nonumber
	\\
	&=-\left(\Anti\pnu\right)_{jk}\bom_{ik}=-\bom_{ik}\left(\Anti\pnu\right)_{jk}=\bom_{ik}\left(\Anti\pnu\right)_{kj} \nonumber
	\\
	&=(\bom\!\cdot\!\Anti\pnu)_{ij}\qquad\forall i,j\in\left\lbrace 1,2,3\right\rbrace .\label{new id}
	\end{align}
	\subsection{Derivation of consistent coupling mixed boundary conditions in the micromorphic model}\label{Derivation new bc}
	
	In this section, for the convenience of the reader, we give the details
	of the first variation of the curvature part of the accounted action
	functional over the space 
	$
	\mathscr{H}^{\,\sharp}(\Omega)\defi\left\{ \pP\in H^{1}(\Omega,\bR^{3\times3})\;\left. \right| \;\left.\pP\times\pnu\right|_{\Gamma}=0\right\}, 
	$
	i.e.,
	\begin{equation}
		\delta\sF_{\text{curv}}\left[\pP,\delta\pP\right]=\left.\mathrm{\frac{d}{dt}}\right|_{t=0}\int_{\Omega}W_{\text{curv}}\left(\text{D}\boldsymbol{P}+t\,\text{D}\delta\boldsymbol{P}\right)\text{dV}=0.
	\end{equation}
	\subsubsection{Normal and tangential decomposition of $\boldsymbol{\mathfrak{m}}\!\cdot\!\boldsymbol{\nu}$}
	
	Starting from the classical triple product relation
	\begin{equation}\label{cross product}
		\boldsymbol{a}\times\left(\boldsymbol{b}\times\boldsymbol{c}\right)=\boldsymbol{b}\left\langle \boldsymbol{a},\boldsymbol{c}\right\rangle -\boldsymbol{c}\left\langle \boldsymbol{a},\boldsymbol{b}\right\rangle ,\qquad\textrm{for}\qquad\boldsymbol{a},\boldsymbol{b},\boldsymbol{c}\in\bR^{3},
	\end{equation}
	let us consider two vectors $\pu,\pnu\in\bR^{3}$ with $\left\| \pnu\right\|=1 $. Then,
	\begin{equation}
		\left\{ \begin{aligned}\pu & =\pu-\left\langle \pu,\pnu\right\rangle \pnu+\left\langle \pu,\pnu\right\rangle \pnu=\underbrace{\overbrace{\left\langle \pnu,\pnu\right\rangle }^{=\,1}\pu-\left\langle \pu,\pnu\right\rangle \pnu}_{=\,\pnu\times\left(\pu\times\pnu\right)}+\underbrace{\left\langle \pu,\pnu\right\rangle \pnu}_{=\,\left(\pnu\otimes\pnu\right)\!\cdot\!\pu}\\
			\pu & =\left(\id-\pnu\otimes\pnu\right)\!\cdot\!\pu+\left(\pnu\otimes\pnu\right)\!\cdot\!\pu
		\end{aligned}
		\right.\qquad\Longrightarrow\qquad\pnu\times\left(\pu\times\pnu\right)=\left(\id-\pnu\otimes\pnu\right)\!\cdot\!\pu.
	\end{equation}
	We want to give the equivalent formulation in terms of $\Anti\pnu$ for $\pnu\in\bR^{3}$ with $\left\| \pnu\right\|=1 $:
		\begin{align}\label{antirel}
		\pnu\times\left(\pu\times\pnu\right) & =-\pnu\times\left(\pnu\times \pu\right) = -\pnu\times\left((\Anti\pnu)\!\cdot\!\pu\right)=-\Anti\pnu\!\cdot\!\left((\Anti\pnu)\!\cdot\!\pu\right) \\
		& =\left(\Anti\pnu\right)^{T}\!\cdot\!\left((\Anti\pnu)\!\cdot\!\pu\right) =\underbrace{\left(\left(\Anti\pnu\right)^{T}\!\cdot\!\Anti\pnu\right)}_{\underset{\left\| \pnu\right\|=1}{=}\,\id-\pnu\otimes\pnu}\!\cdot\!\,\pu=-\left(\Anti\pnu\right)^{2}\!\cdot\!\pu\,. \nonumber
		\end{align}
    We show now that 
	\begin{equation}
		\left[\pnu\times\left(\pu\times\pnu\right)\right]\times\pnu=\pu\times\pnu\label{eq:identity rotation}
	\end{equation}
	using the operator $\Anti$. Indeed
	\begin{align}
		\left[\pnu\times\left(\pu\times\pnu\right)\right]\times\pnu&=-\pnu\times\left[\pnu\times\left(\pu\times\pnu\right)\right]\overset{\eqref{antirel}}{=}\pnu\times \left[\left(\Anti\pnu\right)^{2}\!\cdot\!\pu\right] = (\Anti\pnu)\!\cdot\!\left[\left(\Anti\pnu\right)^{2}\!\cdot\!\pu\right] \notag\\
		&= \left(\Anti\pnu\right)^{3}\!\cdot\!\pu \overset{\eqref{eq:exponentAnti}}{=}-\left\| \pnu\right\|^2(\Anti\pnu)\!\cdot\!\pu = -(\Anti\pnu)\!\cdot\!\pu=-\pnu\times\pu=\pu\times\pnu\,.
		\end{align}
	The identity (\ref{eq:identity rotation}) can be obtained also directly
	as follows
	
	\begin{align}
		\left[\pnu\times\left(\pu\times\pnu\right)\right]\times\pnu & =\left[\pu-\left\langle \pu,\pnu\right\rangle \pnu\right]\times\pnu=\pu\times\pnu-\left\langle \pu,\pnu\right\rangle \underbrace{\pnu\times\pnu}_{=\,0}=\pu\times\pnu.\label{eq:}\\
		\nonumber 
	\end{align}
	Thus, considering $\pP\in L^{2}(\Gamma,\bR^{3\times3})$,
	since $\pnu\in L^{\infty}(\Gamma,\bR^{3})$, we have the $L^2-$
	decomposition
	\begin{align}\label{eq:final}
		\pP & =\pP\!\cdot\!\left(\id-\pnu\otimes\pnu\right)+\pP\!\cdot\!\left(\pnu\otimes\pnu\right)\ =\pP\!\cdot\!\left(\left(\Anti\pnu\right)^{T}\!\cdot\!\Anti\pnu\right)+\pP\!\cdot\!\left(\id-\left(\Anti\pnu\right)^{T}\!\cdot\!\Anti\pnu\right)\nonumber \\
		& =-\,\pP\!\cdot\!\left(\Anti\pnu\right)^{2}+\pP\!\cdot\!\left(\id+\left(\Anti\pnu\right)^{2}\right).
	\end{align}
	\subsubsection{Variation}
	
	We have
	
	\begin{align}
		\delta\sF_{\text{curv}}\left[\pP,\delta\pP\right]&=\left.\mathrm{\frac{d}{dt}}\right|_{t=0}\int_{\Omega}W_{\text{curv}} \left(\text{D}\boldsymbol{P}+t\,\text{D}\delta\boldsymbol{P}\right)\text{dV}\nonumber \\*
		& =\int_{\Omega}\,\mu\,L_{\text{c}}^{2}\,\Bigg(a_{1}\,\langle\text{D}\!\left(\text{dev}\,\text{sym}\,\boldsymbol{P}\right),\text{D}\,\delta\boldsymbol{P}\rangle+a_{2}\,\langle\text{D}\!\left(\text{skew}\,\boldsymbol{P}\right),\text{D}\,\delta\boldsymbol{P}\rangle+\frac{2}{9}\,a_{3}\langle\text{D}\!\left(\text{tr}\!\left(\boldsymbol{P}\right)\id\right),\text{D}\,\delta\boldsymbol{P}\rangle\Bigg)\,\text{dV}\nonumber\\
		& =\int_{\Omega}\,\mu\,L_{\text{c}}^{2}\,\Bigg(\langle a_{1}\,\text{D}\left(\text{dev}\,\text{sym}\,\boldsymbol{P}\right)+a_{2}\,\text{D}\!\left(\text{skew}\,\boldsymbol{P}\right)+\frac{2}{9}\,a_{3}\text{D}\!\left(\text{tr}\!\left(\boldsymbol{P}\right)\id\right),\text{D}\,\delta\boldsymbol{P}\rangle\Bigg)\,\text{dV}\,.\label{eq:first_varia_MM}
	\end{align}
	The product rule implies that 
	\begin{align}
		& \text{div}\left[\left(\,\mu\,L_{\text{c}}^{2}\,a_{1}\,\text{D}\!\left(\text{dev}\,\text{sym}\,\boldsymbol{P}\right)+a_{2}\,\text{D}\!\left(\text{skew}\,\boldsymbol{P}\right)+\frac{2}{9}\,a_{3}\text{D}\!\left(\text{tr}\!\left(\boldsymbol{P}\right)\id\right)\right)\! :\!\delta\boldsymbol{P}\right]\nonumber \\*
		& \hspace{3cm}=\,\mu\,L_{\text{c}}^{2}\,\langle\text{DIV}\left[a_{1}\,\text{D}\left(\text{dev}\,\text{sym}\,\boldsymbol{P}\right)+a_{2}\,\text{D}\!\left(\text{skew}\,\boldsymbol{P}\right)+\frac{2}{9}\,a_{3}\text{D}\!\left(\text{tr}\!\left(\boldsymbol{P}\right)\id\right)\right],\delta\boldsymbol{P}\rangle_{\bR^{3\times3}}\label{eq:first_varia_MM_2}\\*
		& \hspace{4cm}+\,\mu\,L_{\text{c}}^{2}\,\left(\langle a_{1}\,\text{D}\!\left(\text{dev}\,\text{sym}\,\boldsymbol{P}\right)+a_{2}\,\text{D}\!\left(\text{skew}\,\boldsymbol{P}\right)+\frac{2}{9}\,a_{3}\text{D}\!\left(\text{tr}\!\left(\boldsymbol{P}\right)\id\right),\text{D}\,\delta\boldsymbol{P}\rangle_{\bR^{3\times3\times3}}\right)\,.\nonumber 
	\end{align}
	Thus, remembering the definition of $\gbm\in\bR^{3\times3\times3}$ given in \eqref{eq:m_third_MM}, we obtain
	\begin{align}
		\delta\sF_{\text{curv}}\left[\pP,\delta\pP\right] & =\int_{\Omega}\text{div}\left[\gbm\! : \!\delta\pP\right]\,\text{dV}-\int_{\Omega}\langle\DIV\gbm,\delta\pP\rangle_{\bR^{3\times3}}\,\text{dV}.
	\end{align}
	The therm $\int_{\Omega}\langle\DIV\gbm,\delta\pP\rangle\,\text{dV}$
	will contribute a term in the bulk equilibrium equations, while thanks
	to the divergence theorem the term $\int_{\Omega}\text{div}\left[\gbm\! : \!\delta\pP\right]\,\text{dV}$
	can be written as 
	\begin{align}
		\int_{\Omega}\text{div}\left[\gbm\! : \!\delta\pP\right]\,\text{dV}=\int_{\partial\Omega}\langle\gbm\!\cdot\!\boldsymbol{\nu},\delta\pP\rangle_{\bR^{3\times3}}\,\text{ds}\,.\label{eq:first_varia_MM_4}
	\end{align}
		Recall that
		\begin{equation}\label{pat}
		\left.\pP\times\pnu\right|_{\Gamma}=0\quad\Longrightarrow\quad\left.\delta\pP\times\pnu\right|_{\Gamma}=0\quad\Longrightarrow\quad\left.\pnu\times\left(\delta\pP\times\pnu\right)\right|_{\Gamma}=\left.\delta\pP\!\cdot\!\left(\id-\pnu\otimes\pnu\right)\right|_{\Gamma}=0,
		\end{equation}
		so that it holds
		\begin{equation}
		\delta \pP=\delta \pP \!\cdot\! (\nu\otimes\nu) \quad\text{ along $\Gamma$.}
		\end{equation}
		Hence, for all tensor fields $\boB$ we have
		\begin{align}\label{eq:ast}
		\!\!\skalarProd{\boB}{\delta \pP}_{\bR^{3\times3}}&=\skalarProd{\boB}{\delta \pP \!\cdot\! (\pnu\otimes\pnu)}_{\bR^{3\times3}}\!\!\!\overset{\|\pnu\|=1}{=}\skalarProd{\boB}{\delta \pP \!\cdot\! (\pnu\otimes\pnu)^2}_{\bR^{3\times3}} \!\!= \skalarProd{\boB\!\cdot\!(\pnu\otimes \pnu)}{\delta \pP \!\cdot\! (\pnu\otimes\pnu)}_{\bR^{3\times3}}
		\end{align}
		along $\Gamma$.
		Returning to the right-hand side of \eqref{eq:first_varia_MM_4} we conclude:
		\begin{align}
		\int_{\partial\Omega}\skalarProd{\boldsymbol{\mathfrak{m}}\!\cdot\!\pnu}{\delta \pP}_{\bR^{3\times3}}\,\text{ds}& = \int_{\Gamma}\skalarProd{\boldsymbol{\mathfrak{m}}\!\cdot\!\pnu}{\delta \pP}_{\bR^{3\times3}}\,\text{ds} + \int_{\partial\Omega\backslash\overline{\Gamma}}\skalarProd{\boldsymbol{\mathfrak{m}}\!\cdot\!\pnu}{\delta \pP}_{\bR^{3\times3}}\,\text{ds} \notag\\
		& \!\!\!\overset{\eqref{eq:ast}}{=}\!\!\!
		\int_{\Gamma}\bigl\langle\left(\gbm\!\cdot\!\pnu\right)\!\cdot\!\left(\pnu\otimes\pnu\right),\underbrace{\delta\pP\!\cdot\!(\pnu\otimes\pnu)}_{\text{free}}\bigr\rangle_{\bR^{3\times3}}\text{ds}+\int_{\partial\Omega\setminus\overline{\Gamma}}\langle\gbm\!\cdot\!\pnu,\underbrace{\delta\pP}_{\text{free}}\rangle_{\bR^{3\times3}}\,\text{ds}.
		\end{align}
		Since  $\pnu\otimes\pnu=\left(\id-(\Anti\pnu)^{T}\!\cdot\!\Anti\pnu\right)=\left(\id+(\Anti\pnu)^{2}\right)$
		we can also rewrite the first term on the right-hand side of the last equation
		in terms of the $\Anti$ operator:
		\begin{align}
		\int_{\Gamma}\bigl\langle\left(\gbm\!\cdot\!\pnu\right)\!\cdot\!\left(\pnu\otimes\pnu\right),\underbrace{\delta\pP\!\cdot\!(\pnu\otimes\pnu)}_{\text{free}}\bigr\rangle_{\bR^{3\times3}}\text{ds} = \int_{\Gamma}\langle\left(\gbm\!\cdot\!\pnu\right)\!\cdot\!\left(\id+(\text{Anti}\,\pnu)^2\right),\underbrace{\delta\pP\!\cdot\!\left(\id+(\text{Anti}\,\pnu)^2\right)}_{\text{free}}\,\rangle_{\bR^{3\times3}}\,\text{ds}.
		\end{align}
	Note that 
	\begin{align}
		\left(\id-\pnu\otimes\pnu\right)\!\cdot\!\left(\pnu\otimes\pnu\right)&=\boldsymbol{0}\,, & \textrm{but}& & \left(\id-\text{Anti}\,\pnu\right)\!\cdot\!\text{Anti}\,\pnu&\neq\boldsymbol{0}.
	\end{align}
	Thus 
	\begin{equation}
	\delta\sF_{\,\text{curv}}\left[\pP,\delta\pP\right]=\int_{\Gamma}\bigl\langle\left(\gbm\!\cdot\!\pnu\right)\!\cdot\!\left(\pnu\otimes\pnu\right),\underbrace{\delta\pP\!\cdot\!(\pnu\otimes\pnu)}_{\text{free}}\bigr\rangle_{\bR^{3\times3}}\text{ds}+\int_{\partial\Omega\setminus\overline{\Gamma}}\langle\gbm\!\cdot\!\pnu,\underbrace{\delta\pP}_{\text{free}}\rangle_{\bR^{3\times3}}\,\text{ds}
	-\int_{\Omega}\langle\DIV\gbm,\delta\pP\rangle_{\bR^{3\times3}}\,\text{dV}.
	\end{equation}
	The variations $\delta\pP$ in the bulk and at the boundary are independent, therefore $\delta\sF_{\,\text{curv}}\left[\pP,\delta\pP\right]=0$ implies
	\begin{equation}
	\left.\left(\gbm\!\cdot\!\pnu\right)\!\cdot\!\left(\pnu\otimes\pnu\right)\,\right|_{\boldsymbol{\Gamma}} =0\qquad\textrm{and}\qquad\left.\gbm\!\cdot\!\pnu\,\right|_{\partial\Omega\setminus\overline{\Gamma}}.\nonumber 
	\end{equation}
	Hence summarizing, the natural boundary conditions for the mixed problem	for the classical micromorphic model  are
	\[
	\begin{array}{cccccccc}
		\textrm{Dirichlet} &  &  &  &  &  &  & \textrm{Neumann}
		\\
		\\
		\left.\pP\times\pnu\right|_{\Gamma}=0, &  &  &  &  &  &  & \left.(\gbm\ \!\cdot\! \pnu)\!\cdot\!\left(\pnu\otimes\pnu\right)\right|_{\Gamma}=0
		\overset{\eqref{conto second gradient}}{\Leftrightarrow}
		\left.(\gbm \!\cdot\! \boldsymbol{\nu})\!\cdot\! \pnu\right|_{\Gamma}=0,
		\\
		\\
		&  &  &  &  &  &  & \left.\gbm\!\cdot\!\pnu\,\right|_{\partial\Omega\setminus\overline{\Gamma}}=0.
	\end{array}
	\]
	\subsection{Boundary conditions in the relaxed micromorphic model}
	
	For the convenience of the reader we repeat the above reasoning for the relaxed micromorphic model.
	We know that
	\[
	\div(\bou\times \bov)=\left\langle \curl \bou,\bov\right\rangle -\left\langle \curl \bov,\bou\right\rangle \, ,
	\]
	and with
	\begin{equation}
		\boldsymbol{m}\defi\mu\frac{L^{2}}{2}\,\bL\,\Curl\pP\in\bR^{3\times3} \, ,
	\end{equation}
	we have
	\begin{align}
		\!\!\!\!\mu\frac{L^{2}}{2}\sum_{i=1}^{3}\left\langle \left(\bL\,\Curl\pP\right)_{i},\left(\Curl\delta\pP\right)_{i}\right\rangle _{\bR^{3}} & =\sum_{i=1}^{3}\left\langle \left(\boldsymbol{m}\right)_{i},\curl\!\left(\delta\pP\right)_{i}\right\rangle _{\bR^{3}}
		=\sum_{i=1}^{3}\left(-\,\div[\left(\boldsymbol{m}\right)_{i}\times\left(\delta\pP\right)_{i}]+\left\langle \curl\!\left(\boldsymbol{m}\right)_{i},\left(\delta\pP\right)_{i}\right\rangle _{\bR^{3}}\right).
	\end{align}
	Thus,
	\begin{align}
		\delta\sF_{\,\text{curv}}^{\text{relax}}\left[\pP,\delta\pP\right] &=\mu\frac{L^{2}}{2}\int_\Omega \left\langle \bL\,\Curl\pP,\Curl\delta\pP\right\rangle _{\bR^{3\times3}} \,\text{dV} 
		= \mu\frac{L^{2}}{2}\sum_{i=1}^{3}\int_\Omega\left\langle \left(\bL\,\Curl\pP\right)_{i},\left(\Curl\delta\pP\right)_{i}\right\rangle _{\bR^{3}} \,\text{dV} \nonumber
		\\
		&=\sum_{i=1}^{3}  \left(\int_{\Omega}-\,\div[\left(\boldsymbol{m}\right)_{i}\times\left(\delta\pP\right)_{i}]\,\text{dV}+\int_{\Omega}\left\langle \curl\!\left(\boldsymbol{m}\right)_{i},\left(\delta\pP\right)_{i}\right\rangle _{\bR^{3}}\text{dV}\right)\nonumber
		\\
		& =\sum_{i=1}^{3}\int_{\Omega}-\,\div[\left(\boldsymbol{m}\right)_{i}\times\left(\delta\pP\right)_{i}]\,\text{dV}+\int_{\Omega}\left\langle \Curl\boldsymbol{m},\delta\pP\right\rangle _{\bR^{3\times3}}\text{dV} 
	\end{align}
	and
	\begin{align}
	-\sum_{i=1}^{3}\int_{\Omega}\div & \left[\left(\boldsymbol{m}\right)_{i}\times\left(\delta\pP\right)_{i}\right]\text{dV}\nonumber
	\\*
	& =-\sum_{i=1}^{3}\int_{\partial\Omega}\left\langle \left(\boldsymbol{m}\right)_{i}\times\left(\delta\pP\right)_{i},\pnu\right\rangle _{\bR^{3}}\text{ds}=\sum_{i=1}^{3}\int_{\partial\Omega}\left\langle \left(\boldsymbol{m}\right)_{i}\times\pnu,\left(\delta\pP\right)_{i}\right\rangle _{\bR^{3}}\text{ds}\nonumber
	\\
	& =\int_{\partial\Omega}\left\langle \boldsymbol{m}\times\pnu,\delta\pP\right\rangle _{\bR^{3\times3}}\text{ds}\overset{\eqref{new id}}{=}\int_{\partial\Omega}\left\langle \boldsymbol{m}\!\cdot\!\Anti\pnu,\delta\pP\right\rangle _{\bR^{3\times3}}\text{ds}=\int_{\Gamma\cup\partial\Omega\setminus\overline{\Gamma}}\left\langle \boldsymbol{m}\!\cdot\!\Anti\pnu,\delta\pP\right\rangle _{\bR^{3\times3}}\text{ds} \nonumber
	\\
	&=\int_{\Gamma}\left\langle \boldsymbol{m}\!\cdot\!\Anti\pnu,\delta\pP\right\rangle _{\bR^{3\times3}}\text{ds} + \int_{\partial\Omega\setminus\overline{\Gamma}}\left\langle \boldsymbol{m}\!\cdot\!\Anti\pnu,\delta\pP\right\rangle _{\bR^{3\times3}}\text{ds}\notag
	\\
	&\underset{\eqref{power}}{\overset{\mathclap{\eqref{eq:ast}}}{=}} \underbrace{\int_{\Gamma}\left\langle \boldsymbol{m}\!\cdot\!\Anti\pnu,\delta\pP\!\cdot\!\left(\pnu\otimes\pnu\right)\right\rangle _{\bR^{3\times3}}\text{ds}}_{=\,0,\;\textrm{because}\;\boldsymbol{m}\!\cdot\!\Anti\pnu\,\bot\,\delta\pP\!\cdot\!\left(\pnu\otimes\pnu\right)}
		+ \int_{\partial\Omega\setminus\overline{\Gamma}}\left\langle \boldsymbol{m}\!\cdot\!(-\Anti\pnu)^3,\delta\pP\right\rangle _{\bR^{3\times3}}\text{ds}
	\\
	&= \int_{\partial\Omega\setminus\overline{\Gamma}}\left\langle \boldsymbol{m}\!\cdot\!\Anti\pnu,\delta\pP\!\cdot\!(-\Anti\pnu)^2\right\rangle _{\bR^{3\times3}}\text{ds}
	\,= \int_{\partial\Omega\setminus\overline{\Gamma}}\bigl\langle\boldsymbol{m}\!\cdot\!\Anti\pnu,\underbrace{\delta\pP\!\cdot\!\left(\id-\pnu\otimes\pnu\right)}_{\textrm{free}}\bigr\rangle_{\bR^{3\times3}}\text{ds}\,,\notag
	\end{align}
	where in the last step we have again used that $(-\Anti\pnu)^2=\id-\pnu\otimes\pnu$. Thus 
	\begin{equation}
	\delta\sF_{\,\text{curv}}^{\text{relax}}\left[\pP,\delta\pP\right]=\int_{\partial\Omega\setminus\overline{\Gamma}}\left\langle \boldsymbol{m}\!\cdot\!\Anti\pnu,\delta\pP\!\cdot\!\left(\id-\pnu\otimes\pnu\right)\right\rangle _{\bR^{3\times3}}\text{ds}
	+\int_{\Omega}\left\langle \Curl\boldsymbol{m},\delta\pP\right\rangle _{\bR^{3\times3}}\text{dV}.
	\end{equation}
	The variations $\delta\pP$ in the bulk and at the boundary are independent, therefore $\delta\sF_{\,\text{curv}}^{\text{relax}}\left[\pP,\delta\pP\right]=0$ implies
	\begin{equation}
		\left.\boldsymbol{m}\!\cdot\!\Anti\pnu\right|_{\partial\Omega\setminus\overline{\Gamma}}\overset{\eqref{new id}}{=}\left.\boldsymbol{m}\times\pnu\right|_{\partial\Omega\setminus\overline{\Gamma}}=0.
	\end{equation}
	Hence summarizing, the natural boundary conditions for the mixed problem
	for the relaxed micromorphic model  are
	\[
	\begin{array}{cccccccc}
		\textrm{Dirichlet} &  &  &  &  &  &  & \textrm{Neumann}\\
		\\
		\left.\pP\times\pnu\right|_{\Gamma}=0 &  &  &  &  &  &  & \left.\boldsymbol{m}\times\pnu\right|_{\partial\Omega\setminus\overline{\Gamma}}=0
	\end{array}
	\]
	\subsection{Derivation of mixed boundary conditions for a linear second gradient material}
	
	The second gradient linearized elasticity model (without mixed terms) reads (see for example \cite{madeo2016new,ghiba2017variant,munch2017modified,eremeyev2018note,eremeyev2019characterization,eremeyev2021nonlinear,rahali2020surface,eremeyev2021weak,eremeyev2020well,eremeyev2020weak})
	\begin{align}
		\sF:\bH & \fr\bR^{+}, & \pu\mapsto \sF\left[\pu\right] & =\underbrace{\int_{\Omega}\frac{1}{2}\left\langle \bC\,\sym\text{D}\pu,\sym\text{D}\pu\right\rangle \text{dV}}_{(\textrm{I})}+\underbrace{\int_{\Omega}\frac{1}{2}\left\langle \bG\,\text{D}^{2}\pu,\text{D}^{2}\pu\right\rangle \text{dV}}_{(\textrm{II})},
	\end{align}
	where
	\begin{align*}
		\bH\defi\left\{ u\in H_{0,\Gamma}^{1}\left(\Omega,\bR^{3}\right)\cap H^{\,2}\!\left(\Omega,\bR^{3}\right)\;\left.\right|\;\bG\,\text{D}^{2}\pu\in H\!\left(\Div\!\left(\DIV\right);\Omega,\bR^{3}\otimes\Sym\!(3)\right)\right\} 
	\end{align*}
	and
	\begin{align}
		H\!\left(\Div(\DIV\!);\Omega,\bR^{3}\otimes\Sym\!(3)\right) & \defi\left\{ \boldsymbol{U}\in H\!\left(\DIV;\Omega,\bR^{3}\otimes\Sym\!(3)\right)\;\left.\right|\;\Div\DIV\boldsymbol{U}\in L^{2}\left(\Omega,\bR^{3}\right)\right\} \\
		& \;=\left\{ \boldsymbol{U}\in H\!\left(\DIV;\Omega,\bR^{3}\otimes\Sym\!(3)\right)\;\left.\right|\;\DIV\boldsymbol{U}\in H\!\left(\Div;\Omega,\bR^{3\times3}\right)\right\} .\nonumber 
	\end{align}
	In order to evaluate the variational derivative
	$\delta\mathcal{I}\left[\pu,\delta\pu\right]$, with $\delta\pu\in\bH$
	we study the two terms $(\textrm{I})$ and $(\textrm{II})$ separately. To make the calculations easier, we do them assuming the boundary and the involved fields regular. Concerning
	part $\left(\textrm{I}\right)$ we have
	\begin{align*}
		\delta\int_{\Omega}\frac{1}{2} & \left\langle \bC\,\sym\text{D}\pu,\sym\text{D}\pu\right\rangle \text{dV}=-\int_{\Omega}\left\langle \Div\!\left[\bC\,\sym\text{D}\pu\right],\delta\pu\right\rangle \text{dV}+\underbrace{\int_{\partial\Omega\setminus\overline{\Gamma}}\left\langle \left(\bC\,\sym\text{D}\pu\right)\!\cdot\!\pnu,\delta\pu\right\rangle _{\bR^{3}}\text{ds}}_{(a)}
	\end{align*}
	and we obtain the normal stress on $\partial\Omega\setminus\overline{\Gamma}$
	\begin{equation}
		\left.\boldsymbol{\sigma}\!\cdot\!\pnu\,\right|_{\partial\Omega\setminus\overline{\Gamma}}=\left.(\bC\,\sym\text{D}\pu)\!\cdot\!\pnu\,\right|_{\partial\Omega\setminus\overline{\Gamma}}.
	\end{equation}
	Concerning part $\left(\textrm{II}\right)$ we obtain
	\begin{align}
		\delta\int_{\Omega}\frac{1}{2} \left\langle \bG\,\text{D}^{2}\pu,\text{D}^{2}\pu\right\rangle \text{dV}
		&=\int_{\Omega}\left\langle \bG\,\text{D}^{2}\pu,\text{D}^{2}\delta\pu\right\rangle \text{dV}=\int_{\Omega}\left\langle \gbm,\text{D}^{2}\delta\pu\right\rangle \text{dV}
		=-\int_{\Omega}\left\langle \DIV\gbm,\text{D}\delta\pu\right\rangle \text{dV}+\int_{\partial\Omega}\left\langle \gbm\!\cdot\!\pnu,\text{D}\delta\pu\right\rangle _{\bR^{3}}\,\text{ds}\nonumber
		\\
		& =\int_{\Omega}\left\langle \Div\DIV\gbm,\delta\pu\right\rangle \text{dV}
		-\underbrace{\int_{\partial\Omega\setminus\overline{\Gamma}}\left\langle \left(\DIV\gbm\right)\!\cdot\!\pnu,\delta\pu\right\rangle _{\bR^{3}}\text{ds}}_{(b)}
		+\int_{\partial\Omega}\left\langle \gbm\!\cdot\!\pnu,\text{D}\delta\pu\right\rangle _{\bR^{3}}\text{ds}.
	\end{align}
	Using the surface differential operators \cite{madeo2016new}, we can develop further the
	term in duality with $\left.\text{D}\delta \pu\right|_{\partial\Omega}$.
	Indeed, remarking that
	\begin{equation}
		\left.\text{D}\delta \pu\right|_{\partial\Omega}=\left(\left(\id-\pnu\otimes\pnu\right)+\left(\pnu\otimes\pnu\right)\right)\left.\text{D}\delta \pu\right|_{\partial\Omega},
	\end{equation}
	we have
	\begin{align}
		\int_{\partial\Omega}\left\langle \gbm\!\cdot\!\pnu,\text{D}\delta\pu\right\rangle _{\bR^{3\times3}}\text{ds} & =\int_{\partial\Omega}\left\langle \gbm\!\cdot\!\pnu,\text{D}\delta\pu\!\cdot\!\left(\left(\id-\pnu\otimes\pnu\right)+\left(\pnu\otimes\pnu\right)\right)\right\rangle _{\bR^{3\times3}}\text{ds}\nonumber\\
		& =\int_{\partial\Omega}\left\langle \gbm\!\cdot\!\pnu,\text{D}\delta\pu\!\cdot\!\left(\id-\pnu\otimes\pnu\right)\right\rangle _{\bR^{3\times3}}\text{ds}+\int_{\partial\Omega}\left\langle \gbm\!\cdot\!\pnu,\text{D}\delta\pu\!\cdot\!\left(\pnu\otimes\pnu\right)\right\rangle _{\bR^{3\times3}}\text{ds}\nonumber\\
		& =\int_{\partial\Omega}\left\langle -\Div_{\partial\Omega,\top}\left[\gbm\!\cdot\!\pnu\right],\delta\pu\right\rangle _{\bR^{3}}\text{ds}+\int_{\partial\Omega}\left\langle \gbm\!\cdot\!\pnu,\text{D}\delta\pu\!\cdot\!\left(\pnu\otimes\pnu\right)\right\rangle _{\bR^{3\times3}}\text{ds}\\
		& =\underbrace{\int_{\partial\Omega\setminus\overline{\Gamma}}\left\langle -\Div_{\partial\Omega,\top}\left[\gbm\!\cdot\!\pnu\right],\delta\pu\right\rangle _{\bR^{3}}\text{ds}}_{(c)}
		+\int_{\partial\Omega}\left\langle \gbm\!\cdot\!\pnu,\text{D}\delta\pu\!\cdot\!\left(\pnu\otimes\pnu\right)\right\rangle _{\bR^{3\times3}}\text{ds}.\nonumber
	\end{align}
	Now, we develop further the term $\int_{\partial\Omega}\left\langle \gbm\!\cdot\!\pnu,\text{D}\delta\pu\!\cdot\!\left(\pnu\otimes\pnu\right)\right\rangle _{\bR^{3\times3}}\text{ds}$. We have,
	\begin{align} \label{conto second gradient}
		\int_{\partial\Omega}\left\langle \gbm\!\cdot\!\pnu,\text{D}\delta\pu\!\cdot\!\left(\pnu\otimes\pnu\right)\right\rangle _{\bR^{3\times3}}\text{ds}\quad & \overset{\mathclap{\norm{\pnu}=1}}{=}\;\,
		\int_{\partial\Omega}\left\langle \gbm\!\cdot\!\pnu,\text{D}\delta\pu\!\cdot\!\left(\pnu\otimes\pnu\right)^2
		\right\rangle _{\bR^{3\times3}}\text{ds}
		=\int_{\partial\Omega}\left\langle \left(\gbm\!\cdot\!\pnu\right)\!\cdot\!\left(\pnu\otimes\pnu\right),\text{D}\delta\pu\!\cdot\!\left(\pnu\otimes\pnu\right)\right\rangle _{\bR^{3\times3}}\text{ds}\nonumber
		\\
		& =\int_{\partial\Omega}\left\langle \left(\left(\gbm\!\cdot\!\pnu\right)\!\cdot\!\pnu\right)\otimes\pnu,\left(\text{D}\delta\pu\!\cdot\!\pnu\right)\otimes\pnu\right\rangle _{\bR^{3\times3}}\text{ds} \nonumber
		=\int_{\partial\Omega}\left\langle \left(\gbm\!\cdot\!\pnu\right)\!\cdot\!\pnu,\text{D}\delta\pu\!\cdot\!\pnu\right\rangle _{\bR^{3}}\smash{\underbrace{\left\langle \pnu,\pnu\right\rangle _{\bR^{3}}}_{=\,\left\Vert \pnu\right\Vert _{\bR^{3}}^{2}\equiv\,1}}\!\text{ds} \nonumber
		\\
		& =\int_{\partial\Omega}\left\langle \left(\gbm\!\cdot\!\pnu\right)\!\cdot\!\pnu,\text{D}\delta\pu\!\cdot\!\pnu\right\rangle _{\bR^{3}}\text{ds}. 
		\end{align}
	Summing up all contributions (a), (b) and (c), we obtain
	\begin{equation}
		\int_{\partial\Omega\setminus\overline{\Gamma}}\left\langle \left(\bC\,\sym\text{D}\pu\right)\!\cdot\!\pnu-\left(\DIV\gbm\right)\!\cdot\!\pnu-\Div_{\partial\Omega,\top}\left[\gbm\!\cdot\!\pnu\right],\delta\pu\right\rangle _{\bR^{3}}\text{ds}.
	\end{equation}
	Thus 
	\begin{align}
	\delta\sF\!\left[\pu,\delta\pu\right]&=
	\int_{\partial\Omega\setminus\overline{\Gamma}}\left\langle \left(\bC\,\sym\text{D}\pu\right)\!\cdot\!\pnu-\left(\DIV\gbm\right)\!\cdot\!\pnu-\Div_{\partial\Omega,\top}\left[\gbm\!\cdot\!\pnu\right],\delta\pu\right\rangle _{\bR^{3}}\text{ds} \nonumber
	\\
	&\qquad-\int_{\Omega}\left\langle \Div\!\left[\bC\,\sym\text{D}\pu\right],\delta\pu\right\rangle \text{dV}+\int_{\Omega}\left\langle \Div\DIV\gbm,\delta\pu\right\rangle \text{dV} \nonumber
	\\
	&
	\qquad\qquad+\int_{\partial\Omega}\left\langle \left(\gbm\!\cdot\!\pnu\right)\!\cdot\!\pnu,\text{D}\delta\pu\!\cdot\!\pnu\right\rangle _{\bR^{3}}\text{ds}.
	\end{align}
	The variations $\delta\pu$ in the bulk and at the boundary are independent, therefore $\delta\sF\!\left[\pu,\delta\pu\right]=0$ implies
	\begin{subequations}
		\begin{align} 
			\left.\boldsymbol{\eta}\right|_{\partial\Omega\setminus\overline{\Gamma}} & \defi\bigl((\bC\,\sym\text{D}\pu)\!\cdot\!\pnu\vphantom{\Div\gbm\!\cdot\!\pnu}\;-\,(\DIV\gbm)\!\cdot\!\pnu\bigr)\bigr|_{\partial\Omega\setminus\overline{\Gamma}}\;-\;\Div_{\partial\Omega,\top}\left.\left[\gbm\!\cdot\!\pnu\right]\right|_{\partial\Omega\setminus\overline{\Gamma}}\nonumber\\
			& =\left.\left(\boldsymbol{\sigma}\!\cdot\!\pnu-(\DIV\gbm)\!\cdot\!\pnu\right)\right|_{\partial\Omega\setminus\overline{\Gamma}}\left.\,-\;\Div_{\partial\Omega,\top}\left(\gbm\!\cdot\!\pnu\right)\right|_{\partial\Omega\setminus\overline{\Gamma}}\nonumber\\
			& =\left.(\bosigma-\DIV\gbm)\!\cdot\!\pnu\right|_{\partial\Omega\setminus\overline{\Gamma}}\left.\,-\;\Div_{\partial\Omega,\top}\left(\gbm\!\cdot\!\pnu\right)\right|_{\partial\Omega\setminus\overline{\Gamma}}=0,\label{bc second gradient1}
			\\
			\left.\boldsymbol{b}\right|_{\partial\Omega} & \defi\left.\left(\gbm\!\cdot\!\pnu\right)\!\cdot\!\pnu\,\right|_{\partial\Omega}=0 .
			\label{bc second gradient2}
		\end{align}
	\end{subequations}
	\subsection{Relaxed micromorphic model as a particular case of the general micromorphic model}
	
	The relaxed micromorphic model can also be written formally  as the classical micromorphic model since $\Curl\pP$ is related to $\text{D}\pP$ via  a linear map. This means in particular that there exists a linear map 
	\begin{equation}
	\mathcal{L}:\bR^{3\times 3\times 3}\fr\bR^{3\times 3}\qquad\textrm{such that}\qquad  \mathcal{L}\,\text{D}\pP=\Curl\pP,
	\end{equation}
	allowing us to rewrite the curvature term in the energy density of the relaxed micromorphic model as follows
	\begin{equation}
	\left\langle \bL\,\Curl\pP,\Curl\pP\right\rangle =\left\langle \bL\,\mathcal{L}\,\text{D}\pP,\mathcal{L}\text{D}\pP\right\rangle =\bigl\langle\underbrace{\mathcal{L}^{T}\,\bL\,\mathcal{L}}_{=:\,\widehat{\bL}}\,\text{D}\pP,\text{D}\pP\bigr\rangle.
	\end{equation}
	Such a map $\mathcal{L}$ can be given explicitly. Indeed, introducing the Levi-Civita third order tensor $\LC$, the classical $\curl$ operator can be written in terms of the Jacobian as follows
	\begin{equation}
	\curl\pu=\boldsymbol{\epsilon}:\text{D}\pu^T
	\end{equation}
	where the contraction operator $:$ is defined component-wise
	in the following manner:
	\[
	\left(\curl\pu\right)_{k}=\epsilon_{kij}(\text{D}\pu^T)_{ij}=\epsilon_{kij}u_{j,i}\qquad\forall k\in\left\{ 1,2,3\right\} .
	\]
	Hence, from the definition \eqref{Curl definition} we obtain 
	\begin{align}
	\Curl\pP&=\begin{pmatrix}\curl\left(\pP\right)_{1}\\
	\curl\left(\pP\right)_{2}\\
	\curl\left(\pP\right)_{3}
	\end{pmatrix}=\begin{pmatrix}\boldsymbol{\epsilon}:\left(\text{D}\left(\pP\right)_{1}\right)^{T}\\
	\boldsymbol{\epsilon}:\left(\text{D}\left(\pP\right)_{2}\right)^{T}\\
	\boldsymbol{\epsilon}:\left(\text{D}\left(\pP\right)_{3}\right)^{T}
	\end{pmatrix}
	\end{align}
	giving component-wise
	\[
	\left(\Curl\pP\right)_{\alpha\beta}=\epsilon_{\beta ij}P_{\alpha j,i}.
	\]
	
	\begin{rem}
		The positive-definiteness of the bilinear form $\left\langle \bL\,\Curl\pP,\Curl\pP\right\rangle$ in terms of $\Curl\pP$ does not imply the positive definiteness of $\bigl\langle\widehat{\bL}\,\text{D}\pP,\text{D}\pP\bigr\rangle$ in terms of $\text{D}\pP$ but only its positive semi-definiteness. Thus, proving an existence result using the formulation of the relaxed micromorphic model as a classical one, is not straightforward but needs new function spaces and new coercive inequalities \cite{neff2015poincare,lewintan2021p,lewintan2020korn,lewintan2019ne,neff2011canonical}.
	\end{rem}
	Denoting by $\gbm\defi\widehat{\bL}\,\text{D}\pP$ the relative third order stress tensor and remembering that $\bom=\bL\,\Curl\pP$, we would like to understand if the following relation holds
	\begin{equation}
	\gbm\!\cdot\!\pnu=\bom\times\pnu \in \bR^{3\times3}, \qquad\gbm\in\bR^{3\times3\times3},\quad\bom\in\bR^{3\times3}.
	\end{equation} 
	Considering for simplicity the constitutive tensor $\bL$ trivial in the component-wise calculation, we obtain
	
	\medskip
	
	\begin{align}
	\left\langle \smash{\overbrace{\mathcal{L}^{T}\left(\mathcal{L}\,\text{D}\pP\right)}^{=\,\gbm}},\text{D}\pP\right\rangle  
	& 
	=\left\langle \mathcal{L}\,\text{D}\pP,\mathcal{L}\,\text{D}\pP\right\rangle =\left\langle \Curl\pP,\Curl\pP\right\rangle =P_{\alpha j,i}\epsilon_{\beta ij}P_{\alpha n,m}\epsilon_{\beta mn}
	\\
	& 
	=P_{\alpha j,i}\epsilon_{ij\beta}P_{\alpha n,m}\epsilon_{\beta mn}=\left(P_{\alpha j,i}\epsilon_{ij\beta}\epsilon_{\beta mn}\right)P_{\alpha n,m} =\underbrace{\left(P_{\alpha j,i}\epsilon_{ij\beta}\epsilon_{\beta mn}\right)}_{=\,\gbm\,=\,\bom\!\cdot\!\epsilon}P_{\alpha n,m},\nonumber
	\end{align}
	
	\noindent hence
	\begin{equation}
	\left(\gbm\!\cdot\!\pnu\right)_{\alpha\beta}=M_{\alpha\beta\gamma}\nu_{\gamma}=\left(P_{\alpha j,i}\epsilon_{\beta ij}\epsilon_{\beta q\gamma}\right)\nu_{\gamma}=m_{\alpha\beta}\epsilon_{\beta q\gamma}v_{\gamma}=\left(\bom\times\pnu\right)_{\alpha\beta},
	\end{equation}
	i.e.,
	\begin{align}
	\gbm\!\cdot\!\pnu & =\left(\widehat{\bL}\,\text{D}\pP\right)\!\cdot\!\pnu=\left((\mathcal{L}^{T}\,\bL\,\mathcal{L})\,\text{D}\pP\right)\!\cdot\!\pnu=\left((\mathcal{L}^{T}\,\bL)\,(\mathcal{L}\,\text{D}\pP)\right)\!\cdot\!\pnu\\
	& =\big(\mathcal{L}^{T}\,\smash{\underbrace{\bL\,\Curl\pP}_{=\,\bom}}\big)\!\cdot\!\pnu=\left(\mathcal{L}^{T}\bom\right)\!\cdot\!\pnu=\bom\times\pnu,\qquad\textrm{as claimed.} \nonumber
	\end{align}
	We finally also point out that
	\begin{align*}
	\DIV\left(\mathcal{L}^{T}\left(\mathcal{L}\,\text{D}\pP\right)\right) & =\left(\mathcal{L}^{T}\left(\mathcal{L}\,\text{D}\pP\right)\right)_{\alpha mn,n}=\gbm_{\alpha mn,n}=\left(\bom\!\cdot\!\epsilon\right)_{\alpha mn,n}=\left(P_{\alpha j,in}\epsilon_{ij\beta}\epsilon_{\beta mn}\right)_{\alpha m}=\\
	& =\left(-\,P_{\alpha j,in}\epsilon_{ij\beta}\epsilon_{\beta nm}\right)_{\alpha m}=-\,\Curl\Curl\pP.
	\end{align*}
	since in 3 dimensions the Levi-Civita tensors is symmetric for a cyclical permutation of the indexes while antisymmetric for an anticyclic permutation.
	\subsection{Operatorial structure of the anisotropic micromorphic model with consistent coupling boundary conditions}
	
	In Romano et al. \cite{romano2016micromorphic} the authors present the classical linear micromorphic model in the operatorial form which is very useful when we deal with the full anisotropic micromorphic model. In this Section, we want to present the anisotropic micromorphic model with consistent coupling boundary condition in this operatorial form.
	Setting $\bH_{1}\defi H^{1}(\Omega,\bR^{3})\times H^{1}(\Omega,\bR^{3\times3})$,
	let us introduce the formal operator
	\begin{equation}
	\begin{array}{rcl}
	\sA_{1}:\bH_{1} & -\!\!-\!\!\!\longrightarrow & \bH_{2}\\
	\\
	\begin{pmatrix}\pu\\
	\pP
	\end{pmatrix} & \longmapsto & \begin{pmatrix}\sym\text{D} & -\sym\\
	\skew\text{D} & -\skew\\
	0 & \sym\\
	0 & \text{D}
	\end{pmatrix}\!\cdot\!\begin{pmatrix}\pu\\
	\pP
	\end{pmatrix}=\begin{pmatrix}\sym\!\left(\text{D}\pu-\pP\right)\\
	\skew\!\left(\text{D}\pu-\pP\right)\\
	\sym\pP\\
	\text{D}\pP
	\end{pmatrix}
	\end{array}
	\end{equation}
	where
	\begin{equation}
	\bH_{2}\defi L^{2}(\Omega,\Sym(3))\times L^{2}(\Omega,\so)\times H^{1}(\Omega,\Sym(3))\times L^{2}(\Omega,\bR^{3\times3\times3}).
	\end{equation}
	We introduce also 
	\begin{equation*}
	\begin{array}{rcl}
	\sA_{2}:\bH_{2} & -\!\!-\!\!\!\longrightarrow & \bH_{2}\\
	\\
	\begin{pmatrix}\sym\!\left(\text{D}\pu-\pP\right)\\
	\skew\!\left(\text{D}\pu-\pP\right)\\
	\sym\pP\\
	\text{D}\pP
	\end{pmatrix} & \longmapsto & {\displaystyle \frac{1}{2}}\begin{pmatrix}\bC_{\textrm{e}} & 0 & 0 & 0\\
	0 & \bC_{\textrm{c}} & 0 & 0\\
	0 & 0 & \bC_{\textrm{micro}} & 0\\
	0 & 0 & 0 & \mu L_{\textrm{c}}^{2}\widehat{\bL}
	\end{pmatrix}\!\cdot\!\begin{pmatrix}\sym\!\left(\text{D}\pu-\pP\right)\\
	\skew\!\left(\text{D}\pu-\pP\right)\\
	\sym\pP\\
	\text{D}\pP
	\end{pmatrix}={\displaystyle \frac{1}{2}}\begin{pmatrix}\bC_{\textrm{e}}\,\sym\!\left(\text{D}\pu-\pP\right)\\
	\bC_{\textrm{c}}\,\skew\!\left(\text{D}\pu-\pP\right)\\
	\bC_{\textrm{micro}}\,\sym\pP\\
	\mu L_{\textrm{c}}^{2}\widehat{\bL}\,\text{D}\pP
	\end{pmatrix}.
	\end{array}
	\end{equation*}
	The anisotropic potential energy density for the micromorphic model without mixed terms  can be rewritten as follows 
	\[
	W\left(\text{D}\boldsymbol{u},\boldsymbol{P},\text{D}\boldsymbol{P}\right)=\left\langle \frac{1}{2}\begin{pmatrix}\bC_{\textrm{e}} & 0 & 0 & 0\\
	0 & \bC_{\textrm{c}} & 0 & 0\\
	0 & 0 & \bC_{\textrm{micro}} & 0\\
	0 & 0 & 0 & \mu L_{\textrm{c}}^{2}\widehat{\bL}
	\end{pmatrix}\!\cdot\!\begin{pmatrix}\sym\!\left(\text{D}\pu-\pP\right)\\
	\skew\!\left(\text{D}\pu-\pP\right)\\
	\sym\pP\\
	\text{D}\pP
	\end{pmatrix},\begin{pmatrix}\sym\!\left(\text{D}\pu-\pP\right)\\
	\skew\!\left(\text{D}\pu-\pP\right)\\
	\sym\pP\\
	\text{D}\pP
	\end{pmatrix}\right\rangle,
	\]
	where
	\[
	\begin{cases}
	\bC_{\textrm{e}},\bC_{\textrm{micro}}:\Sym(3)\fr\Sym(3) & \textrm{classical \ensuremath{4^{th}}order elasticity tensors, }\\
	\bC_{\textrm{c}}:\so\fr\so & \textrm{ \ensuremath{4^{th}} order coupling tensor},\\
	\widehat{\mathbb{L}}:\bR^{3\times3\times3}\fr\bR^{3\times3\times3} & \textrm{ \ensuremath{6^{th}} order curvature tensor},
	\end{cases}
	\]
	Since 
	\begin{equation}
	L^{2}\!\left(\Omega,\Sym(3)\right)\times L^{2}\!\left(\Omega,\so\right)\simeq L^{2}\!\left(\Omega,\Sym(3)\oplus\so\right)=L^{2}\!\left(\Omega,\bR^{3\times3}\right),
	\end{equation}
	we can identify 
	\begin{equation}
	\begin{pmatrix}
	\bC_{\textrm{e}}\,\sym\!\left(\text{D}\pu-\pP\right)
	\\
	\bC_{\textrm{c}}\,\skew\!\left(\text{D}\pu-\pP\right)
	\\
	\bC_{\textrm{micro}}\,\sym\pP
	\\
	\mu L_{\textrm{c}}^{2}\,\widehat{\bL}\,\text{D}\pP
	\end{pmatrix}
	\longmapsto
	\begin{pmatrix}
	\bC_{\textrm{e}}\,\sym\!\left(\text{D}\pu-\pP\right)+\bC_{\textrm{c}}\,\skew\!\left(\text{D}\pu-\pP\right)
	\\
	\bC_{\textrm{micro}}\,\sym\pP
	\\
	\mu\,L_{\textrm{c}}^{2}\,\widehat{\bL}\,\text{D}\pP
	\end{pmatrix}
	=
	\begin{pmatrix}
	\boldsymbol{\sigma}
	\\
	\bC_{\textrm{micro}}\,\sym\pP
	\\
	\gbm
	\end{pmatrix}
	\end{equation}
	and improving the regularity, i.e., asking for 
	\begin{align}
	\begin{pmatrix}\boldsymbol{\sigma}\\
	\bC_{\textrm{micro}}\,\sym\pP\\
	\gbm
	\end{pmatrix}&\in H\!\left(\Div;\Omega,\bR^{3\times3}\right)\times H^{1}\!\left(\Omega,\Sym\left(3\right)\right)\times H\!\left(\DIV;\Omega,\bR^{3\times3\times3}\right)=:\bH_{3}
	\end{align}
	we can write also the ensuing PDE system in strong form, via the equilibrium operator
	\[
	\begin{array}{rcl}
	\sA_{3}:\bH_{3} & -\!\!-\!\!\!\longrightarrow & L^2(\Omega,\bR^3)\times L^2(\Omega,\bR^{3\times 3})\\
	\\
	\begin{pmatrix}\bosigma\\
	\bC_{\textrm{micro}}\,\sym\pP\\
	\gbm
	\end{pmatrix} & \longmapsto & \begin{pmatrix}\Div & 0 & 0\\
	\id & -\id & \DIV
	\end{pmatrix}\!\cdot\!\begin{pmatrix}\bosigma\\
	\bC_{\textrm{micro}}\,\sym\pP\\
	\gbm
	\end{pmatrix}=\begin{pmatrix}\Div\bosigma\\
	\\
	\bosigma-\bC_{\textrm{micro}}\,\sym\pP+\DIV\gbm
	\end{pmatrix}.
	\end{array}
	\]
	Assuming the partition $\Gamma\cup(\partial\Omega\setminus\overline{\Gamma})\subseteq\partial\Omega$,
	the boundary condition operator for regular fields is
	\begin{equation*}
	\begin{array}{rcl}
	\sB:\bH_{3} & -\!\!-\!\!\!\longrightarrow & H^{-\frac{1}{2\mathstrut}}\left(\partial\Omega\setminus\overline{\Gamma},\bR^{3}\right)\times H^{-\frac{1}{2\mathstrut}}\left(\partial\Omega\setminus\overline{\Gamma},\bR^{3\times 3}\right)\times H^{-\frac{1}{2\mathstrut}}\left(\Gamma,\bR^{3\times3}\right)\\
	\\
	\begin{pmatrix}\bosigma\\
	\bC_{\textrm{micro}}\,\sym\pP\\
	\gbm
	\end{pmatrix} & \longmapsto & 
	\begin{pmatrix}\left.\bosigma\!\cdot\!\pnu\,\right|_{\partial\Omega\setminus\overline{\Gamma}}\\
	\left.\gbm\!\cdot\!\pnu\,\right|_{\partial\Omega\setminus\overline{\Gamma}}\\
	\left.\left(\gbm\!\cdot\!\pnu\right)\!\cdot\!\left(\pnu\otimes\!\pnu\,\right)\right|_{\Gamma}
	\end{pmatrix}.
	\end{array}
	\end{equation*}

\subsubsection{Generalized consistent coupling boundary conditions}\label{GC}
	Recently, generalizations of theorem 1 have been proved in \cite{lewintan2020korn}. Let us consider the larger Sobolev spaces 
	\begin{align*}
       H(\sym\Curl;\Omega,\bR^{3\times3}) & =\left\{ \pP\in L^{2}(\Omega,\bR^{3\times3})\;\left.\right|\;\sym\Curl\pP\in L^{2}(\Omega,\bR^{3\times3})\right\} ,\\
       H(\dev\sym\Curl;\Omega,\bR^{3\times3}) & =\left\{ \pP\in L^{2}(\Omega,\bR^{3\times3})\;\left.\right|\;\dev\sym\Curl\pP\in L^{2}(\Omega,\bR^{3\times3})\right\} ,
    \end{align*}
    equipped respectively with the norms
    \begin{align*}
       \left\Vert \pP\right\Vert _{H\left(\sym\Curl\!\right)}^{2} & =\left\Vert \pP\right\Vert _{L^{2}}^{2}+\left\Vert \sym\Curl\pP\right\Vert _{L^{2}}^{2},\\
       \left\Vert \pP\right\Vert _{H\left(\dev\sym\Curl\!\right)}^{2} & =\left\Vert \pP\right\Vert _{L^{2}}^{2}+\left\Vert \dev\sym\Curl\pP\right\Vert _{L^{2}}^{2}.
    \end{align*}
    In theorem 3.3 \cite{lewintan2020korn} it is proved that
    \begin{align*}
       H_{0,\Gamma}(\sym\Curl;\Omega,\bR^{3\times3}) & \defi\left\{ \pP\in H(\sym\Curl;\Omega,\bR^{3\times3})\;\left.\right|\;\sym\!(\!\left.\pP\times\pnu\,\right|_{\Gamma})=0\right\} ,\\
       H_{0,\Gamma}(\dev\sym\Curl;\Omega,\bR^{3\times3}) & \defi\left\{ \pP\in H(\dev\sym\Curl;\Omega,\bR^{3\times3})\;\left.\right|\;\dev\sym\!(\!\left.\pP\times\pnu\,\right|_{\Gamma})=0\right\} ,
    \end{align*}
    are closed subspaces of $H(\sym\Curl;\Omega,\bR^{3\times3})$
    and $H(\dev\sym\Curl;\Omega,\bR^{3\times3})$ respectively
    and there exist $c_{1},c_{2}>0$ such that the following inequalities hold
    \begin{align*}
       \left\Vert \pP\right\Vert _{L^{2}}^{2} & \leqslant c_{1}\left(\left\Vert \sym\pP\right\Vert _{L^{2}}^{2}+\left\Vert \sym\Curl\pP\right\Vert _{L^{2}}^{2}\right) & \forall\pP\in H_{0,\Gamma}(\sym\Curl;\Omega,\bR^{3\times3}),\\
       \left\Vert \pP\right\Vert _{L^{2}}^{2} & \leqslant c_{2}\left(\left\Vert \sym\pP\right\Vert _{L^{2}}^{2}+\left\Vert \dev\sym\Curl\pP\right\Vert _{L^{2}}^{2}\right) & \forall\pP\in H_{0,\Gamma}(\dev\sym\Curl;\Omega,\bR^{3\times3}).
   \end{align*}
   Moreover, the authors proved that the accounted boundary conditions allow to control $\skew\pP$ with $\sym\Curl\pP$ and $\dev\sym\Curl\pP$ respectively in the introduced closed subspaces. The proposed boundary conditions were derived as usual via integration by parts. We can use these results to further weaken the consistent coupling boundary condition. Indeed, it is straight-forward to prove that the two following boundary-values problems are well posed: 
   
   \medskip
   
   Problem 1: find $(\pu,\pP)\in H^{1}(\Omega,\bR^{3})\times H^{1}(\Omega,\bR^{3\times3})$ such that
   \[
   \left.
         \begin{aligned}
             \Div\bosigma & =0 & \textrm{in} & \;\Omega
             \\
             \quad\bosigma-2\mu_{\text{micro}}\,\text{sym}\,\pP-\lambda_{\text{micro}}\text{tr}(\pP)\id+\,\DIV\gbm & =0 & \textrm{in} & \;\Omega\\
   \\
   \left.\pu\right|_{\Gamma} & =\widehat{\pu} & \textrm{on} & \;\Gamma\\
   \sym\!(\!\left.\pP\times\pnu\right|_{\Gamma}) & =\sym\!(\!\left.\text{D}\pu\times\pnu\right|_{\Gamma}) & \textrm{on} & \;\Gamma\\
   \\
   {\color{red}\left.\bosigma\!\cdot\!\pnu\,\right|_{\partial\Omega\setminus\overline{\Gamma}}}\, & {\color{red}=0} & \textrm{on} & \;\partial\Omega\setminus\overline{\Gamma}\\
   {\color{red}\left.\gbm\!\cdot\!\pnu\,\right|_{\partial\Omega\setminus\overline{\Gamma}}}\, & {\color{red}=0} & \textrm{on} & \;\partial\Omega\setminus\overline{\Gamma}\\
   {\color{red}\left.\left(\gbm\!\cdot\!\pnu\right)\!\cdot\!\left(\pnu\otimes\pnu\right)\,\right|_{{\color{black}{\boldsymbol{\Gamma}}}}}\, & {\color{red}=0} & \textrm{on} & \;\Gamma\\
   {\color{red}\skew\!\left[\!\left.\left(\gbm\!\cdot\!\pnu\right)\!\cdot\!\left(\id-\pnu\otimes\pnu\right)\,\right|_{{\color{black}{\boldsymbol{\Gamma}}}}\right]}\, & {\color{red}=0} & \textrm{on} & \;\Gamma
   \end{aligned}
   \right\} 
   \] 
   where $\bosigma$ and $\gbm$ are as in \eqref{eq:bos} and \eqref{eq:m_third_MM},
   \medskip
   
   Problem 2: find $(\pu,\pP)\in H^{1}(\Omega,\bR^{3})\times H^{1}(\Omega,\bR^{3\times3})$ such that
   \[
   \left.
         \begin{aligned}
             \Div\bosigma & =0 & \textrm{in} & \;\Omega
             \\
             \quad\bosigma-2\mu_{\text{micro}}\,\text{sym}\,\pP-\lambda_{\text{micro}}\text{tr}(\pP)\id+\,\DIV\gbm & =0 & \textrm{in} & \;\Omega\\
   \\
   \left.\pu\right|_{\Gamma} & =\widehat{\pu} & \textrm{on} & \;\Gamma\\
   \dev\sym\!(\!\left.\pP\times\pnu\right|_{\Gamma}) & =\dev\sym\!(\!\left.\text{D}\pu\times\pnu\right|_{\Gamma}) & \textrm{on} & \;\Gamma\\
   \\
   {\color{red}\left.\bosigma\!\cdot\!\pnu\,\right|_{\partial\Omega\setminus\overline{\Gamma}}}\, & {\color{red}=0} & \textrm{on} & \;\partial\Omega\setminus\overline{\Gamma}\\
   {\color{red}\left.\gbm\!\cdot\!\pnu\,\right|_{\partial\Omega\setminus\overline{\Gamma}}}\, & {\color{red}=0} & \textrm{on} & \;\partial\Omega\setminus\overline{\Gamma}\\
   {\color{red}\left.\left(\gbm\!\cdot\!\pnu\right)\!\cdot\!\left(\pnu\otimes\pnu\right)\,\right|_{{\color{black}{\boldsymbol{\Gamma}}}}}\, & {\color{red}=0} & \textrm{on} & \;\Gamma\\
   {\color{red}\skew\!\left[\!\left.\left(\gbm\!\cdot\!\pnu\right)\!\cdot\!\left(\id-\pnu\otimes\pnu\right)\,\right|_{{\color{black}{\boldsymbol{\Gamma}}}}\right]}\, & {\color{red}=0} & \textrm{on} & \;\Gamma\\
   {\color{red}\tr\!\left[\!\left.\left(\gbm\!\cdot\!\pnu\right)\!\cdot\!\left(\id-\pnu\otimes\pnu\right)\,\right|_{{\color{black}{\boldsymbol{\Gamma}}}}\right]}\, & {\color{red}=0} & \textrm{on} & \;\Gamma
        \end{aligned}
   \right\} 
   \]
   Extactly as in \eqref{eq:Korn gen}, we can prove that the following spaces 
   \begin{align*}
   \sH_{\sym}^{\sharp}\!(\Omega) & \defi H^{1}(\Omega,\bR^{3\times3})\cap H_{0,\Gamma}(\sym\Curl;\Omega,\bR^{3\times3}),\\
   \sH_{\dev\sym}^{\sharp}\!(\Omega) & \defi H^{1}(\Omega,\bR^{3\times3})\cap H_{0,\Gamma}(\dev\sym\Curl;\Omega,\bR^{3\times3}),
   \end{align*}
   equipped with the norm $\left\Vert \pP\right\Vert _{\sharp}^{2}=\left\Vert \sym\pP\right\Vert _{L^2}^{2}+\left\Vert \D\!\pP\right\Vert _{L^2}^{2}$
   are Hilbert spaces (i.e., we can not have generalized rigid body motion $\boA\in\so$) and that $\left\Vert \cdot\right\Vert _{\sharp}^{2}$ is
   equivalent to $\left\Vert \cdot\right\Vert _{1,2,\Omega}^{2}$. The well-posedness of the associated homogeneous problems
   \begin{equation}
		\int_{\Omega}W(\text{D}\boldsymbol{u},\pP,\text{D}\pP)\,\text{dV}\longrightarrow\min,\qquad\left(\pu,\pP\right)\in H_{0,\Gamma}^{1}(\Omega,\bR^{3})\times\sH_{\sym}^{\sharp}\!(\Omega)
   \end{equation}
   and
   \begin{equation}
		\int_{\Omega}W(\text{D}\boldsymbol{u},\pP,\text{D}\pP)\,\text{dV}\longrightarrow\min,\qquad\left(\pu,\pP\right)\in H_{0,\Gamma}^{1}(\Omega,\bR^{3})\times\sH_{\dev\sym}^{\sharp}\!(\Omega)
   \end{equation}
   respectively, follows as before.
 }
\end{document}